\documentclass[english,10pt,a4paper]{article}
\usepackage[utf8]{inputenc}
\usepackage[T1]{fontenc}
\usepackage[english]{babel}
\usepackage{hyperref}
\usepackage{amsmath}
\usepackage{amssymb}
\usepackage{amsthm}
\usepackage{textcomp}
\theoremstyle{definition}
\newtheorem{Theorem}{Theorem}[section]
\newtheorem{Definition}{Definition}[section]
\newtheorem{Lemma}{Lemma}[section]
\newtheorem*{Proof}{Proof}
\newtheorem{Proposition}{Proposition}[section]
\newtheorem{Corollary}{Corollary}[section]

\newtheorem{Remark}{Remark}
\newtheorem{Example}{Example}
\usepackage{tikz-cd}
\usepackage{authblk}
\usepackage[numbers,sort&compress]{natbib}
\usepackage{shuffle}
\usepackage{appendix}
\usepackage{tocbibind}
\usepackage{tikz}
\usetikzlibrary{decorations.pathreplacing}

\begin{document}
	\title{The Graded Dual of a Combinatorial Hopf Algebra on Partition Diagrams}
	\author{Lingxiao Hao \thanks{Corresponding author: 21110180006@m.fudan.edu.cn}}	
	\author{Shenglin Zhu}
	\affil{School of Mathematical Sciences, Fudan University, Shanghai, China} 
	\date{\today}
	\maketitle
	\begin{abstract}
		John M. Campbell constructed a Combinatorial Hopf Algebra (CHA) \text{ParSym} on partition diagrams by lifting the CHA structure of \text{NSym}, which is the Hopf algebra of noncommutative symmetric functions. In this paper, we construct the graded dual \text{ParQSym} of \text{ParSym}. Similar to the construction of CHA \text{QSym} of quasisymmetric functions, the CHA structure of  \text{ParQSym} is described in an explicit, combinatorial way. In addition, we present some subcoalgebras and Hopf subalgebras of \text{ParQSym}, gradings and filtrations of \text{ParSym} and \text{ParQSym}, and certain bases of \text{ParSym} and \text{ParQSym} by analogy with the distinguished bases of \text{NSym} and \text{QSym}.

\text{Keywords: } Partition diagram; Combinatorial Hopf algebra; noncommutative symmetric function;  quasisymmetric function
	\end{abstract}
	\section{Introduction}
	 A
	{\it composition} $\alpha = (\alpha_{1}, \alpha_{2}, \ldots , \alpha_{p})$ of a positive integer $n$ is a finite sequence $\alpha_{1}, \alpha_{2}, \ldots , \alpha_{p}$ of positive
	integers such that $|\alpha| := \sum_{i=1}^{p} \alpha_{i} = n$, where $l(\alpha) := p$ is called the {\it length} of $\alpha$. When $n=0$, the empty composition is denoted by $\left( \right) $.

        Let $\mathbb{K}$ be a field, and $\mathbb{K} [[X]]$ be the ring of formal power series in countably many commuting variables $X =\left\lbrace  x_{1}, x_{2}, x_{3},\cdots \right\rbrace $.  {\it Quasisymmetric functions}, introduced by Gessel~\cite{ref36}, are bounded-degree formal power series in $\mathbb{K}$ such that the coefficient of $x_{1}^{\alpha_{1}}x_{2}^{\alpha_{2}}\cdots x_{p}^{\alpha_{p}}$
	is equal to the coefficient of
	$x_{i_{1}}^{\alpha_{1}}x_{i_{2}}^{\alpha_{2}}\cdots x_{i_{p}}^{\alpha_{p}}$, for any non-negative integers $\alpha_{1}, \alpha_{2}, \ldots , \alpha_{p}$ and any strictly increasing sequence of distinct indices $i_{1} < i_{2} <\cdots< i_{p}$. Accordingly, the ring \text{QSym}  of all quasisymmetric functions has the {\it monomial quasisymmetric functions}  $$M_{\alpha}=\mathop{\sum}\limits_{i_{1} < i_{2} <\cdots< i_{p}}x_{i_{1}}^{\alpha_{1}}x_{i_{2}}^{\alpha_{2}}\cdots x_{i_{p}}^{\alpha_{p}},$$for all compositions $\alpha=(\alpha_{1}, \alpha_{2}, \ldots , \alpha_{p})$, as $\mathbb{K}$-basis elements.

	Aguiar, Bergeron, and Sottile~\cite{ref4} introduced the category of  {\it Combinatorial Hopf Algebra}, whose terminal object is \text{QSym}, equipped with a canonical character $\zeta_{\text{QSym}}:\text{QSym}\rightarrow\mathbb{K}$. The graded dual of \text{QSym} is the Hopf algebra \text{NSym} of non-commutative symmetric functions~\cite{ref35,ref10}. The underlying algebra of the bialgebra NSym is  $$\text{NSym}=\mathbb{K}\left\langle H_{1},H_{2},\cdots\right\rangle .$$ For any composition 
	$\alpha = (\alpha_{1}, \alpha_{2}, \ldots , \alpha_{p})$, $$H_{\alpha}:=H_{\alpha_{1}}H_{\alpha_{2}}\cdots H_{\alpha_{p}}.$$ Then $\left\lbrace H_{\alpha} \mid \alpha \text{ is a composition}\right\rbrace $ is a basis of \text{NSym}.

 As discussed by Campbell~\cite{ref1}, the construction of a CHA on partition diagrams is motivated by the wealth of literature in physics and algebraic combinatorics related to Hopf algebras with bases indexed by
 graphs. Recently there have been many discussions on this subject, as one can see in~\cite{uniform block permutations, planar trees, diagrams, Ordered forests, families of trees, planar binary trees, combinatorial structures, subgraph, Rota-Baxter,solomon2024}.
 
 John M. Campbell~\cite{ref1} has constructed a Combinatorial Hopf Algebra \text{ParSym} on partition diagrams, which is a lifting of \text{NSym}. The Hopf subalgebra structures of \text{ParSym} spanned by planar diagrams, by matchings forms, by perfect matchings forms, and by partial permutations are all illustrated also in~\cite{ref1}. 

 In this paper, we define the graded dual \text{ParQSym}  of \text{ParSym}. Similar to that of the CHA \text{QSym} of quasisymmetric functions, and based on partition diagrams, the CHA structure of \text{ParQSym} is described in an explicit, combinatorial way. 
	
        Throughout this article, $\mathbb{K}$ is a field. We first recall the CHA structure of \text{ParSym}~\cite{ref1} in Section~\ref{Knowledge needed}.

Section~\ref{Dual} is the core part of this paper. In this section, we define a CHA \text{ParQSym}, which is the graded dual of \text{ParSym}, and construct some subcoalgebras and Hopf subalgebras of \text{ParQSym}. The definition of \text{ParQSym} relies on partition diagrams~\cite{partition algebra}.

 Section~\ref{$R$-basis of ParSym and $L$-basis of ParQSym} is devoted to constructing the $R$-basis of \text{ParSym} and $L$-basis of \text{ParQSym}, by analogy with the  {\it non-commutative ribbon functions} basis of \text{NSym} and the  {\it fundamental quasi-symmetric functions} basis of \text{QSym}.

 In Section~\ref{Other gradings}, we compare different gradings and filtrations of \text{ParSym} and \text{ParQSym}.

 Inspired by Ricky Ini Liu and Michael Tang~\cite{ref8}, we define a linear map $\eta:\text{QSym}\rightarrow\mathbb{K}$ in section~\ref{Infinitorial Hopf algebra}, to make \text{ParQSym} an infinitorial Hopf algebra.

A deconcatenation basis of \text{ParQSym}, and a method of constructing deconcatenation bases can be found in Section~\ref{Deconcatenation basis}. Similar to the  {\it enriched monomial functions} basis of \text{QSym}, a basis $\left\lbrace \eta_{\pi}\right\rbrace$ of \text{ParQSym} is defined in Section~\ref{The enriched monomial basis}.

 In Section~\ref{The enriched q-monomial basis}, for any $q$ such that $q+1$ is invertible in $\mathbb{K}$, we define a basis $\left\lbrace \eta_{\pi}^{(q)}\right\rbrace$ of \text{ParQSym} by analogy with the  {\it enriched q-monomial functions} basis of \text{QSym}, and give its graded dual. We will also give the formulas for the transformation between the new bases and the bases defined in other sections.

	\section{ParSym}\label{Knowledge needed}
	\numberwithin{equation}{section}
	\setcounter{equation}{0}
	In this section, we recall the CHA structure of \text{ParSym} given by John M. Campbell~\cite{ref1}.

	For a set $S$, a  {\it set partition} of $S$ is a collection of non-empty, pairwise disjoint subsets of $S$ whose union is $S$. An element in a set partition is called a  {\it block}. For instance, $\left\lbrace \left\lbrace 1 \right\rbrace ,\left\lbrace 2,3 \right\rbrace  \right\rbrace $ is a set partition of $\left\lbrace 1,2 ,3 \right\rbrace $. Letting $\pi$ denote a set partition of $\left\lbrace 1,2,\ldots, k,1',2',\ldots, k'\right\rbrace $, this set partition may be denoted with a graph $G$ formed by placing the elements in
	$\left\lbrace 1,2,\ldots, k\right\rbrace$ and $\left\lbrace 1',2',\ldots, k'\right\rbrace $, respectively, into top and bottom rows so
	that the connected components of $G$ agree with the elements in $\pi$. Two
	graphs $G$ and $G'$ on $\left\lbrace 1,2,\ldots, k,1',2',\ldots, k'\right\rbrace $ are considered to be equivalent if the components of these graphs are the same. Also, the set partition
	$\pi$ may be identified with any graph $G$ whose  corresponding set partition is  $\pi$. A {\it partition diagram}, introduced by Halverson \cite{partition algebra}, is the equivalence class of a graph $G$ of the specified form and may be identified with and denoted by $G$ or its corresponding set partition $\pi$.
	The {\it order} of $\pi$ and $G$ is $k$.
	
	\begin{Example}\label{example 1}
		The set partition $\left\lbrace \left\lbrace 1 \right\rbrace, \left\lbrace 2 \right\rbrace, \left\lbrace 4 \right\rbrace, \left\lbrace 3, 1',2'\right\rbrace ,\left\lbrace 3', 4'\right\rbrace \right\rbrace $ may be illustrated as below.
		
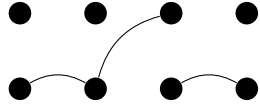
\begin{center}
			\begin{tikzpicture}
				[
				dot/.style={circle, fill=black, inner sep=3pt}, 
				bend angle=30
				]
				\node[dot] (A) at (0,0) {}; 
				\node[dot] (B) at (1,0) {}; 
				\node[dot] (C) at (2,0) {};
				\node[dot] (D) at (3,0) {};
				\node[dot] (A') at (0,1) {}; 
				\node[dot] (B') at (1,1) {}; 
				\node[dot] (C') at (2,1) {}; 
				\node[dot] (D') at (3,1) {};
				\draw (A) to[bend left] (B);
				\draw (B) to[bend left] (C');
				\draw (C) to[bend left] (D);
			\end{tikzpicture}
		\end{center}
	\end{Example}

	Let $A_{i}$ be the set of all partition diagrams of order $i$.
	\text{ParSym} is defined as
	$$\text{ParSym}:=\mathop{\oplus}\limits_{i\geqslant 0}\text{ParSym}_{i},$$where $$\text{ParSym}_{i} =\text{span}_{\mathbb{K}}\left\lbrace H_{\pi}:\pi\in A_{i}\right\rbrace.$$

	For any partition diagrams $\pi$ and $\rho$, $\pi \otimes \rho$ denotes the partition diagram obtained by placing $\rho$ to the right of $\pi$. By convention, the empty diagram $\emptyset$ is the unit of $\otimes$. A partition diagram $\pi$ is {\it $ \otimes $-irreducible}~\cite{ref1}
	if it cannot be expressed as $\rho \otimes \sigma$ for any non-empty partition diagrams $\rho$ and $\sigma$. The partition diagram in Example~\ref{example 1} is $ \otimes $-irreducible.

	The empty partition diagram $\emptyset$ is also $ \otimes $-irreducible. There is an important lemma related to $ \otimes $-irreducible partition diagrams:
	
	\begin{Lemma}\label{length unique}\cite[Lemma 1]{ref1}
		Every non-empty partition diagram $\pi$ can be uniquely written in the form $$\pi=\pi_{1} \otimes \pi_{2} \otimes \cdots \otimes \pi_{l},$$ where  $l\in \mathbb{N}$, and  $\pi_{1},\ldots,\pi_{l} $ are all  non-empty $  \otimes $-irreducible.
	\end{Lemma}
	We call $l$ in the above lemma the  {\it length} of $\pi$ and denote it by $l(\pi)$. By convention, we set $l(\emptyset) = 0$. Then $\pi$ is $ \otimes $-irreducible if and only if $l(\pi)\leq 1$.
	\begin{Definition}\label{ParSym product}
		The product on \text{ParSym} is given by
		\begin{equation*}
			H_{\rho}H_{\sigma}=H_{\rho \otimes \sigma}.
		\end{equation*}
	\end{Definition}

	In addition to the product $ \otimes $, there is another binary operation $\bullet$ we need to introduce before presenting the coproduct. The partition diagram $\pi \bullet \rho$, for any partition diagrams $\pi$ and $\rho$, is the graph obtained by placing $\rho$ to the right of $\pi$ and joining $\pi$ and $\rho$ with an edge between the bottom-right node of $\pi$ and the bottom-left node of $\rho$ \cite{ref1}. By convention, the empty diagram $\emptyset$ is the unit of $\bullet$, that is, $\emptyset \bullet  \rho = \rho$ and $\pi  \bullet \emptyset = \pi$.
	\begin{Example}
		The partition diagram in Example~\ref{example 1} can be represented as $\pi \bullet \rho$ in two ways:	
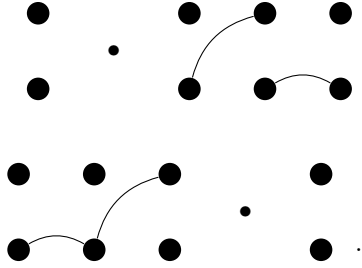
\begin{center}
			\begin{tikzpicture}
				[
				dot/.style={circle, fill=black, inner sep=3pt}, 
				bend angle=30
				]
				\node[dot] (A) at (0,0) {}; 
				\node[dot] (B) at (2,0) {}; 
				\node[dot] (C) at (3,0) {};
				\node[dot] (D) at (4,0) {};
				\node[dot] (A') at (0,1) {}; 
				\node[dot] (B') at (2,1) {}; 
				\node[dot] (C') at (3,1) {}; 
				\node[dot] (D') at (4,1) {};
				\draw (B) to[bend left] (C');
				\draw (C) to[bend left] (D);
				\node at (1,0.5) {$\bullet$};
			\end{tikzpicture}
		\end{center}
		or
		\begin{center}
			\begin{tikzpicture}
				[
				dot/.style={circle, fill=black, inner sep=3pt}, 
				bend angle=30
				]
				\node[dot] (A) at (0,0) {}; 
				\node[dot] (B) at (1,0) {}; 
				\node[dot] (C) at (2,0) {};
				\node[dot] (D) at (4,0) {};
				\node[dot] (A') at (0,1) {}; 
				\node[dot] (B') at (1,1) {}; 
				\node[dot] (C') at (2,1) {}; 
				\node[dot] (D') at (4,1) {};
				\draw (A) to[bend left] (B);
				\draw (B) to[bend left] (C');
				\node at (3,0.5) {$\bullet$};
				\node at (4.5,0) {.};
			\end{tikzpicture}
		\end{center}
	\end{Example}
	A partition diagram is called  {\it $ \bullet $-irreducible} if it cannot be written as $\pi \bullet \rho$ of two non-empty partition diagrams $\pi$ and $\rho$ \cite{ref1}. Notice that $\emptyset$ is also $ \bullet $-irreducible.
	\begin{Example}
		Here is a $ \bullet $-irreducible partition diagram
		
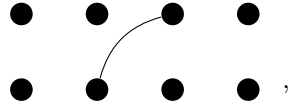
\begin{center}
			\begin{tikzpicture}
				[
				dot/.style={circle, fill=black, inner sep=3pt}, 
				bend angle=30
				]
				\node[dot] (A) at (0,0) {}; 
				\node[dot] (B) at (1,0) {}; 
				\node[dot] (C) at (2,0) {};
				\node[dot] (D) at (3,0) {};
				\node[dot] (A') at (0,1) {}; 
				\node[dot] (B') at (1,1) {}; 
				\node[dot] (C') at (2,1) {}; 
				\node[dot] (D') at (3,1) {};
				\draw (B) to[bend left] (C');
				\node at (3.5,0) {,};
			\end{tikzpicture}
		\end{center}
		which is $ \otimes $-reducible as it can be written as 
		
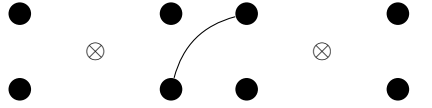
\begin{center}
			\begin{tikzpicture}
				[
				dot/.style={circle, fill=black, inner sep=3pt}, 
				bend angle=30
				]
				\node[dot] (A) at (0,0) {}; 
				\node[dot] (B) at (2,0) {}; 
				\node[dot] (C) at (3,0) {};
				\node[dot] (D) at (5,0) {};
				\node[dot] (A') at (0,1) {}; 
				\node[dot] (B') at (2,1) {}; 
				\node[dot] (C') at (3,1) {}; 
				\node[dot] (D') at (5,1) {};
				\draw (B) to[bend left] (C');
				\node at (1,0.5) {$\otimes$};
				\node at (4,0.5) {$\otimes$};
				\node at (5.5,0) {.};
			\end{tikzpicture}
		\end{center}		
	\end{Example} 
	\begin{Remark}
		For any partition diagrams $\pi$, $\rho$ and $\sigma$, we have that
		\begin{equation*}
			\begin{aligned}
				(\pi\otimes\rho)\otimes\sigma&=\pi\otimes(\rho\otimes\sigma),\\
				(\pi\otimes\rho)\bullet\sigma&=\pi\otimes(\rho\bullet\sigma),\\
				(\pi\bullet\rho)\otimes\sigma&=\pi\bullet(\rho\otimes\sigma),\\
				(\pi\bullet\rho)\bullet\sigma&=\pi\bullet(\rho\bullet\sigma).\\
			\end{aligned}
		\end{equation*}
		
	\end{Remark}

	The coproduct of \text{ParSym} is given as follows: 
	\begin{Definition}\label{ParSym coproduct}\cite{ref1}
		For any $ \otimes $-irreducible partition diagram $\sigma$, 
		\begin{equation}
			\bigtriangleup H_{\sigma}=\mathop{\sum}\limits_{\pi,\rho}H_{\pi} \otimes  H_{\rho},
		\end{equation}where the sum is over all $ \otimes $-irreducible partition diagrams $\pi,\rho $ such that $\pi  \bullet  \rho=\sigma$.
		Let $\bigtriangleup$ be the linear map compatible with the product of
		\text{ParSym} (that is, for any $\pi=\pi_{1} \otimes \pi_{2} \otimes \cdots \otimes \pi_{l}$, let $\bigtriangleup H_{\pi}=(\bigtriangleup H_{\pi_{1}})(\bigtriangleup H_{\pi_{2}})\cdots(\bigtriangleup H_{\pi_{l}})$).
	\end{Definition}
	With both the product and the coproduct,  \text{ParSym} is a bialgebra. The following theorems will be used later.
	
	\begin{Theorem}[\protect{\cite[Theorem 5-6]{ref1}}]\label{irr}
		Let $\rho$ and $\pi$ be two partition diagrams. Then $\rho \bullet \pi$ is $ \otimes $-irreducible
		if and only if $\rho$ and $\pi$ are $ \otimes $-irreducible.
	\end{Theorem}
	\begin{Theorem}\label{irr2}
		Every non-empty partition diagram $\pi$ can be uniquely written as $$\pi=\pi_{1} \bullet \pi_{2} \bullet \cdots \bullet \pi_{l},$$ where  $l\in \mathbb{N}$, and $\pi_{1},\ldots,\pi_{l} $ are all non-empty  $\bullet $-irreducible.
	\end{Theorem}
	
	\section{A dual CHA construction}\label{Dual}

        In this section, we will define the graded dual $\text{ParQSym}$ of \text{ParSym}.
        Let $A_{i}$ be the set of all partition diagrams of order $i$ and for a compositions $\alpha=(\alpha_{1}, \alpha_{2}, \ldots , \alpha_{p})$, $$M_{\alpha}=\mathop{\sum}\limits_{i_{1} < i_{2} <\cdots< i_{p}}x_{i_{1}}^{\alpha_{1}}x_{i_{2}}^{\alpha_{2}}\cdots x_{i_{p}}^{\alpha_{p}}.$$

        Take $$\text{ParQSym}=\mathop{\oplus}\limits_{n\geqslant 0}\text{ParQSym}_{i},$$where $\text{ParQSym}_{i} =\text{span}_{\mathbb{K}}\left\lbrace M_{\pi}:\pi\in A_{i}\right\rbrace $. Define a bilinear form $\left\langle \cdot,\cdot\right\rangle :\text{ParQSym}  \otimes  \text{ParSym} \rightarrow \mathbb{K}$ such that
	$\left\langle M_{\rho},H_{\pi}\right\rangle=\delta_{\rho,\pi}= M_{\rho}(H_{\pi})$, and we can view $M_{\rho} $  as an element in $\text{Hom}_{\mathbb{K}}(\text{ParSym},\mathbb{K})$. For each natural number $i$, if we restrict all $M_{\rho}\in \text{ParQSym}_{i}$ to $\text{ParSym}_{i}$, we will find that $\text{ParQSym}_{i}=\text{ParSym}_{i}^{\ast}$ as $\dim\text{ParSym}_{i}<\infty$. Thus, \text{ParQSym} is the graded dual of \text{ParSym}.

	Let $\left\langle M_{\pi} \otimes  M_{\rho},H_{\pi'} \otimes  H_{\rho'}\right\rangle :=\left\langle M_{\pi},H_{\pi'}\right\rangle\left\langle M_{\rho},H_{\rho'}\right\rangle$, for any partition diagrams $\pi$, $\pi'$, $\rho$ and $\rho'$.
	
	\subsection{The coalgebra structure on partition diagrams}

 	We need to define a coproduct $\bigtriangleup$ on $\text{ParQSym} $ such that \begin{equation*}
		\begin{split}
			\left\langle \bigtriangleup M_{\sigma},H_{\pi} \otimes  H_{\rho}\right\rangle &=M_{\sigma}(H_{\pi} H_{\rho})\\&=M_{\sigma}(H_{\pi \otimes \rho})\\&=\delta_{\sigma,\pi \otimes \rho}.
		\end{split}
	\end{equation*} 
	One can see that the coefficient of $M_{\pi} \otimes  M_{\rho} $ in $\bigtriangleup M_{\sigma}$ is non-zero if and only if $\pi  \otimes  \rho=\sigma$. With this observation,  we can define the coproduct of \text{ParQSym}.  
	\begin{Definition}[Coproduct on \text{ParQSym}]\label{coproduct}
		For any partition diagram $\sigma$, the coproduct is defined as $$\bigtriangleup M_{\sigma}=\mathop{\sum}\limits_{\pi  \otimes  \rho=\sigma}M_{\pi} \otimes  M_{\rho},$$where the sum ranges over all factorizations of $\sigma$ into $\pi  \otimes  \rho$.
	\end{Definition}

 	It is easy to see $\bigtriangleup$ is coassociative with the associativity of $ \otimes $, since both $(id \otimes  \bigtriangleup )(\bigtriangleup M_{G})$ and $(\bigtriangleup \otimes  id )(\bigtriangleup M_{G})$ can be rewritten as $$\mathop{\sum}\limits_{\pi  \otimes  \rho \otimes \sigma=G}M_{\pi} \otimes  M_{\rho} \otimes  M_{\sigma}. $$

 	By defining a counit morphism $\varepsilon: \text{ParQSym}\rightarrow \mathbb{K}$ so that $\varepsilon(M_{\emptyset}) = 1_{\mathbb{K}}$ and $\varepsilon(M_{\rho})
	= 0$ for any non-empty partition diagram $\rho$, we 
	obtain a coalgebra structure on \text{ParQSym}. The unique group-like element is $M_{\emptyset}$.

	Recall that the dual basis of the complete homogeneous basis $\left\lbrace H_{\alpha}\right\rbrace$ of the CHA \text{NSym} is the monomial quasisymmetric basis $\left\lbrace M_{\alpha}\right\rbrace$ of the CHA \text{QSym}. The coproduct is given by $\bigtriangleup M_{\alpha}=\mathop{\sum}\limits_{\beta \cdot \gamma=\alpha}M_{\beta} \otimes  M_{\gamma}$, where ``$\cdot$'' denotes the concatenation operation for integer compositions. Our comultiplication operation on $\text{ParQSym} $ may be seen as an analogue of this concatenation operation.

	\subsection{The algebra structure on partition diagrams}
	We want to define a product $ \star $ on $\text{ParQSym} $ such that  for any partition diagrams $G_{1}$, $G_{2}$, $\pi=\pi_{1} \otimes \pi_{2} \otimes \cdots \otimes \pi_{l}$, where $\pi_{i} $ is $  \otimes $-irreducible for all $1\leq i\leq l$, \begin{equation*}
		\begin{split}
			\left\langle M_{G_{1}} \star  M_{G_{2}},H_{\pi}\right\rangle &=\left\langle M_{G_{1}} \otimes  M_{G_{2}},\bigtriangleup H_{\pi}\right\rangle \\&=\left\langle M_{G_{1}} \otimes  M_{G_{2}},(\bigtriangleup H_{\pi_{1}})(\bigtriangleup H_{\pi_{2}})\cdots(\bigtriangleup H_{\pi_{l}})\right\rangle \\&=\left\langle M_{G_{1}} \otimes  M_{G_{2}},(\sum H_{\rho_{1}} \otimes  H_{\sigma_{1}})\cdots(\sum H_{\rho_{l}} \otimes  H_{\sigma_{l}})\right\rangle \\&=\left\langle M_{G_{1}} \otimes  M_{G_{2}},\sum H_{\rho_{1} \otimes \rho_{2} \otimes \cdots \otimes \rho_{l}} \otimes  H_{\sigma_{1} \otimes \sigma_{2} \otimes \cdots \otimes \sigma_{l}}\right\rangle \\&=\sum\delta_{G_{1},\rho_{1} \otimes \rho_{2} \otimes \cdots \otimes \rho_{l}}\delta_{G_{2},\sigma_{1} \otimes \sigma_{2} \otimes \cdots \otimes \sigma_{l}},
		\end{split}
	\end{equation*} where  $ \rho_{i},\sigma_{i}$ are $  \otimes $-irreducible, $\rho_{i} \bullet \sigma_{i}=\pi_{i}$ for all $1\leq i\leq l$.  We can find that $\left\langle M_{G_{1}} \star  M_{G_{2}},H_{\pi}\right\rangle$ is non-zero only when there exists an integer $l$ and a pair of ordered sequences of $  \otimes $-irreducible (not necessarily non-empty) partition diagrams $\left( \rho_{1},\ldots,\rho_{l}\right)$ and $\left( \sigma_{1},\ldots,\sigma_{l}\right)$ such that $G_{1}=\rho_{1} \otimes \rho_{2} \otimes \cdots \otimes \rho_{l}$, $G_{2}=\sigma_{1} \otimes \sigma_{2} \otimes \cdots \otimes \sigma_{l}$, and $\pi=\rho_{1} \bullet \sigma_{1} \otimes \rho_{2} \bullet \sigma_{2} \otimes \cdots \otimes \rho_{l} \bullet \sigma_{l}$. 
	From now on, the superscript ``$\,\hat{}\; $'' on a sequence of partition diagrams is used to 
	remove any empty partition diagrams and keep the order of non-empty ones. For instance, if $\rho$, $\sigma$, and $\pi$ are non-empty partition diagrams, then $$(\emptyset,\rho,\sigma,\emptyset,\pi,\emptyset,\emptyset)\hat{}=(\rho,\sigma,\pi).$$
	Now we can define the product of $\text{ParQSym} $ as follows.

	\begin{Definition}\label{product} (Product on \text{ParQSym})
		For any non-empty partition diagrams $\rho=\rho_{1} \otimes \rho_{2} \otimes  \cdots  \otimes \rho_{n}$, $\sigma=\sigma_{1} \otimes \sigma_{2} \otimes  \cdots  \otimes \sigma_{m}$, where $ \rho_{i}$, $\sigma_{j}$ are non-empty $  \otimes $-irreducible partition diagrams for $1\leq i\leq n$, $1\leq j\leq m$ (with Lemma~\ref{length unique}, every non-empty partition diagram can be uniquely written in this form), $$M_{\rho}  \star  M_{\sigma}=\mathop{\sum}M_{\rho'_{1} \bullet \sigma'_{1} \otimes \cdots \otimes \rho'_{k} \bullet \sigma'_{k}},$$ where the sum is over all $ k\in \mathbb{N}$ and all pairs of ordered sequences $\left( \rho'_{1},\ldots,\rho'_{k}\right)$ and $\left( \sigma'_{1},\ldots,\sigma'_{k}\right)$ such that $\left( \rho'_{1},\ldots,\rho'_{k}\right)\hat{}=\left( \rho_{1},\ldots,\rho_{n}\right) $ and $\left( \sigma'_{1},\ldots,\sigma'_{k}\right)\hat{}=\left( \sigma_{1},\ldots,\sigma_{m}\right) $, and such that $\rho'_{s} \bullet \sigma'_{s}$ is non-empty for all $1\leq s\leq k$.
		By convention, we set $M_{\rho}  \star  M_{\emptyset}=M_{\rho}=M_{\emptyset}  \star  M_{\rho}$ for every partition diagram $\rho$, thus $M_{\emptyset} $ is the unit.
	\end{Definition}
	There is another way to understand this product. For any non-empty partition diagrams $\rho=\rho_{1} \otimes \rho_{2} \otimes  \cdots  \otimes \rho_{n}$, $\sigma=\sigma_{1} \otimes \sigma_{2} \otimes  \cdots  \otimes \sigma_{m}$, where $ \rho_{i}$, $\sigma_{j}$ are non-empty $  \otimes $-irreducible partition diagrams for $1\leq i\leq n$, $1\leq j\leq m$,$$M_{\rho}  \star  M_{\sigma}=\mathop{\sum}M_{\pi_{1}  \divideontimes_{1} \pi_{2}  \divideontimes_{2} \cdots  \divideontimes_{n+m-1} \pi_{n+m}}. $$The sum is over possible expressions $\pi_{1}  \divideontimes_{1} \pi_{2}  \divideontimes_{2} \cdots  \divideontimes_{n+m-1} \pi_{n+m}$ satisfying the following conditions:\\
	(1) $ \pi_{1}\pi_{2}\cdots\pi_{n+m}$ is a word shuffle of  $\rho_{1}\rho_{2} \cdots \rho_{n} $ and $ \sigma_{1}\sigma_{2} \cdots \sigma_{m}$;\\
	(2) $ \divideontimes_{i} \in\left\lbrace  \otimes , \bullet \right\rbrace$ for $1\leq i\leq n$;\\
	(3) only when $\pi_{l}=\rho_{i}$ and $\pi_{l+1}=\sigma_{j}$ for some $1\leq i\leq n$, $1\leq j\leq m$, $ \divideontimes_{l}$ can be $ \bullet $.
	\begin{Proposition}
		The operation $ \star $ is associative.
	\end{Proposition}
	\begin{Proof}
		For any non-empty partition diagrams $\rho=\rho_{1} \otimes \rho_{2} \otimes  \cdots  \otimes \rho_{n}$, $\sigma=\sigma_{1} \otimes \sigma_{2} \otimes  \cdots  \otimes \sigma_{m}$, $\pi=\pi_{1} \otimes \pi_{2} \otimes  \cdots  \otimes \pi_{l}$, where $ \rho_{i}$, $\sigma_{j}$ and $\pi_{t} $ are non-empty $  \otimes $-irreducible partition diagrams for $1\leq i\leq n$, $1\leq j\leq m$, $1\leq t\leq l$,
		from the associativity of $ \bullet $ and \cite[Theorem 5]{ref1}, both $(M_{\rho}  \star  M_{\sigma}) \star  M_{\pi}$ and $M_{\rho}  \star  (M_{\sigma} \star  M_{\pi})$ can be rewritten as $$\mathop{\sum}M_{\rho'_{1} \bullet \sigma'_{1} \bullet \pi'_{1} \otimes \cdots \otimes \rho'_{r} \bullet \sigma'_{r} \bullet \pi'_{r}}, $$
		where the sum is over all $ r\in \mathbb{N}$ and all triples of ordered sequences $\left( \rho'_{1},\ldots,\rho'_{r}\right)$, $\left( \sigma'_{1},\ldots,\sigma'_{r}\right)$ and $\left( \pi'_{1},\ldots,\pi'_{r}\right)$ such that $\left( \rho'_{1},\ldots,\rho'_{r}\right)\hat{}=\left( \rho_{1},\ldots,\rho_{n}\right) $, $\left( \sigma'_{1},\ldots,\sigma'_{r}\right)\hat{}=\left( \sigma_{1},\ldots,\sigma_{m}\right) $, $\left( \pi'_{1},\ldots,\pi'_{r}\right)\hat{}=\left( \pi_{1},\ldots,\pi_{l}\right) $ and such that $\rho'_{s} \bullet \sigma'_{s} \bullet \pi'_{s}$ is non-empty for all $1\leq s\leq r$. $\qedsymbol$
	\end{Proof}
	By defining a unit morphism $\eta:\mathbb{K}\rightarrow \text{ParQSym} $ so that $\varepsilon(1_{\mathbb{K}}) =M_{\emptyset} $,  we 
	obtain an algebra structure on \text{ParQSym}.

        Now we compare the product of the $\left\lbrace M_{\pi}\right\rbrace$-basis of \text{ParQSym} with the product in the $\left\lbrace M_{\alpha}\right\rbrace$-basis in \text{QSym}, which can be defined as follows  \cite[Proposition 5.1.3]{ref2}.
	Fix three disjoint chain posets $(i_{1}<\cdots<i_{l})$, $(j_{1}<\cdots<j_{m})$ and $(k_{1}<k_{2}<\cdots)$, if $\alpha=(\alpha_{1},\alpha_{2},\ldots,\alpha_{l})$, $\beta=(\beta_{1},\ldots,\beta_{m})$ are two compositions then $$M_{\alpha}M_{\beta}=\sum M_{wt(f)},$$ where the sum is over all $p \in \mathbb{N}$ and all  surjective and strictly order-preserving maps $f$ from the disjoint union of two chains to a chain$$(i_{1}<\cdots<i_{l})\sqcup (j_{i}<\cdots<j_{m})\rightarrow (k_{1}<k_{2}<\cdots<k_{p}),$$ and where the composition $wt(f)=(wt_{1}(f),\ldots,wt_{p}(f))$ is given by $$wt_{s}(f)=\mathop{\sum}\limits_{u:f(i_{u})=k_{s}}\alpha_{u}+\mathop{\sum}\limits_{v:f(j_{v})=k_{s}}\beta_{v}. $$

 Our definition of $ \star $ in the basis of $\left\lbrace M_{\pi}  \right\rbrace $ can be seen similarly. Fix three disjoint chain posets $(i_{1}<\cdots<i_{l})$, $(j_{1}<\cdots<j_{m})$ and $(k_{1}<k_{2}<\cdots)$, if $\rho=\rho_{1} \otimes \rho_{2} \otimes  \cdots  \otimes \rho_{l}$, $\sigma=\sigma_{1} \otimes \sigma_{2} \otimes  \cdots  \otimes \sigma_{m}$, where $ \rho_{i}$, $\sigma_{j}$ are non-empty $  \otimes $-irreducible partition diagrams for $1\leq i\leq n$, $1\leq j\leq m$, then $$M_{\rho}M_{\sigma}=\sum M_{wt(f)},$$ where the sum is over all $p \in \mathbb{N}$ and all  surjective and strictly order-preserving maps $f$ from the disjoint union of two chains to a chain$$(i_{1}<\cdots<i_{l})\sqcup (j_{i}<\cdots<j_{m})\rightarrow (k_{1}<k_{2}<\cdots<k_{p}),$$ and where the diagram $wt(f)=wt_{1}(f) \otimes \cdots \otimes  wt_{p}(f)$ is given by
	\begin{equation*}
		wt_{s}(f)=
		\begin{cases}
			\rho_{u} \bullet  \sigma_{v} & \text{if } f(i_{u})=f(j_{v})=k_{s}; \\
			\rho_{u} & \text{if } f(i_{u})=k_{s} \text{ and } \left\lbrace  v:f(j_{v})=k_{s} \right\rbrace =\emptyset; \\
			\sigma_{v} & \text{if } f(j_{v})=k_{s} \text{ and } \left\lbrace  u:f(i_{u})=k_{s} \right\rbrace =\emptyset.
		\end{cases}
	\end{equation*}
	This definition is given by	replacing the integers, ``+'' and ``,''  in compositions by $ \otimes $-irreducible partition diagrams, ``$ \bullet $'' and ``$ \otimes $'' respectively, and then modifying a little bit. A difference that needs to be noticed is that ``$ \bullet $'' is not commutative like ``+''.

	\subsection{ParQSym is a bialgebra}
	With the previously defined product and coproduct, we get a bialgebra $\text{ParQSym}$.
	\begin{Proposition}
		$(\text{ParQSym}, \star ,\bigtriangleup,M_{\emptyset},\varepsilon)$ is a bialgebra.
	\end{Proposition}
	\begin{Proof}
		For any non-empty partition diagrams $\rho=\rho_{1} \otimes \rho_{2} \otimes  \cdots  \otimes \rho_{n}$ and $\sigma=\sigma_{1} \otimes \sigma_{2} \otimes  \cdots  \otimes \sigma_{m}$, where $ \rho_{i}$, $\sigma_{j}$ are non-empty $  \otimes $-irreducible partition diagrams for $1\leq i\leq n$, $1\leq j\leq m$, 	\begin{equation}\label{1}
			\begin{split}
				\bigtriangleup(M_{\rho}  \star  M_{\sigma})=&\bigtriangleup(\mathop{\sum}M_{\rho'_{1} \bullet \sigma'_{1} \otimes \cdots \otimes \rho'_{k} \bullet \sigma'_{k}})\\=&\sum\mathop{\sum}\limits_{0\leq i\leq k}M_{\rho'_{1} \bullet \sigma'_{1} \otimes \cdots \otimes \rho'_{i} \bullet \sigma'_{i}} \otimes  M_{\rho'_{i+1} \bullet \sigma'_{i+1} \otimes \cdots \otimes \rho'_{k} \bullet \sigma'_{k}},\\ 
			\end{split}
		\end{equation}
		where $\left( \rho'_{1},\ldots,\rho'_{k}\right)$ and $\left( \sigma'_{1},\ldots,\sigma'_{k}\right)$ are as in Definition~\ref{product}, and  we set $\rho'_{0}=\sigma'_{0}=\rho'_{n+1}=\sigma'_{m+1}=\emptyset$. On the other hand, \begin{equation*}
			\begin{split}
				&(\bigtriangleup M_{\rho})  \star  \bigtriangleup (M_{\sigma})\\=&(\mathop{\sum}\limits_{0\leq i\leq n}M_{\rho_{1} \otimes \cdots \otimes \rho_{i}} \otimes  M_{\rho_{i+1} \otimes \cdots \otimes \rho_{n}}) \star (\mathop{\sum}\limits_{0\leq j\leq m}M_{\sigma_{1} \otimes \cdots \otimes  \sigma_{j}} \otimes  M_{\sigma_{j+1} \otimes \cdots \otimes \sigma_{m}})\\
				=&\mathop{\sum}\limits_{0\leq i\leq n,0\leq j\leq m}(M_{\rho_{1} \otimes \cdots \otimes \rho_{i}} \star  M_{\sigma_{1} \otimes \cdots \otimes  \sigma_{j}}) \otimes (M_{\rho_{i+1} \otimes \cdots \otimes \rho_{n}} \star  M_{\sigma_{j+1} \otimes \cdots \otimes \sigma_{m}})\\=&\mathop{\sum}\limits_{0\leq i\leq n,0\leq j\leq m}(\sum M_{\rho'_{1} \bullet \sigma'_{1} \otimes \cdots \otimes \rho'_{s} \bullet \sigma'_{s}} \otimes  M_{\rho'_{s+1} \bullet \sigma'_{s+1} \otimes \cdots \otimes \rho'_{s+t} \bullet \sigma'_{s+t}}),
			\end{split}
		\end{equation*}
		where we set $\rho_{0}=\sigma_{0}=\rho_{k+1}=\sigma_{k+1}=\emptyset$, and  for each pair of $(i,j)$, the inner sum is over all $s,t\in \mathbb{N}$ and all quaternions of ordered sequences $\left( \rho'_{1},\ldots,\rho'_{s}\right)$, $\left( \sigma'_{1},\ldots,\sigma'_{s}\right)$, $\left( \rho'_{s+1},\ldots,\rho'_{s+t}\right)$, $\left( \sigma'_{s+1},\ldots,\sigma'_{s+t}\right)$ such that $\left( \rho'_{1},\ldots,\rho'_{s}\right)\hat{}=\left( \rho_{1},\ldots,\rho_{i}\right) $, $\left( \sigma'_{1},\ldots,\sigma'_{s}\right)\hat{}=\left( \sigma_{1},\ldots,\sigma_{j}\right) $, $\left( \rho'_{s+1},\ldots,\rho'_{s+t}\right)\hat{}=\left( \rho_{i+1},\ldots,\rho_{n}\right) $ and $\left( \sigma'_{s+1},\ldots,\sigma'_{s+t}\right)\hat{}=\left( \sigma_{j+1},\ldots,\sigma_{m}\right) $ and such that $\rho'_{l} \bullet \sigma'_{l}$ is non-empty for all $1\leq l\leq s+t$. By concatenating $\left( \rho'_{1},\ldots,\rho'_{s}\right)$ and $\left( \rho'_{s+1},\ldots,\rho'_{s+t}\right)$, $\left( \sigma'_{s+1},\ldots,\sigma'_{s+t}\right)$, $\left( \sigma'_{1},\ldots,\sigma'_{s}\right)$ and $\left( \sigma'_{s+1},\ldots,\sigma'_{s+t}\right)$, we get the same pairs as in (\ref{1}). So we get $$\bigtriangleup(M_{\rho}  \star  M_{\sigma})= (\bigtriangleup M_{\rho})  \star  \bigtriangleup (M_{\sigma}).$$
		In addition, $\varepsilon(a)\varepsilon(b)\neq 0$ if and only if $a=b=M_{\emptyset}$ if and only if $\varepsilon(a \star  b)\neq 0$, thus $$\varepsilon(a)\varepsilon(b)=\varepsilon(a \star  b)$$for any partition diagrams $a,b$. $\qedsymbol$
	\end{Proof}
	With this product and coproduct, the primitive elements of $\text{ParQSym}$ are closely related to non-empty $ \otimes $-irreducible diagrams.
	\begin{Corollary}\label{primitive}
		The set $P(\text{ParQSym})$ of the primitive elements of $\text{ParQSym}$ is  $$P(\text{ParQSym})=\text{span}_{\mathbb{K}}\left\lbrace  M_{\sigma}\mid \sigma \text{ is non-empty}  \otimes -\text{irreducible}\right\rbrace . $$
	\end{Corollary}
	\begin{Proof}
		This follows directly from the definition of the coproduct and the fact that $M_{\emptyset}$ is the unit of the product. $\qedsymbol$
	\end{Proof}

	\subsection{ParQSym is a Hopf algebra}
	We can also assign an antipode to \text{ParQSym} to make it a Hopf algebra. The following lemma is useful to build connected graded Hopf algebras.
	 \begin{Lemma}[\cite{Sweedler}]\label{connected}
		If a bialgebra $H=\oplus_{n\geq 0}H_{n}$ is both a graded algebra and a graded coalgebra of the same grading, and $H_{0}=\mathbb{K}$, then it is a connected graded Hopf algebra.
	\end{Lemma}
	It is clear that $$\bigtriangleup \text{ParQSym}_{n}\subseteq \bigoplus_{i+j=n}\text{ParQSym}_{i} \otimes  \text{ParQSym}_{j},$$ $$\text{ParQSym}_{i} \star  \text{ParQSym}_{j}\subseteq \text{ParQSym}_{i+j},$$ and $\text{ParQSym}_{0}=\mathbb{K}\left\lbrace M_{\emptyset}\right\rbrace \cong\mathbb{K}$, by Lemma~\ref{connected} \text{ParQSym} is a connected Hopf algebra.
	\begin{Corollary}\label{coradical}
		The coradical of the connected Hopf algebra  $\text{ParQSym}$ is $\mathbb{K}\left\lbrace M_{\emptyset} \right\rbrace\cong \mathbb{K}$.
	\end{Corollary}
	There is a general formula for the antipode of a connected graded Hopf algebra $H$, due to Takeuchi~\cite[Lemma 14]{ref5}.  Let $(H, m, u, \bigtriangleup, \varepsilon)$ be an connected graded bialgebra, set 
	\begin{equation*}
		\begin{aligned}
			m^{(-1)}&=u,
			\bigtriangleup^{(-1)}=\varepsilon,\\
			m^{(0)}&=\bigtriangleup^{(0)}=id,\\
			m^{(1)}&=m,
			\bigtriangleup^{(1)}=\bigtriangleup,
		\end{aligned}
	\end{equation*}
	and for any $k\geq 2$, 
	\begin{equation*}
		\begin{aligned}
			m^{(k)}&=m(m^{(k-1)} \otimes  id),\\
			\bigtriangleup^{(k)}&=(\bigtriangleup^{(k-1)} \otimes  id)\bigtriangleup.
		\end{aligned}
	\end{equation*}
	Then the antipode is given by 
	\begin{equation}\label{Antipode}
		S=\mathop{\sum}\limits_{k\geq 0}(-1)^{k}m^{(k-1)}(id-u\varepsilon)^{ \otimes  k}\bigtriangleup^{(k-1)}.
	\end{equation}
	The formula	(\ref{Antipode}) is often referred to as    {\it Takeuchi's formula}.
	By this formula, we can get an antipode as follows.  First, set $S(M_{\emptyset})=M_{\emptyset}$. Then, for any non-empty partition diagram $\pi=\pi_{1} \otimes \pi_{2} \otimes \cdots \otimes \pi_{n}$, where $\pi_{i}$ is non-empty $ \otimes $-irreducible partition diagram for $1\leq i\leq n$, let
	$$S (M_{\pi})=\mathop{\sum}\limits_{k\geq0}(-1)^{k}(\mathop{\sum}\limits_{r\geq 0}M_{\pi_{1}^{1} \bullet \pi_{2}^{1} \bullet \cdots \bullet \pi_{k}^{1} \otimes \pi_{1}^{2} \bullet \pi_{2}^{2} \bullet \cdots \bullet \pi_{k}^{2} \otimes \cdots \otimes \pi_{1}^{r} \bullet \pi_{2}^{r} \bullet \cdots \bullet \pi_{k}^{r}}),$$ where $(\pi_{1}^{1},\pi_{1}^{2},\ldots,\pi_{1}^{r},\pi_{2}^{1},\ldots,\pi_{2}^{r},\ldots,\pi_{k}^{1},\ldots,\pi_{k}^{r})\hat{}=(\pi_{1},\pi_{2},\ldots,\pi_{n})$ and $\pi_{1}^{s} \bullet \pi_{2}^{s} \bullet \cdots \bullet \pi_{k}^{s}$ and $\pi_{i}^{1} \otimes \pi_{i}^{2} \otimes \cdots \otimes \pi_{i}^{r}$ are non-empty for $1\leq s\leq r$, $1\leq i\leq k $.

        By a well-known result that the antipode of a Hopf algebra with a cocommutative coradical is bijective (see \cite[5.2.11]{Montgomery}), the antipode $S$ is bijective. Therefore, $\text{ParQSym}^{cop}$ is a Hopf algebra with antipode $\overline{S}$, where $\text{ParQSym}^{cop}=\text{ParQSym}$ as a vector space with the same product but the opposite coproduct, and $ \overline{S}=S^{-1}$.
	The Hopf algebra $\text{ParQSym}^{cop}=\bigoplus \text{ParQSym}_{i} $ is also a connected graded, then we can use Takeuchi's formula again. $\overline{S}(M_{\emptyset})=M_{\emptyset}$. For any non-empty partition diagram $\pi=\pi_{1} \otimes \pi_{2} \otimes \cdots \otimes \pi_{n}$, where $\pi_{i}$ is non-empty $ \otimes $-irreducible partition diagram for $1\leq i\leq n$,
	$$\overline{S} (M_{\pi})=\mathop{\sum}\limits_{k\geq0}(-1)^{k}(\mathop{\sum}\limits_{r\geq 0}M_{\pi_{k}^{1} \bullet \pi_{k-1}^{1} \bullet \cdots \bullet \pi_{1}^{1} \otimes \pi_{k}^{2} \bullet \pi_{k-1}^{2} \bullet \cdots \bullet \pi_{1}^{2} \otimes \cdots \otimes \pi_{k}^{r} \bullet \pi_{k-1}^{r} \bullet \cdots \bullet \pi_{1}^{r}}),$$ where $(\pi_{1}^{1},\pi_{1}^{2},\ldots,\pi_{1}^{r},\pi_{2}^{1},\ldots,\pi_{2}^{r},\ldots,\pi_{k}^{1},\ldots,\pi_{k}^{r})\hat{}=(\pi_{1},\pi_{2},\ldots,\pi_{n})$ and $\pi_{k}^{s} \bullet \pi_{k-1}^{s} \bullet \cdots \bullet \pi_{1}^{s}$ and $\pi_{i}^{1} \otimes \pi_{i}^{2} \otimes \cdots \otimes \pi_{i}^{r}$ are non-empty for $1\leq s\leq r$, $1\leq i\leq k $.
	\subsection{ParQSym is a Combinatorial Hopf algebra}\label{Combinatorial Hopf algebra}
	A  {\it Combinatorial Hopf Algebra}  is a connected graded Hopf algebra $\mathcal{H}$ over a field $\mathbb{K}$, together with a  {\it character} $\zeta:\mathcal{H}\rightarrow \mathbb{K}$ (i.e. an algebra map).
	In  \cite{ref4}, for any composition $\alpha$,
	\begin{equation*}
		\zeta_{\text{QSym}}(M_{\alpha})=\begin{cases}
			1 & \text{if } \alpha = \left( n\right)\text{ or } \alpha = \left( \right) ,\\
			0 & \text{otherwise}.
		\end{cases}
	\end{equation*}
	Similarly, we define $\zeta_{\text{ParQSym}}$ for any partition diagram $\pi$ as 
	\begin{equation*}
		\zeta_{\text{ParQSym}}(M_{\pi})=\begin{cases}
			1 & \text{if }\pi\text{ is } \otimes -\text{irreducible}  ,\\
			0 & \text{otherwise}.
		\end{cases}		
	\end{equation*}
	\begin{Proposition}
		With the character $\zeta_{\text{ParQSym}}$, $(\text{ParQSym}, \zeta_{\text{ParQSym}})$ is a Combinatorial Hopf Algebra.
	\end{Proposition}
	\begin{Proof}
		For the empty diagram, $$\zeta_{\text{ParQSym}}(M_{\emptyset})=1 .$$ It suffices to prove that for any partition diagrams $\pi$ and $\rho$,
		\begin{equation}\label{character}
			\zeta_{\text{ParQSym}}(M_{\pi})\zeta_{\text{ParQSym}}(M_{\rho}) = \zeta_{\text{ParQSym}}(M_{\pi} \star  M_{\rho}).
		\end{equation} This equation holds when $\pi=\emptyset$ or $\rho=\emptyset$. If $\pi$ and $\rho$ are both non-empty,
		the left side of (\ref{character}) is	\begin{equation*}
			\zeta_{\text{ParQSym}}(M_{\pi})\zeta_{\text{ParQSym}}(M_{\rho}) =\begin{cases}
				1 & \text{if } l(\pi)=l(\rho)=1,\\ 
				0 & \text{otherwise}.
			\end{cases}
		\end{equation*}
		By  Definition~\ref{product},   the partition diagrams appearing in $M_{\pi} \star  M_{\rho} $ have lengths no less than $\max\left\lbrace l(\pi), l(\rho)\right\rbrace$, so  $\zeta_{\text{ParQSym}}(M_{\pi} \star  M_{\rho})\neq 0$ only when $l(\pi)=l(\rho)=1 $. In this case, $$\zeta_{\text{ParQSym}}(M_{\pi} \star  M_{\rho})=\zeta_{\text{ParQSym}}(M_{\pi \otimes \rho}+M_{\pi \bullet \rho}+M_{\rho \otimes \pi})=\zeta_{\text{ParQSym}}(M_{\pi \bullet \rho})=1.$$
		The last equation is from Theorem~\ref{irr}. Therefore, we have equality (\ref{character}). 	$\qedsymbol$
	\end{Proof}
	A fundamental result on CHAs due to Aguiar, Bergeron, and Sottile, is reformulated below.
	\begin{Proposition}\cite[Theorem 4.1]{ref4}
		For any CHA $(\mathcal{H},\zeta)$, there exists a unique morphism of CHAs
		$$\Psi:(\mathcal{H},\zeta)\rightarrow(\text{QSym},\zeta_{\text{QSym}}).$$
		The map $\Psi$ is given as follows. For $h\in \mathcal{H}_{n}$,$$\Psi(h)=\mathop{\sum}\limits_{\alpha\vDash n}\zeta_{\alpha}(h)M_{\alpha},$$where for $\alpha=(\alpha_{1},\ldots,\alpha_{k})$, $\zeta_{\alpha} $ is the composite$$\mathcal{H}\stackrel{\bigtriangleup^{(k-1)}}{\longrightarrow}\mathcal{H}^{ \otimes  k}\twoheadrightarrow \mathcal{H}_{\alpha_{1}} \otimes \cdots \otimes  \mathcal{H}_{\alpha_{k}}\stackrel{\zeta^{ \otimes  k}}{\longrightarrow}\mathbb{K}, $$ with the unlabeled map being the tensor product of the canonical projections onto the
		homogeneous components $\mathcal{H}_{\alpha_{i}} $ for $1\leq i\leq k$.
	\end{Proposition}

 Using the above proposition, we can get the unique morphism of CHAs $\Psi_{PQ}: \text{ParQSym}\rightarrow \text{QSym}$ as follows:
	$\Psi_{PQ}(M_{\emptyset})=M_{()}$. For any non-empty partition diagrams $\rho=\rho_{1} \otimes \rho_{2} \otimes  \cdots  \otimes \rho_{n}$,  where $ \rho_{i}$ is a non-empty $  \otimes $-irreducible partition diagram with order $\alpha_{i}$  for $1\leq i\leq n$, and let $\alpha_{\rho}=(\alpha_{1},\ldots,\alpha_{n})$,
	\begin{equation}\label{quotient map}
		\Psi_{PQ}(M_{\rho})=M_{\alpha_{\rho}}.
	\end{equation}
	 $\Psi_{PQ}$ is surjective and can be seen as the quotient map by identifying all $  \otimes $-irreducible partition diagrams of the same order.

	John M. Campbell~\cite[Section 3.4]{ref1} defined an injective graded Hopf algebra map 
	$\Phi: \text{NSym}\rightarrow \text{ParSym}$    so that 
	\begin{equation}\label{49}
		\Phi(H_{(n)})=H_{\begin{tikzpicture}
				[
				dot/.style={circle, fill=black, inner sep=3pt}, 
				bend angle=30, 
				decoration={brace, mirror, amplitude=10pt} 
				]
				\node[dot] (A) at (0,0) {}; 
				\node[dot] (B) at (1,0) {}; 
				\node[dot] (C) at (2,0) {};
				\node[dot] (D) at (3,0) {};
				\node[dot] (E) at (4,0) {}; 
				\node[dot] (A') at (0,1) {}; 
				\node[dot] (B') at (1,1) {}; 
				\node[dot] (C') at (2,1) {}; 
				\node[dot] (D') at (3,1) {};
				\node[dot] (E') at (4,1) {};
				\draw (A) to[bend left] (B);
				\draw (B) to[bend left] (C);
				\draw (D) to[bend left] (E);
				\node at (2.5,0.5) {$\cdots$};
				
		\end{tikzpicture}}=H_{\pi_{(n)}},
	\end{equation}
	where $$\pi_{(n)}=\left\lbrace \left\lbrace 1\right\rbrace ,\left\lbrace 2\right\rbrace ,\ldots,\left\lbrace n\right\rbrace ,\left\lbrace 1',2',\ldots,n'\right\rbrace \right\rbrace$$ as a set partition, 
	with (\ref{49}) extended linearly multiplicatively.

 For any composition $\alpha=(\alpha_{1},\alpha_{2},\ldots,\alpha_{k})$, let
	\begin{equation*}
		\begin{aligned}
			\pi_{\alpha}&=\pi_{(\alpha_{1})} \otimes \pi_{(\alpha_{1})} \otimes \cdots \otimes \pi_{(\alpha_{k})}\\
			&=	\begin{tikzpicture}
				[baseline=(current bounding box.east),
				dot/.style={circle, fill=black, inner sep=3pt}, 
				bend angle=30, 
				decoration={brace, mirror, amplitude=10pt} 
				]
				\node[dot] (A) at (0,0) {}; 
				\node[dot] (B) at (1,0) {}; 
				\node[dot] (C) at (2,0) {};
				\node[dot] (D) at (3,0) {};
				\node[dot] (E) at (4,0) {};
				\node[dot] (F) at (5,0) {}; 
				\node[dot] (G) at (6,0) {}; 
				\node[dot] (H) at (7,0) {};
				\node[dot] (A') at (0,1) {}; 
				\node[dot] (B') at (1,1) {}; 
				\node[dot] (C') at (2,1) {}; 
				\node[dot] (D') at (3,1) {};
				\node[dot] (E') at (4,1) {};
				\node[dot] (F') at (5,1) {}; 
				\node[dot] (G') at (6,1) {}; 
				\node[dot] (H') at (7,1) {};
				\draw (A) to[bend left] (B);
				\draw (C) to[bend left] (D);
				\node at (1.5,0.5) {$\cdots$};
				\node at (4.5,0.5) {$\cdots$};
				\node at (5.5,0.5) {$\cdots$};
				\node at (6.5,0.5) {$\cdots$};
				
				\draw [decorate] 
				(A.south) -- (D.south) node [midway, below=10pt] {$\pi_{\alpha_{1}}$};
				\draw [decorate] 
				(E.south) -- (F.south) node [midway, below=10pt] {$\pi_{\alpha_{2}}$};
				\draw [decorate] 
				(G.south) -- (H.south) node [midway, below=10pt] {$\pi_{\alpha_{k}}$};
			\end{tikzpicture}\\
		\end{aligned}
	\end{equation*}
	Then we have  $\Phi(H_{\alpha})=H_{\pi_{\alpha}}$.
	\text{NSym} is the graded dual of \text{ParSym}, so there is a bilinear form $\left\langle \cdot,\cdot\right\rangle :\text{QSym}  \otimes  \text{NSym} \rightarrow \mathbb{K}$ such that
	$$\left\langle M_{\alpha},H_{\beta}\right\rangle=\delta_{\alpha,\beta}$$for all compositions $\alpha$ and $\beta$. Then we have for any compositions $\alpha$ and $\beta$, $$\alpha_{\pi_{\beta}}=\beta,\left\langle M_{\pi_{\beta}}, \Phi(H_{\alpha})\right\rangle =\left\langle \Psi_{PQ}(M_{\pi_{\beta}}),H_{\alpha}\right\rangle=\delta_{\beta,\alpha} .$$
	
	\subsection{Subcoalgebras and Hopf subalgebras of ParQSym}\label{Subcoalgebras and Hopf subalgebras}
	By analogy with  \cite{ref1}, we find several subcoalgebras and Hopf subalgebras of \text{ParQSym}.
		Partition diagrams that can be expressed as planar graphs are called  {\it planar diagrams}.
	It is easy to see the following two lemmas.
	\begin{Lemma}\label{planar1}
		For a partition diagram $\pi=\pi_{1} \otimes \pi_{2} \otimes \cdots \otimes \pi_{l}$ ($\pi_{1}$, $\pi_{2}\cdots\pi_{l}$ are not necessarily $ \otimes $-irreducible),  $\pi$ is a planar diagram if and only if $\pi_{1}$, $\pi_{2} $,$\cdots$, $\pi_{l} $ are all planar diagrams.
	\end{Lemma}
	\begin{Lemma}\label{planar2}
		If $\rho$ and $\sigma$ are planar diagrams, then $\rho \bullet \sigma$ is planar. 
	\end{Lemma}
	\begin{Theorem}\label{planar} (Planar Subalgebra)
		The graded vector subspace of \text{ParQSym} spanned by planar diagrams
		constitutes a Hopf subalgebra of \text{ParQSym}.
	\end{Theorem}
	\begin{Proof}
		The target vector subspace is closed under the coproduct, thus it is a subcoalgebra of \text{ParQSym} by Lemma~\ref{planar1}.

	For any non-empty partition diagrams $\rho=\rho_{1} \otimes \rho_{2} \otimes  \cdots  \otimes \rho_{n}$, $\sigma=\sigma_{1} \otimes \sigma_{2} \otimes  \cdots  \otimes \sigma_{m}$, where $ \rho_{i}$, $\sigma_{j}$ are non-empty $  \otimes $-irreducible partition diagrams for $1\leq i\leq n$, $1\leq j\leq m$, $$M_{\rho}  \star  M_{\sigma}=\mathop{\sum}M_{\rho'_{1} \bullet \sigma'_{1} \otimes \cdots \otimes \rho'_{k} \bullet \sigma'_{k}},$$ where the sum is over all $ k\in \mathbb{N}$ and all pairs of ordered sequences $\left( \rho'_{1},\ldots,\rho'_{k}\right)$ and $\left( \sigma'_{1},\ldots,\sigma'_{k}\right)$ such that $\left( \rho'_{1},\ldots,\rho'_{k}\right)\hat{}=\left( \rho_{1},\ldots,\rho_{n}\right) $ and \linebreak[4]$\left( \sigma'_{1},\ldots,\sigma'_{k}\right)\hat{}=\left( \sigma_{1},\ldots,\sigma_{m}\right) $, and such that $\rho'_{s} \bullet \sigma'_{s}$ is non-empty for all $1\leq s\leq k$. If $\rho$ and $\sigma$ are planar, then $ \rho_{i}$, $\sigma_{j}$ are planar for $1\leq i\leq n$, $1\leq j\leq m$ by Lemma~\ref{planar1}. Adding that $\emptyset$ is also planar, we have that all $ \rho'_{i}$ and $\sigma'_{j} $ are planar. It follows that all $\rho'_{i} \bullet  \sigma'_{i}$ are planar by Lemma~\ref{planar2}. Hence, all $\rho'_{1} \bullet \sigma'_{1} \otimes \cdots \otimes \rho'_{k} \bullet \sigma'_{k}$ are planar by Lemma~\ref{planar1}.
		We have that the target vector subspace is closed under the product $ \star $ and contains the unit $\emptyset$, thus it is a subalgebra of \text{ParQSym}. Therefore, it is a Hopf subalgebra of \text{ParQSym}.$\qedsymbol$
	\end{Proof}
		A block in a set partitions is  {\it propagating} if it has at least one upper vertex and at least one lower vertex. The  {\it propagation number} of a partition diagram $\pi$ refers to the number of components in $\pi$ that contain at least one upper vertex and at least one lower vertex.
	\begin{Theorem}
		The graded vector subspace of \text{ParQSym} spanned by partition diagrams with propagation number 0 forms a Hopf subalgebra.
	\end{Theorem}
	\begin{Proof}
		For a partition diagram $\pi=\pi_{1} \otimes \pi_{2} \otimes \cdots \otimes \pi_{l}$, it is easy to see the propagation number of $\pi$ is 0 if and only if those of $\pi_{1}$, $\pi_{2} $,$\cdots$, $\pi_{l} $ are 0. If $\rho$ and $\sigma$ contain no propagating blocks, then $\rho \bullet \sigma $  does not contain any propagating blocks by the definition of $ \bullet $. 
		In a similar process to Theorem~\ref{planar}, we have that the target vector subspace is a Hopf subalgebra of \text{ParQSym}.$\qedsymbol$
	\end{Proof}
	\begin{Theorem}
		The graded vector subspace of \text{ParQSym} spanned by partition diagrams whose upper nodes are all isolated (including $\emptyset$) forms a Hopf subalgebra.
	\end{Theorem}
	\begin{Proof}
		For a partition diagram $\pi=\pi_{1} \otimes \pi_{2} \otimes \cdots \otimes \pi_{l}$, it is easy to see the upper nodes of $\pi$ are all isolated if and only if those of $\pi_{1}$, $\pi_{2} $,$\cdots$, $\pi_{l} $ are all isolated. If the upper nodes of $\rho$ and $\sigma$ are all isolated, then the upper nodes of $\rho \bullet \sigma $ are all isolated by the definition of $ \bullet $. 
		Similar to Theorem~\ref{planar},  the target vector subspace is a Hopf subalgebra of \text{ParQSym}.$\qedsymbol$
	\end{Proof}

 A  {\it matching} is a partition diagram $\pi$ such that all blocks in $\pi$	are of size at most two.
		A  {\it perfect matching} is a matching such that each block is of size two. 
		A  {\it permuting diagram} is a partition diagram $\pi$ such that each block of $\pi$ is of size two and propagating. A permuting diagram can be written as $\left\lbrace\left\lbrace 1,p(1)'\right\rbrace ,\left\lbrace 2,p(2)'\right\rbrace ,\ldots,\left\lbrace k,p(k') \right\rbrace  \right\rbrace $ for some permutation $p$.
	%
		A  {\it partial permutation} is a partition diagram $\pi$ such that each block of $\pi$ is of size one or two and such that each block of size two in $\pi$ is propagating.
	The above notions can all form subcoalgebras of \text{ParQSym}.
	\begin{Theorem}
		The graded vector subspace of \text{ParQSym} spanned by matchings (perfect matchings,  permuting diagrams, or partial permutations) forms a subcoalgebra.
	\end{Theorem}
	\begin{Proof}
		Every block in $\pi=\rho \otimes \sigma$ comes from either $\rho$ or $\sigma$. By the definition of the coproduct, we have that the target vector spaces are closed under the coproduct. $\qedsymbol$
	\end{Proof}
	
	\section{\texorpdfstring{$R$-basis of ParSym and $L$-basis of ParQSym}{R-basis of ParSym and L-basis of ParQSym}}\label{$R$-basis of ParSym and $L$-basis of ParQSym}
	In this section, we define the $R$-basis of \text{ParSym} and the $L$-basis of \text{ParQSym}, by analogy with the $R$-basis of \text{NSym} and the $L$-basis of \text{QSym}.

	For any compositions $\alpha=(\alpha_{1},\alpha_{2},\ldots,\alpha_{n})$ and $\beta=(\beta_{1},\beta_{2},\ldots,\beta_{m})$, define  $$\alpha\cdot\beta=(\alpha_{1},\ldots,\alpha_{n},\beta_{1},\ldots,\beta_{m}),$$ $$\alpha\odot\beta=(\alpha_{1},\alpha_{2},\ldots,\alpha_{n}+\beta_{1},\beta_{2},\ldots,\beta_{m}).$$ If $\alpha$ and $\beta$ are both compositions of $k\in \mathbb{N}$, say that $\alpha$  {\it refines} $\beta$ or $\beta$  {\it coarsens} $\alpha$ if $\alpha$ can be written as a concatenation
	\begin{equation}\label{refine part}
		\alpha=\alpha_{\beta}^{(1)}\cdot\alpha_{\beta}^{(2)}\cdot\cdots\cdot\alpha_{\beta}^{(l(\beta))}
	\end{equation}
	where $\alpha_{\beta}^{(i)}\vDash \beta_{i}$ for each $i$.

 Now we define ``refinement'' of partition diagrams.
	\begin{Definition}
		For non-empty partition diagrams $\pi=\pi_{1}  \divideontimes^{\pi}_{1} \pi_{2}  \divideontimes^{\pi}_{2} \cdots  \divideontimes^{\pi}_{n-1} \pi_{n}$ and $\sigma=\sigma_{1}  \divideontimes^{\sigma}_{1} \sigma_{2}  \divideontimes^{\sigma}_{2} \cdots  \divideontimes^{\sigma}_{m-1} \sigma_{m}$, where $ \divideontimes^{\pi}_{i},  \divideontimes^{\sigma}_{j}\in\left\lbrace  \otimes, \bullet \right\rbrace$, and where $\pi_{1},\ldots,\pi_{n}$ and $\sigma_{1},\ldots,\sigma_{m}$ are non-empty $ \otimes $-irreducible and $ \bullet $-irreducible (every non-empty partition diagram can be uniquely written in this form by Theorems~\ref{length unique} and \ref{irr2}. To see this, one can use \ref{length unique} first and then use \ref{irr2} to every $ \otimes $-irreducible part), we say $\pi\sim \sigma$ if the ordered sequences $(\pi_{1},\pi_{2},\ldots,\pi_{n})$ and $(\sigma_{1},\sigma_{2},\ldots,\sigma_{m}) $ are the same. Let $$S(\pi):=\left\lbrace i\in [n-1]\mid\divideontimes^{\pi}_{i} \text{ is } \otimes  \right\rbrace .$$ We say $\sigma$  {\it coarsens} $\pi$ or $\pi$  {\it refines} $\sigma$, denoted $\pi\leq\sigma$, if $\pi\sim\sigma$ and $S(\sigma)\subseteq S(\pi)$. From the definition of the length of a partition diagram, it is obvious that $l(\pi)=1+\#S(\pi)$. We say that $\emptyset$  {\it coarsens} $\emptyset$ and $\emptyset$  {\it refines} $\emptyset$ by convention.
	\end{Definition}
	The following simple lemmas will be useful later.
	\begin{Lemma}\label{poset}
		The relation $\sim$ is an equivalence relation,
		and the refinement relation defined above is a partial order.
	\end{Lemma}
	\begin{Proof}
		Directly from the definitions of the two relations.$\qedsymbol$
	\end{Proof}
	\begin{Lemma}\label{repeat}
		If $\rho'\leq\rho$ and $\sigma' \leq\sigma$, then $\rho' \bullet \sigma'\leq\rho \bullet \sigma$, and $\rho' \otimes \sigma'\leq\rho \otimes \sigma$.
	\end{Lemma}
	\begin{Proof}
		For $\rho=\rho_{1}  \divideontimes^{\rho}_{1} \rho_{2}  \divideontimes^{\rho}_{2} \cdots  \divideontimes^{\rho}_{n-1} \rho_{n}$ and $\sigma=\sigma_{1}  \divideontimes^{\sigma}_{1} \sigma_{2}  \divideontimes^{\sigma}_{2} \cdots  \divideontimes^{\sigma}_{m-1} \sigma_{m}$, where $ \divideontimes^{\rho}_{i},  \divideontimes^{\sigma}_{j}\in\left\lbrace  \otimes , \bullet \right\rbrace$, and $\rho_{1},\rho_{2},\ldots,\rho_{n}$ and $ \sigma_{1},\sigma_{2},\ldots,\sigma_{m}$ are non-empty $ \otimes $-irreducible and $ \bullet $-irreducible, we have
		$$S(\rho \bullet \sigma) =S(\rho)\cup \left\lbrace i+n|i\in S(\sigma)\right\rbrace , $$
		$$S(\rho' \bullet \sigma') =S(\rho')\cup \left\lbrace i+n|i\in S(\sigma')\right\rbrace , $$
		$$S(\rho \otimes \sigma) =S(\rho)\cup\left\lbrace n\right\rbrace \cup \left\lbrace i+n|i\in S(\sigma)\right\rbrace , $$
		$$S(\rho' \otimes \sigma') =S(\rho')\cup\left\lbrace n\right\rbrace \cup \left\lbrace i+n|i\in S(\sigma')\right\rbrace . $$
		Then
		$$S(\rho)\subseteq S(\rho'),S(\sigma)\subseteq S(\sigma')\Rightarrow S(\rho \otimes \sigma)\subseteq S(\rho' \otimes \sigma'), S(\rho \bullet \sigma)\subseteq S(\rho' \bullet \sigma').  $$
		$\qedsymbol$
	\end{Proof}
	\subsection{\texorpdfstring{$R$-basis of ParSym}{R-basis of ParSym}}
	
	There is another $\mathbb{K}$-basis of \text{NSym} called the  {\it non-commutative ribbon functions} $\left\lbrace R_{\alpha}\right\rbrace $ \cite{ref2}:$$R_{\alpha}=\mathop{\sum}\limits_{\beta\text{ coarsens } \alpha}(-1)^{l(\beta)-l(\alpha)}H_{\beta}, $$
	for any composition $\alpha$. The product   is given by $$ R_{\alpha}R_{\beta}=R_{\alpha\cdot\beta}+R_{\alpha\odot\beta}$$if $\alpha$ and $\beta$ are non-empty, and $R_{()}$ is the multiplicative identity \cite[Theorem 5.4.10]{ref2}.
	By analogy with $\left\lbrace R_{\alpha}\right\rbrace $, we define the following basis of $\text{ParSym} $:
	Let $$R_{\pi}:=\mathop{\sum}\limits_{\sigma\text{ coarsens } \pi}(-1)^{l(\sigma)-l(\pi)}H_{\sigma}, $$then  $\left\lbrace R_{\pi}\right\rbrace_{\pi\in A_{i}, i\geq 0}$ forms a basis of $\text{ParSym} $.
	\begin{Proposition}
		The product of the $R$-basis of \text{ParSym} is given by $$R_{\pi}R_{\sigma}=R_{\pi \otimes \sigma}+R_{\pi \bullet \sigma}$$if $\pi$ and $\sigma$ are non-empty and $R_{\emptyset}$ is the multiplicative identity.
	\end{Proposition} 
	\begin{Proof}
		Consider the empty diagram,
		$$R_{\emptyset}=H_{\emptyset},$$ so we get a unit.
		For any non-empty partition diagrams $\pi$ and $\sigma$, 
		\begin{equation*}
			\begin{split}
				l(\pi \otimes \sigma)&=1+\#S(\pi \otimes \sigma)\\
				&=1+(\#S(\pi)+1+\#S(\sigma))\\
				&=1+(l(\pi)-1+1+l(\sigma)-1)\\
				&=l(\pi)+l(\sigma),\\
				l(\pi \bullet \sigma)&=1+\#S(\pi \bullet \sigma)\\
				&=1+(\#S(\pi)+\#S(\sigma))\\
				&=1+(l(\pi)-1+l(\sigma)-1)\\
				&=l(\pi)+l(\sigma)-1.\\
			\end{split}
		\end{equation*}
		\begin{equation*}
			\begin{split}
				&R_{\pi}R_{\sigma}\\&=(\mathop{\sum}\limits_{\pi' \text{ coarsens } \pi}(-1)^{l(\pi')-l(\pi)}H_{\pi'})(\mathop{\sum}\limits_{\sigma'\text{ coarsens } \sigma}(-1)^{l(\sigma')-l(\sigma)}H_{\sigma'})\\
				&=\mathop{\sum}\limits_{\pi'\text{ coarsens } \pi,\sigma'\text{ coarsens } \sigma}(-1)^{(l(\pi')+l(\sigma'))-(l(\pi)+l(\sigma))}H_{\pi'}H_{\sigma'}\\&=\mathop{\sum}\limits_{\pi'\text{ coarsens } \pi,\sigma'\text{ coarsens } \sigma}(-1)^{(l(\pi' \otimes \sigma'))-(l(\pi \otimes \sigma))}H_{\pi' \otimes \sigma'}\\
				&=\mathop{\sum}\limits_{\rho\text{ coarsens } \pi \otimes  \sigma}(-1)^{l(\rho)-l(\pi \otimes \sigma)}H_{\rho}-\mathop{\sum}\limits_{\rho'\text{ coarsens } \pi \bullet  \sigma}(-1)^{l(\rho')-l(\pi \otimes \sigma)}H_{\rho'}\\
				&=R_{\pi \otimes \sigma}+\mathop{\sum}\limits_{\rho'\text{ coarsens } \pi \bullet  \sigma}(-1)^{l(\rho')-l(\pi \bullet \sigma)}H_{\rho'}\\
				&=R_{\pi \otimes \sigma}+R_{\pi \bullet \sigma}.
			\end{split}
		\end{equation*}
		$\qedsymbol$
	\end{Proof}
	When $\sigma$ is $ \otimes $-irreducible, we have that $$R_{\sigma}=H_{\sigma}$$and 
	$$	\bigtriangleup R_{\sigma}=	\bigtriangleup H_{\sigma}\\=\mathop{\sum}\limits_{\sigma_{1} \bullet \sigma_{2}=\sigma}H_{\sigma_{1}} \otimes  H_{\sigma_{2}}\\\\=\mathop{\sum}\limits_{\sigma_{1} \bullet \sigma_{2}=\sigma}R_{\sigma_{1}} \otimes  R_{\sigma_{2}}.$$
	The coproduct can be defined recursively by the above equation and 
	$$\bigtriangleup R_{\pi \otimes \sigma}=(\bigtriangleup R_{\pi})(\bigtriangleup R_{\sigma})-\bigtriangleup R_{\pi \bullet \sigma}.$$
	\subsection{\texorpdfstring{$L$-basis of ParQSym}{L-basis of ParQSym}}
	Dual to the $\left\lbrace R_{\alpha}\right\rbrace$, there is another $\mathbb{K}$-basis of \text{QSym}  called the  {\it fundamental quasi-symmetric functions} $\left\lbrace L_{\alpha}\right\rbrace $. $$L_{\alpha}=\mathop{\sum}\limits_{\beta\leq \alpha}M_{\beta}, $$
	for any composition $\alpha$. 
	The coproduct   is given by $$ \bigtriangleup L_{\alpha}=\mathop{\sum}\limits_{\beta\cdot\gamma=\alpha\text{ or } \beta\odot\gamma= \alpha}L_{\beta} \otimes  L_{\gamma},$$see \cite[Proposition 5.2.15]{ref2}.
	By analogy with $\left\lbrace L_{\alpha}\right\rbrace $, we define the following basis of $\text{ParQSym} $.

	Let $$L_{\pi}:=\mathop{\sum}\limits_{\sigma\leq \pi}M_{\sigma}, $$then  $\left\lbrace L_{\pi}\right\rbrace_{\pi\in A_{i}, i\geq 0}$ forms a basis of $\text{ParQSym} $.\\
	
	\begin{Proposition}
		The coproduct of $L$-basis of \text{ParQSym} is given by $$\bigtriangleup L_{\pi}=\mathop{\sum}\limits_{\rho \otimes \sigma=\pi\text{ or } \rho \bullet \sigma= \pi}L_{\rho} \otimes  L_{\sigma}.$$
	\end{Proposition}
	\begin{Proof}
		Consider the empty diagram:
		$$\bigtriangleup L_{\emptyset}=\bigtriangleup M_{\emptyset}=M_{\emptyset} \otimes  M_{\emptyset}=L_{\emptyset} \otimes  L_{\emptyset}.$$For any non-empty partition diagram $\pi$,
		\begin{equation*}
			\begin{split}
				\bigtriangleup L_{\pi}&=\mathop{\sum}\limits_{\sigma\leq\pi}\bigtriangleup M_{\sigma}\\
				&=\mathop{\sum}\limits_{\pi\sim\sigma,S(\pi)\subseteq S(\sigma)}\bigtriangleup M_{\sigma}\\
				&=\mathop{\sum}\limits_{\substack{i\in S(\sigma)\\\pi_{1}  \divideontimes_{1} \cdots  \divideontimes_{i-1} \pi_{i} \otimes \pi_{i+1}  \divideontimes_{i+1} \cdots  \divideontimes_{n-1} \pi_{n}=\sigma\\\pi\sim\sigma,S(\pi)\subseteq S(\sigma), \divideontimes_{j} \in\left\lbrace \bullet , \otimes  \right\rbrace }} M_{\pi_{1}  \divideontimes_{1} \cdots  \divideontimes_{i-1} \pi_{i}} \otimes  M_{\pi_{i+1}  \divideontimes_{i+1} \cdots  \divideontimes_{n-1} \pi_{n}}\\
				&=\mathop{\sum}\limits_{\substack{i\in S(\pi)\\\pi_{1}  \divideontimes_{1} \cdots  \divideontimes_{i-1} \pi_{i} \otimes \pi_{i+1}  \divideontimes_{i+1} \cdots  \divideontimes_{n-1} \pi_{n}=\sigma\\\pi\sim\sigma,S(\pi)\subseteq S(\sigma), \divideontimes_{j} \in\left\lbrace \bullet , \otimes  \right\rbrace }} M_{\pi_{1} \divideontimes_{1} \cdots \divideontimes_{i-1} \pi_{i}} \otimes  M_{\pi_{i+1} \divideontimes_{i+1} \cdots \divideontimes_{n-1} \pi_{n}}+\\&\mathop{\sum}\limits_{\substack{i\in S(\sigma)\setminus S(\pi)\\\pi_{1} \divideontimes_{1} \cdots \divideontimes_{i-1} \pi_{i} \otimes \pi_{i+1} \divideontimes_{i+1} \cdots \divideontimes_{n-1} \pi_{n}=\sigma\\\pi\sim\sigma,S(\pi)\subseteq S(\sigma), \divideontimes_{j} \in\left\lbrace \bullet , \otimes  \right\rbrace }} M_{\pi_{1} \divideontimes_{1} \cdots \divideontimes_{i-1} \pi_{i}} \otimes  M_{\pi_{i+1} \divideontimes_{i+1} \cdots \divideontimes_{n-1} \pi_{n}}\\
				&=\mathop{\sum}\limits_{\pi\sim\sigma_{1} \otimes \sigma_{2},S(\pi)\subseteq S(\sigma_{1} \otimes \sigma_{2}) } M_{\sigma_{1}} \otimes  M_{\sigma_{2}}+\mathop{\sum}\limits_{S(\pi)\subseteq S(\sigma_{1'} \bullet \sigma_{2'})} M_{\sigma_{1'}} \otimes  M_{\sigma_{2'}}\\
				&=\mathop{\sum}\limits_{\substack{\pi_{1} \otimes \pi_{2}=\pi\\\pi_{1}\sim\sigma_{1},S(\pi_{1})\subseteq S(\sigma_{1}),\\\pi_{2}\sim\sigma_{2},S(\pi_{2})\subseteq S(\sigma_{2})}} M_{\sigma_{1}} \otimes  M_{\sigma_{2}}+\mathop{\sum}\limits_{\substack{\pi_{1'} \bullet \pi_{2'}=\pi\\\pi_{1'}\sim\sigma_{1'},S(\pi_{1'})\subseteq S(\sigma_{1'}),\\\pi_{2'}\sim\sigma_{2'},S(\pi_{2'})\subseteq S(\sigma_{2'})}} M_{\sigma_{1'}} \otimes  M_{\sigma_{2'}}\\
				&=\mathop{\sum}\limits_{\pi_{1} \otimes \pi_{2}=\pi} L_{\pi_{1}} \otimes  L_{\pi_{2}}+\mathop{\sum}\limits_{\pi_{1'} \bullet \pi_{2'}=\pi} L_{\pi_{1'}} \otimes  L_{\pi_{2'}}.\\
			\end{split}
		\end{equation*}
		$\qedsymbol$
	\end{Proof}
	For any non-empty partition diagram $\pi=\pi_{1} \divideontimes_{1} \pi_{2} \divideontimes_{2} \cdots \divideontimes_{n-1} \pi_{n}$, where $ \divideontimes_{i} \in\left\lbrace  \otimes , \bullet \right\rbrace $ for any $i \in [n-1]$, and $\pi_{1},\pi_{2},\ldots,\pi_{n}$ are non-empty $ \otimes $-irreducible and $ \bullet $-irreducible, the above proposition means $$	\bigtriangleup L_{\pi}=L_{\emptyset} \otimes  L_{\pi}+\mathop{\sum}\limits_{1\leq i\leq n-1}L_{\pi_{1} \divideontimes_{1} \pi_{2} \divideontimes_{2} \cdots \divideontimes_{i-1} \pi_{i}}  \otimes  L_{\pi_{i+1} \divideontimes_{i+1} \cdots \divideontimes_{n-1} \pi_{n}}+L_{\pi} \otimes  L_{\emptyset}$$
	\begin{Proposition}
		$\left\{L_{\pi}\right\}$ in \text{ParQSym} is the dual basis of $\left\{R_{\pi}\right\}$ in $\text{ParSym}$.
	\end{Proposition}
	\begin{Proof}
		For any partition diagrams $\pi$ any $\rho$,
		\begin{equation*}
			\begin{split}
				\left\langle L_{\rho},R_{\pi}\right\rangle&=	\left\langle\mathop{\sum}\limits_{\rho'\leq \rho}M_{\rho'},\mathop{\sum}\limits_{\pi'\geq \pi}H_{\pi'}\right\rangle\\&=\mathop{\sum}\limits_{\substack{\rho'\leq \rho\\\pi'\geq \pi}}(-1)^{l(\pi')-l(\pi)}\delta_{\rho',\pi'},
			\end{split}
		\end{equation*}
		$$\left\langle L_{\rho},R_{\emptyset}\right\rangle=\mathop{\sum}\limits_{\rho'\leq \rho}\delta_{\rho',\emptyset}=\delta_{\rho,\emptyset}$$
		$$\left\langle L_{\emptyset},R_{\pi}\right\rangle=\mathop{\sum}\limits_{\pi'\text{ coarsens } \pi}(-1)^{l(\pi')-l(\pi)}\delta_{\emptyset,\pi'}=\delta_{\emptyset,\pi'}. $$
		For any non-empty $\rho$ and $\pi$, if $\rho=\rho_{1} \divideontimes^{\rho}_{1} \cdots \divideontimes^{\rho}_{n-1} \rho_{n}$, $\pi=\pi_{1} \divideontimes^{\pi}_{1} \cdots \divideontimes^{\pi}_{m-1} \pi_{m}$, where  $\rho_{1}\cdots\rho_{n},\pi_{1},\pi_{2},\ldots,\pi_{m}$ are non-empty $ \otimes $-irreducible and $ \bullet $-irreducible,and $ \divideontimes^{\rho}_{i},  \divideontimes^{\pi}_{j}\in\left\lbrace \bullet , \otimes  \right\rbrace $ for all possible $i$and $j$, then $\rho'$ and $\pi' $ can also be expressed by  $\rho_{1} \divideontimes^{\rho'}_{1} \cdots \divideontimes^{\rho'}_{n-1} \rho_{n}$ and  $\pi_{1} \divideontimes^{\pi'}_{1} \cdots \divideontimes^{\pi'}_{m-1}\pi_{m}$ respectively, where $ \divideontimes^{\rho'}_{i},  \divideontimes^{\pi'}_{j}\in\left\lbrace \bullet , \otimes  \right\rbrace $ for all possible $i$ and $j$. If there exists $\rho'$ and $\pi'$ such that $\delta_{\rho',\pi'}\neq 0 $, then $n=m$ and $\rho_{i}=\pi_{i}$ for $1\leq i\leq n$. In this case,
		\begin{equation*}
			\begin{split}
				\left\langle L_{\rho},R_{\pi}\right\rangle&=\mathop{\sum}\limits_{\substack{\rho\sim\rho',S(\rho)\subseteq S(\rho')\\\pi'\sim\pi,S(\pi')\subseteq S(\pi)}}(-1)^{\#S(\pi')-\#S(\pi)}\delta_{\rho',\pi'}\\&=\mathop{\sum}\limits_{\rho\sim\pi'\sim\pi,S(\rho)\subseteq S(\pi')\subseteq S(\pi)}(-1)^{\#S(\pi')-\#S(\pi)}.
			\end{split}
		\end{equation*}
		If $\left\langle L_{\rho},R_{\pi}\right\rangle\neq 0$, then $S(\rho)\subseteq S(\pi)$, and 
		\begin{equation*}
			\begin{split}
				\left\langle L_{\rho},R_{\pi}\right\rangle&=(-1)^{\#S(\rho)-\#S(\pi)}\mathop{\sum}\limits_{\rho\sim\pi'\sim\pi,S(\rho)\subseteq S(\pi')\subseteq S(\pi)}(-1)^{\#S(\pi')-\#S(\rho)}\\&=(-1)^{\#S(\rho)-\#S(\pi)}\mathop{\sum}\limits_{S\subseteq S(\pi)\setminus S(\rho)}(-1)^{\#S}\\&=\delta_{S(\pi)\setminus S(\rho),\emptyset}\\&=\delta_{\rho,\pi}.
			\end{split}
		\end{equation*} 
		We conclude that $$\left\langle L_{\rho},R_{\pi}\right\rangle=\delta_{\rho,\pi}.$$ $\qedsymbol$
	\end{Proof}
	For any composition $\alpha=(\alpha_{1},\alpha_{2},\ldots,\alpha_{n})$, let $$D(\alpha):=\left\lbrace \alpha_{1}, \alpha_{1}+\alpha_{2},\ldots,\alpha_{1}+\alpha_{2}+\cdots+\alpha_{n-1}\right\rbrace  .$$
	For any labeled linear order of  integers $\omega=(\omega_{1}<_{\omega}\omega_{2}<_{\omega}\cdots<_{\omega}\omega_{n})$, let $$ Des(\omega):=\left\lbrace i|\omega_{i}>_{\mathbb{Z}}\omega_{i+1}\right\rbrace .$$ 
	The product of the L-basis of \text{QSym} is given as follows:
	\begin{Proposition}[\protect{\cite[Proposition 5.2.15]{ref2}}]
		For any compositions $\alpha$ and $\beta$, $$L_{\alpha}L_{\beta}=\mathop{\sum}\limits_{\omega\in \omega_{\alpha}\shuffle\omega_{\beta}}L_{\gamma(\omega)}$$
		where $\omega_{\alpha}$ is any labeled linear order with underlying set $\left\lbrace 1, 2, . . . , |\alpha|\right\rbrace $ such that $Des (\omega_{\alpha}) = D (\alpha)$, $\omega_{\beta}$ is any labeled linear order with underlying set $\left\lbrace |\alpha| + 1, |\alpha| + 2, . . . , |\alpha| + |\beta|\right\rbrace $ such that $Des (\omega_{\beta}) =
		D (\beta)$,  and $\gamma(\omega)$ is the unique composition of $|\alpha| + |\beta|$ with $D(\gamma(\omega)) = Des(\omega)$.
	\end{Proposition}
	Similarly, we have a product on the basis $\left\lbrace L_{\pi} \right\rbrace $ for \text{ParQSym}.
	\begin{Proposition}
		The product of the basis $\left\lbrace L_{\pi} \right\rbrace $ for \text{ParQSym}  is given by
		\begin{equation}\label{Lproduct1}
			L_{\rho} \star  L_{\sigma}=\mathop{\sum}\limits_{k\geq 0}L_{\rho_{1} \bullet \sigma_{1} \otimes \rho_{2} \bullet \sigma_{2} \otimes \cdots \otimes \rho_{k} \bullet \sigma_{k}}
		\end{equation}for any non-empty partition diagrams $\rho$ and $\sigma$,where the sum is over all $ k\in \mathbb{N}$ and all pairs of ordered sequences $\left( \rho_{1},\ldots,\rho_{k}\right)$ and $\left( \sigma_{1},\ldots,\sigma_{k}\right)$ of (possibly $ \otimes $-reducible) partition diagrams such that $\rho$ can be expressed as $\rho_{1} \divideontimes^{\rho}_{1} \cdots \divideontimes^{\rho}_{k-1} \rho_{k}$  and $\sigma$ can be expressed as $\sigma_{1} \divideontimes^{\sigma}_{1} \cdots \divideontimes^{\sigma}_{k-1} \sigma_{k}$ ($ \divideontimes^{\rho}_{i},  \divideontimes^{\sigma}_{j} \in\left\lbrace  \otimes , \bullet \right\rbrace$), and such that $\sigma_{1},\ldots,\sigma_{k-1}$ and $\rho_{2},\ldots,\rho_{k}$ are non-empty. Accordingly, we have that $L_{\rho} \star  L_{\emptyset}=L_{\emptyset} \star  L_{\rho}=L_{\rho}$ for any partition diagram $\rho$.
	\end{Proposition} 
	\begin{Proof}
		Consider the empty diagram: $$	L_{\emptyset}=\mathop{\sum}\limits_{\rho'\leq \emptyset}M_{\rho'} =M_{\emptyset}, $$so $L_{\emptyset} $ is the unit of $ \star $.\\
		For any non-empty partition diagrams $\rho$ and $\sigma$,
		\begin{equation}\label{Lproduct2}
			\begin{split}
				L_{\rho}  \star  L_{\sigma}&=(\mathop{\sum}\limits_{\rho'\leq \rho}M_{\rho'}) \star (\mathop{\sum}\limits_{\sigma'\leq \sigma}M_{\sigma'}) \\
				&=\mathop{\sum}\limits_{\substack{\rho'\leq \rho\\\sigma'\leq \sigma}}(\mathop{\sum}\limits_{k\geq 0}M_{\rho'_{1} \bullet \sigma'_{1} \otimes \rho'_{2} \bullet \sigma'_{2} \otimes \cdots \otimes \rho'_{k} \bullet \sigma'_{k}})
			\end{split}
		\end{equation}
		where for any pair $(\rho',\sigma')$, the sum is over all $ k\in \mathbb{N}$ and all pairs of ordered sequences $\left( \rho'_{1},\ldots,\rho'_{k}\right)$ and $\left( \sigma'_{1},\ldots,\sigma'_{k}\right)$ of $ \otimes $-irreducible (possibly empty) partition diagrams such that $\rho'_{1} \otimes \cdots \otimes \rho'_{k}=\rho'$ and $\sigma'_{1} \otimes \cdots \otimes \sigma'_{k}=\sigma', $ and such that $\rho'_{s} \bullet \sigma'_{s}\neq\emptyset$ for all $1\leq s\leq k$.
		
		We need to show that every $\rho'_{1} \bullet \sigma'_{1} \otimes \rho'_{2} \bullet \sigma'_{2} \otimes \cdots \otimes \rho'_{k} \bullet \sigma'_{k} $ in (\ref{Lproduct2}) refines some $\rho_{1} \bullet \sigma_{1} \otimes \rho_{2} \bullet \sigma_{2} \otimes \cdots \otimes \rho_{k} \bullet \sigma_{k}$ in (\ref{Lproduct1}), and every partition diagram that refines some $\rho_{1} \bullet \sigma_{1} \otimes \rho_{2} \bullet \sigma_{2} \otimes \cdots \otimes \rho_{k} \bullet \sigma_{k}$ in (\ref{Lproduct1}) appears in (\ref{Lproduct2}).
		
		Ignore all empty partitions in $ \rho'_{2},\ldots,\rho'_{k}$ and $ \sigma'_{2},\ldots,\sigma'_{k-1}$, and rewrite $\rho'_{1} \bullet \sigma'_{1} \otimes \rho'_{2} \bullet \sigma'_{2} \otimes \cdots \otimes \rho'_{k} \bullet \sigma'_{k} $ as $\rho''_{1} \bullet \sigma''_{1} \otimes \rho''_{2} \bullet \sigma''_{2} \otimes \cdots \otimes \rho''_{k'} \bullet \sigma''_{k'}$, where $ \rho''_{2},\ldots,\rho''_{k'}$ and $ \sigma''_{1},\ldots,\sigma''_{k'-1}$ are non-empty (possibly $ \otimes $-reducible), and $\rho''_{1} \otimes \cdots \otimes \rho''_{k'}=\rho'$ and $\sigma''_{1} \otimes \cdots \otimes \sigma''_{k'}=\sigma'$. For example, if $\sigma'_{1}=\emptyset$, replace $\rho'_{1} \otimes \rho'_{2}$ by $\rho''_{1}$. 
		
		To get $\rho$ from $\rho'$, we need to change some $ \otimes $s in $\rho'$ to $ \bullet $, each of which may be between $\rho''_{i}$ and $\rho''_{i+1}$  or ``inside'' $\rho''_{j}$, for  some $1\leq i, j\leq k'$, which means there exist $\rho_{1},\ldots,\rho_{k'}$ such that $\rho''_{i}\leq \rho_{i}$ for all  $1\leq i\leq k'$ and $\rho$ can be expressed as $\rho_{1} \divideontimes_{1} \cdots \divideontimes_{k-1} \rho_{k'}$,$ \divideontimes_{i} \in\left\lbrace \otimes , \bullet  \right\rbrace $. (We get $ \rho_{i}$ from $\rho''_{i} $ by replacing necessary $ \otimes $s to $ \bullet $.) 
		
		Similarly, there exist $\sigma_{1},\ldots,\sigma_{k'}$ such that $\sigma''_{i}\leq \sigma_{i}$ for all  $1\leq i\leq k'$ and $\sigma$ can be expressed as $\sigma_{1} \divideontimes_{1} \cdots \divideontimes_{k-1} \sigma_{k'}$,  $ \divideontimes_{i} \in\left\lbrace \otimes , \bullet  \right\rbrace $.
		Then we have that $\rho'_{1} \bullet \sigma'_{1} \otimes \rho'_{2} \bullet \sigma'_{2} \otimes \cdots \otimes \rho'_{k} \bullet \sigma'_{k}=\rho''_{1} \bullet \sigma''_{1} \otimes \rho''_{2} \bullet \sigma''_{2} \otimes \cdots \otimes \rho''_{k'} \bullet \sigma''_{k'} $ refines  $\rho_{1} \bullet \sigma_{1} \otimes \rho_{2} \bullet \sigma_{2} \otimes \cdots \otimes \rho_{k} \bullet \sigma_{k'} $ through Lemma~\ref{repeat}.
		
		Then we prove another direction. Every partition diagram which refines $\rho_{1} \bullet \sigma_{1} \otimes \rho_{2} \bullet \sigma_{2} \otimes \cdots \otimes \rho_{k} \bullet \sigma_{k}$ can be expressed as $\rho''_{1} \divideontimes_{1} \sigma''_{1} \otimes \rho''_{2} \divideontimes_{2}  \sigma''_{2} \otimes \cdots \otimes \rho''_{k} \divideontimes_{k-1}  \sigma''_{k}$, where $ \divideontimes_{i} \in\left\lbrace \otimes , \bullet  \right\rbrace $,$\rho''_{i}\leq\rho_{i}$ and $\sigma''_{i}\leq\sigma_{i}$, for all $1\leq i\leq k$. 
		
		If  $\rho$ can be expressed as $\rho_{1} \divideontimes_{1} \cdots \divideontimes_{k-1} \rho_{k}$, $ \divideontimes_{i} \in\left\lbrace \otimes , \bullet  \right\rbrace $, then $\rho_{1} \otimes \cdots \otimes \rho_{k} $ must refine $\rho$ according to the definition of refinement. The partition diagram $\rho''_{1} \otimes \cdots \otimes \rho''_{k}$ refines $\rho_{1} \otimes \cdots \otimes \rho_{k} $ through Lemma~\ref{repeat}. Therefore, $\rho''_{1} \otimes \cdots \otimes \rho''_{k}$ refines $\rho$ through Lemma~\ref{poset}. Similarly, $\sigma''_{1} \otimes \cdots \otimes \sigma''_{k} $ refines $\sigma$. 
		
		Since $\rho_{1},\ldots,\rho_{k}$ and $\sigma_{1},\ldots,\sigma_{k} $ might be $ \otimes $-reducible, add some $\emptyset$s (if necessary) to rewrite $\rho''_{1} \divideontimes_{1} \sigma''_{1} \otimes \rho''_{2} \divideontimes_{1} \sigma''_{2} \otimes \cdots \otimes \rho''_{k} \divideontimes_{k-1} \sigma''_{k}$ as $\rho'_{1} \bullet \sigma'_{1} \otimes \rho'_{2} \bullet \sigma'_{2} \otimes \cdots \otimes \rho'_{k'} \bullet \sigma'_{k'} $, where $\rho'_{1},\ldots,\rho'_{k'} $ and $\sigma'_{1},\ldots,\sigma'_{k'}$ are $ \otimes $-irreducible, and $ \rho'_{1} \otimes \cdots \otimes \rho'_{k'}=\rho''_{1} \otimes \cdots \otimes \rho''_{k} $, and $\sigma'_{1} \otimes \cdots \otimes \sigma'_{k'}=\sigma''_{1} \otimes \cdots \otimes \sigma''_{k} $. We still have that $ \rho'_{1} \otimes \cdots \otimes \rho'_{k'}$ refines $\rho$, and that $\sigma'_{1} \otimes \cdots \otimes \sigma'_{k'}$ refines $\sigma$.

		$\qedsymbol$
	\end{Proof}
	There is another way to understand this product: For any non-empty partition diagrams $\rho=\rho_{1} \divideontimes^{\rho}_{1} \rho_{2} \divideontimes^{\rho}_{2} \cdots  \divideontimes^{\rho}_{n-1} \rho_{n}$, $\sigma=\sigma_{1} \divideontimes^{\sigma}_{1} \sigma_{2} \divideontimes^{\sigma}_{2} \cdots  \divideontimes^{\sigma}_{m-1} \sigma_{m}$, where $ \divideontimes^{\rho}_{i}, \divideontimes^{\sigma}_{j} \in\left\lbrace \otimes , \bullet  \right\rbrace $, $ \rho_{i}$ and $\sigma_{j}$ are non-empty $  \otimes $-irreducible and $ \bullet $-irreducible partition diagrams for $1\leq i\leq n$, $1\leq j\leq m$,$$L_{\rho}  \star  L_{\sigma}=\mathop{\sum}L_{\pi_{1} \divideontimes_{1} \pi_{2} \divideontimes_{2} \cdots \divideontimes_{n+m-1} \pi_{n+m}} $$ where the sum  over all possible word shuffles $ \pi_{1}\pi_{2}\cdots\pi_{n+m}$ of  $\rho_{1}\rho_{2} \cdots \rho_{n} $ and $ \sigma_{1}\sigma_{2} \cdots \sigma_{m}$. For each $ \pi_{1}\pi_{2}\cdots\pi_{n+m}$, the expression $\pi_{1} \divideontimes_{1} \pi_{2} \divideontimes_{2} \cdots \divideontimes_{n+m-1} \pi_{n+m}$ is unique: \\
	if $\pi_{l}=\rho_{i},\pi_{l} =\sigma_{j}$ for some $1\leq i\leq n,1\leq j\leq m$, then $ \divideontimes_{l} = \bullet $;\\
	if $\pi_{l}=\sigma_{j},\pi_{l} =\rho_{i}$ for some $1\leq i\leq n,1\leq j\leq m$, then $ \divideontimes_{l} = \otimes $; \\if $\pi_{l}=\rho_{i},\pi_{l} =\rho_{i+1}$ for some $1\leq i\leq n$, then $ \divideontimes_{l} = \divideontimes^{\rho}_{i}$;\\
	if $\pi_{l}=\sigma_{j},\pi_{l} =\sigma_{j+1}$ for some $1\leq j\leq m$, then $ \divideontimes_{l} = \divideontimes^{\sigma}_{j}$.

 For any partition diagram $\pi$,
	$$\varepsilon(L_{\pi})=\mathop{\sum}\limits_{\sigma\leq \pi}\varepsilon(M_{\sigma})=\mathop{\sum}\limits_{\sigma\leq \pi}\delta_{\sigma,\emptyset}=\delta_{\pi,\emptyset}.$$
	By Takeuchi's formula, we can get the antipode as follows:  $S(L_{\emptyset})=L_{\emptyset}$. For any non-empty partition diagram $\pi$,
	$$S (L_{\pi})=\mathop{\sum}\limits_{k\geq0}(-1)^{k}(\mathop{\sum}\limits_{r\geq 0}L_{\pi_{1}^{1} \bullet \pi_{2}^{1} \bullet \cdots \bullet \pi_{k}^{1} \otimes \cdots \otimes \pi_{1}^{r} \bullet \pi_{2}^{r} \bullet \cdots \bullet \pi_{k}^{r}}),$$ where $\pi$ can be expressed as $\pi_{1}^{1} \divideontimes_{1}^{1} \pi_{1}^{2} \divideontimes_{1}^{2} \cdots \divideontimes_{1}^{r-1} \pi_{1}^{r} \divideontimes_{1} \pi_{2}^{1} \divideontimes_{2}^{1} \cdots \divideontimes_{2}^{r-1} \pi_{2}^{r} \divideontimes_{2} \cdots \divideontimes_{k-1} \pi_{k}^{1} \divideontimes_{k}^{1} \cdots \divideontimes_{k}^{r-1} \pi_{k}^{r}$, and  $\pi_{1}^{s} \bullet \pi_{2}^{s} \bullet \cdots \bullet \pi_{k}^{s}$ and $\pi_{i}^{1} \otimes \pi_{i}^{2} \otimes \cdots \otimes \pi_{i}^{r}$ are non-empty, $ \divideontimes_{i}^{s}, \divideontimes_{i} \in\left\lbrace  \otimes , \bullet \right\rbrace $, $\max\left\lbrace t|\pi_{t}^{s}\neq\emptyset\right\rbrace >\min\left\lbrace l|\pi_{l}^{s+1}\neq\emptyset\right\rbrace $ for all possible $s,i$.
	
	\section{Other gradings}\label{Other gradings}
	In this section, we introduce further graded algebra structures on both \text{ParSym} and \text{ParQSym}, and compare them with the gradings by orders of partition diagrams.
	\subsection{Grading by \texorpdfstring{$ \otimes $-irreducibility}{irreducibility}}
	Let $P(n):=\left\lbrace \pi| l(\pi)=n\right\rbrace $ for $n\geq 0$. Let $$\text{ParSym}(i)=\text{span}_{\mathbb{K}}\left\lbrace H_{\pi}:\pi\in P(i)\right\rbrace.$$    $$\text{ParSym}=\mathop{\oplus}\limits_{i\geqslant 0}\text{ParSym}(i).$$ 
	This grading can also make \text{ParSym} a graded algebra by the definition of the product of \text{ParSym} (see Definition~\ref{ParSym product}), but cannot make it a graded coalgebra, as shown in the following example.
	\begin{Example}\label{ex1}
		For $ \otimes $-irreducible and $ \bullet $-irreducible $\pi$ and $\rho$, $$\bigtriangleup H_{\pi \bullet \rho}=H_{\emptyset} \otimes  H_{\pi \bullet \rho}+H_{\pi} \otimes  H_{\rho}+H_{\pi \bullet \rho} \otimes  H_{\emptyset},  $$ $H_{\pi \bullet \rho}\in \text{ParSym}(1)$ by Theorem~\ref{irr}, $H_{\emptyset} \otimes  H_{\pi \bullet \rho}+H_{\pi \bullet \rho} \otimes  H_{\emptyset}\in \text{ParSym}(0) \otimes  \text{ParSym}(1)\bigoplus \text{ParSym}(1) \otimes  \text{ParSym}(0)$, but $$H_{\pi} \otimes  H_{\rho}\in \text{ParSym}(1) \otimes  \text{ParSym}(1).$$ Hence  $$\bigtriangleup(\text{ParSym}(1))\nsubseteq  \text{ParSym}(0) \otimes  \text{ParSym}(1)\bigoplus \text{ParSym}(1) \otimes  \text{ParSym}(0).$$
	\end{Example}

	Let $$\text{ParQSym}(i)=\text{span}_{\mathbb{K}}\left\lbrace M_{\pi}:\pi\in P(i)\right\rbrace. $$Then we have $$ \text{ParQSym}=\mathop{\oplus}\limits_{i\geqslant 0}\text{ParQSym}(i).$$
	This grading can also make \text{ParQSym} a graded coalgebra by the definition of the coproduct of \text{ParQSym}~\ref{coproduct}, but cannot make it a graded algebra.
	\begin{Example}
		For $ \otimes $-irreducible $\pi$ and $\rho$, 
		$$M_{\pi} \star  M_{\rho}=M_{\pi \otimes \rho}+M_{\rho \otimes \pi} +M_{\pi \bullet \rho},$$ $M_{\pi \otimes \rho}+M_{\rho \otimes \pi}\in \text{ParQSym}(2)$, but $M_{\pi \bullet \rho}\in \text{ParQSym}(1) $. Therefore $$ \text{ParQSym}(1) \star  \text{ParQSym}(1)\nsubseteq \text{ParQSym}(2) .$$
	\end{Example}
	
	A graded coalgebra $C=\bigoplus C(i)$ is  {\it strictly graded} if $C(0)\cong \mathbb{K}$ and $P(C) = C(1)$, where $ P(C)$ is set of all primitive elements of $C$  \cite{Sweedler}. Directly from Corollary~\ref{primitive}, we have that \text{ParQSym} is strictly graded by the grading $\text{ParQSym}=\mathop{\oplus}\limits_{i\geqslant 0}\text{ParQSym}(i)$. We can also construct gradings for subcoalgebras of \text{ParQSym} to make them strictly graded.
	\begin{Theorem}
		Let $PQ'$ be any subcoalgebra of \text{ParQSym} (for instance, the subcoalgebras mentioned in Section~\ref{Subcoalgebras and Hopf subalgebras}), then $PQ'$ is strictly graded by the grading $$PQ'=\mathop{\oplus}\limits_{i\geqslant 0}(PQ'\cap\text{ParQSym}(i)).$$
	\end{Theorem}
	\begin{Proof}
		Since $PQ'$ is a coalgebra, we have that $$\bigtriangleup PQ'\subseteq PQ' \otimes  PQ'.$$ Notice that $$\bigtriangleup\text{ParQSym}(n)\subseteq \mathop{\oplus}\limits_{i+j=n}\text{ParQSym}(i) \otimes \text{ParQSym}(j).$$Thus $$\bigtriangleup(PQ'\cap\text{ParQSym}(n))\subseteq\mathop{\oplus}\limits_{i+j=n}(PQ'\cap\text{ParQSym}(i)) \otimes (PQ'\cap\text{ParQSym}(j)).$$ So $PQ'=\mathop{\oplus}\limits_{i\geqslant 0}(PQ'\cap\text{ParQSym}(i))$ is a graded coalgebra. In addition,
		$$PQ'\cap\text{ParQSym}(0)=\mathbb{K}\left\lbrace M_{\emptyset} \right\rbrace\cong \mathbb{K},$$
		and $$P(PQ')=PQ'\cap P(\text{ParQSym})=PQ'\cap \text{ParQSym}(1).$$ Therefore $PQ'=\mathop{\oplus}\limits_{i\geqslant 0}(PQ'\cap\text{ParQSym}(i))$ is strictly graded. $\qedsymbol$ 
	\end{Proof}
	
	\begin{Remark}\label{rmk1}
		For $k\geq 1$, let $\pi_{k}\in A_{k}$ be the partition diagram
		
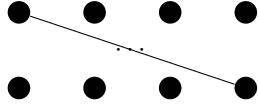
\begin{center}
			\begin{tikzpicture}
				[
				dot/.style={circle, fill=black, inner sep=3pt}, 
				bend angle=30, 
				decoration={brace, mirror, amplitude=10pt} 
				]
				\node[dot] (A) at (0,0) {}; 
				\node[dot] (B) at (1,0) {}; 
				\node[dot] (C) at (2,0) {};
				\node[dot] (D) at (3,0) {};
				\node[dot] (A') at (0,1) {}; 
				\node[dot] (B') at (1,1) {}; 
				\node[dot] (C') at (2,1) {}; 
				\node[dot] (D') at (3,1) {};
				\draw (A') to (D);
				\node at (1.5,0.5) {$\cdots$};
				
			\end{tikzpicture}
                      \end{center}
                      where only one block $ \left\lbrace 1, k'\right\rbrace $ is of size more than one.
                      
		Then $\pi_{k} $ is $ \otimes $-irreducible, that is $$\left\lbrace \pi_{k}|k\geq 1 \right\rbrace \subseteq P(1).$$ Therefore $\text{ParSym}(i)$ and $\text{ParQSym}(i)$ are infinite dimensional, and are not dual to each other, for $i\geq 1$. \text{ParSym} and \text{ParQSym} cannot be seen as (graded) dual by these gradings.
	\end{Remark}
	
	\subsection{Grading by $\otimes$-irreducibility and $\bullet$-irreducibility}
	Let $I_{0}=\left\lbrace \emptyset\right\rbrace $, $$I_{1}=\left\lbrace \pi|\pi\text{ is  non-empty},
	\otimes \text{-irreducible and } \bullet \text{-irreducible}\right\rbrace ,$$ and 
	
	$$I_{n}=\left\lbrace \pi_{1} \divideontimes_{1} \pi_{2} \divideontimes_{2} \cdots \divideontimes_{n-1} \pi_{n}|  \divideontimes_{i} \in\left\lbrace  \otimes ,  \bullet \right\rbrace , \pi_{1},\pi_{2},\ldots,\pi_{n}\in I_{1}\right\rbrace $$ for $n\geq 2$.
	Let $$\text{ParSym}^{(i)}=\text{span}_{\mathbb{K}}\left\lbrace H_{\pi}:\pi\in I_{i}\right\rbrace,$$ and $$\text{ParQSym}^{(i)}=\text{span}_{\mathbb{K}}\left\lbrace M_{\pi}:\pi\in I_{i}\right\rbrace. $$  
	Then $\text{ParSym}=\mathop{\oplus}\limits_{i\geqslant 0}\text{ParSym}^{(i)}$ and $\text{ParQSym}=\mathop{\oplus}\limits_{i\geqslant 0}\text{ParQSym}^{(i) }$.
	By the Definitions~\ref{ParSym product},~\ref{ParSym coproduct},~\ref{coproduct} and~\ref{product}, we have
	$$\bigtriangleup \text{ParSym}^{(n)}\subseteq \bigoplus_{i+j=n}\text{ParSym}^{(i)} \otimes  \text{ParSym}^{(j)},$$ $$\text{ParSym}^{(i)} \text{ParSym}^{(j)}\subseteq \text{ParSym}^{(i+j)},$$  $$\bigtriangleup \text{ParQSym}^{(n)}\subseteq \bigoplus_{i+j=n}\text{ParQSym}^{(i)} \otimes  \text{ParQSym}^{(j)},$$ $$\text{ParQSym}^{(i)} \star  \text{ParQSym}^{(j)}\subseteq \text{ParQSym}^{(i+j)},$$ so $\text{ParSym} $ and $ \text{ParQSym}$ are still graded Hopf algebras by these gradings.
	\begin{Remark}
		$\pi_{k} $ in Remark~\ref{rmk1} is also $ \bullet  $-irreducible for $k\geq 1$, so $\text{ParSym}^{(i)}$ and $\text{ParQSym}^{(i)}$ are infinite dimensional, and are not dual to each other either, for $i\geq 1$. 
		\text{ParSym} and \text{ParQSym} cannot be seen as (graded) dual by these gradings.
	\end{Remark}
	\subsection{Filtration}
	Notice that for any $k\geq0$, $$A_{k}\subseteq \mathop{\cup}\limits_{i\leq k}I_{i},I_{k}\subseteq \mathop{\cup}\limits_{i\leq k}P(i),$$thus$$\text{ParSym}_{k}\subseteq \mathop{\bigoplus}\limits_{i\leq k}\text{ParSym}^{(i)}, \text{ParSym}^{(k)}\subseteq \mathop{\bigoplus}\limits_{i\leq k}\text{ParSym}(i),$$ $$\text{ParQSym}_{k}\subseteq \mathop{\bigoplus}\limits_{i\leq k}\text{ParQSym}^{(i)}, \text{ParQSym}^{(k)}\subseteq \mathop{\bigoplus}\limits_{i\leq k}\text{ParQSym}(i).$$Then we have that $$\mathop{\bigoplus}\limits_{i\leq k}\text{ParSym}_{i} \subseteq \mathop{\bigoplus}\limits_{i\leq k}\text{ParSym}^{(i)} \subseteq \mathop{\bigoplus}\limits_{i\leq k}\text{ParSym}(i), $$ $$\mathop{\bigoplus}\limits_{i\leq k}\text{ParQSym}_{i} \subseteq \mathop{\bigoplus}\limits_{i\leq k}\text{ParQSym}^{(i)} \subseteq \mathop{\bigoplus}\limits_{i\leq k}\text{ParQSym}(i) .$$
	Let $PS_{k}=\mathop{\bigoplus}\limits_{i\leq k}\text{ParSym}_{i}$,$PS^{(k)}=\mathop{\bigoplus}\limits_{i\leq k}\text{ParSym}^{(i)}$,$PS(k)=\mathop{\bigoplus}\limits_{i\leq k}\text{ParSym}(i)$,$PQ_{k}=\mathop{\bigoplus}\limits_{i\leq k}\text{ParQSym}_{i}$,$PQ^{(k)}=\mathop{\bigoplus}\limits_{i\leq k}\text{ParQSym}^{(i)}$,$PQ(k)=\mathop{\bigoplus}\limits_{i\leq k}\text{ParQSym}(i)$, then $$ PS_{k}\subseteq PS^{(k)}\subseteq PS(k),$$ $$ PQ_{k}\subseteq PQ^{(k)}\subseteq PQ(k).$$\\

	The properties of gradings can affect corresponding filtrations.
	\begin{Lemma}
		Let $H=\mathop{\bigoplus}\limits_{n\geq 0}H_{n}$ be a Hopf algebra, $S_{n}=\mathop{\bigoplus}\limits_{i\leq n}H_{i}$, then\\
		(a) if $H=\mathop{\bigoplus}\limits_{n\geq 0}H_{n}$ is a graded coalgebra, then $\left\lbrace S_{n} \right\rbrace $ is a coalgebra filtration;\\
		(b) if $H=\mathop{\bigoplus}\limits_{n\geq 0}H_{n}$ is a graded algebra, then $\left\lbrace S_{n} \right\rbrace $ is an algebra filtration;\\
		(c) if $H=\mathop{\bigoplus}\limits_{n\geq 0}H_{n}$ is a graded Hopf algebra, then $\left\lbrace S_{n} \right\rbrace $ is a Hopf algebra filtration.
	\end{Lemma}
	\begin{Proof}
		By the definition of $S_{n}$, we can see that $$S_{0}\subseteq S_{1}\subseteq\cdots\subseteq\mathop{\cup}\limits_{n\geq 0}S_{n}=H.$$
		(a) if $H=\mathop{\bigoplus}\limits_{n\geq 0}H_{n}$ is a graded coalgebra, then 
		$$\bigtriangleup S_{n}=\mathop{\Sigma}\limits_{i\leq n}\bigtriangleup H_{i}
		\subseteq \mathop{\Sigma}\limits_{i\leq n,k+j=i}H_{k} \otimes  H_{j}
		=\mathop{\Sigma}\limits_{k\leq n,j\leq n-k}H_{k} \otimes  H_{j}
		\subseteq \mathop{\Sigma}\limits_{k\leq n}S_{k} \otimes  S_{n-k}.$$
		(b) if $H=\mathop{\bigoplus}\limits_{n\geq 0}H_{n}$ is a graded algebra, then $$S_{n}S_{m}=\mathop{\Sigma}\limits_{i\leq n,j\leq m}H_{i}H_{j}\subseteq \mathop{\Sigma}\limits_{i\leq n,j\leq m}H_{i+j}\subseteq S_{n+m}.$$
		(c) if $H=\mathop{\bigoplus}\limits_{n\geq 0}H_{n}$ is a graded Hopf algebra, then $S(H_{n})\subseteq H_{n}$. Therefore
		$$S(S_{n})=\mathop{\Sigma}\limits_{i\leq n}S(H_{i})\subseteq \mathop{\Sigma}\limits_{i\leq n}H_{i}=S_{n}.$$Combining (a) and (b), we have that $\left\lbrace S_{n} \right\rbrace $ is a Hopf algebra filtration.$\qedsymbol$
	\end{Proof}
	Using the above lemma, we can find the properties of the six filtrations defined in this section.
	\begin{Proposition}
		Considering the six filtrations defined in this section, we have the following results:\\
		(a)$\left\lbrace PS_{k}\right\rbrace $ and $\left\lbrace PS^{(k)} \right\rbrace $ are Hopf filtrations of \text{ParSym};\\
		(b)$\left\lbrace PQ_{k}\right\rbrace $, $\left\lbrace PQ(k)\right\rbrace $ and $\left\lbrace PQ^{(k)} \right\rbrace $ are Hopf filtrations of \text{ParQSym};\\ (c)$\left\lbrace PS(k)\right\rbrace $ is an algebra filtration of \text{ParSym} but not a coalgebra filtration;\\
		(d)$\left\lbrace PQ(k)\right\rbrace $ is the coradical filtration of \text{ParQSym}.
	\end{Proposition}
	\begin{Proof}
		From the last lemma and the analysis of the gradings, we have \\ (1)$\left\lbrace PS_{k}\right\rbrace $ and $\left\lbrace PS^{(k)} \right\rbrace $ are Hopf filtrations of \text{ParSym} (The Hopf filtration property given in (a) is thus proved);\\
		(2) $\left\lbrace PQ_{k}\right\rbrace $ and $\left\lbrace PQ^{(k)} \right\rbrace $ are Hopf filtrations of \text{ParQSym};\\ (3)$\left\lbrace PS(k)\right\rbrace $ is an algebra filtration of \text{ParSym};\\(4) $\left\lbrace PQ(k)\right\rbrace $ is a coalgebra filtration of \text{ParQSym}. \\Notice that the operation $ \bullet $ can only reduce the length, so for any partition diagrams $\pi$ and $\rho$, $$M_{\pi} \star  M_{\rho}\in PQ(l(\pi)+l(\rho)), S(M_{\pi})\in PQ(l(\pi)).$$  Then $$PQ(n) \star  PQ(m)\subseteq PQ(n+m),S(PQ(n))\subseteq PQ(n)$$for all $n,m\geq 0$. Thus $ \left\lbrace PQ(k)\right\rbrace $ is also a Hopf filtration of \text{ParQSym}. That finishes the proof of the Hopf filtration property given in (b).\\
		Consider Example~\ref{ex1}, let $\pi,\rho\in P(1)$, then $H_{\pi \bullet \rho}\in PS(1)$. However, $$H_{\pi} \otimes  H_{\rho}\notin PS(0) \otimes  PS(1)+PS(1) \otimes  PS(0).$$ So $\bigtriangleup PS(1)\varsubsetneq PS(0) \otimes  PS(1)+PS(1) \otimes  PS(0)$, this is the Hopf filtration property given in (c). \\
		We proceed to prove the filtration property in (d).
		According to Lemma~\ref{connected}, we have that \text{ParQSym} is connected, that is, the coradical of \text{ParQSym} is $\mathbb{K}=PQ(0)$. Then we will prove $PQ(n)=\bigtriangleup^{-1}(\text{ParQSym} \otimes  PQ(n-1)+PQ(0) \otimes \text{ParQSym})$. \\For any partition diagram $\pi=\pi_{1} \otimes \pi_{2}$, if $l(\pi)\leq n$, then $l(\pi_{1})=0$ or $l(\pi_{2})\leq n-1$, which means $\bigtriangleup PQ(n)\subseteq \text{ParQSym} \otimes  PQ(n-1)+PQ(0) \otimes \text{ParQSym}$.\\ If $l(\pi)>n$, then there exist a pair of partition diagrams $\pi_{1}$ and $\pi_{2}$ such that $\pi=\pi_{1} \otimes \pi_{2}$, $l(\pi_{1})>0$ and  $l(\pi_{2})> n-1$, which means $$\bigtriangleup M_{\pi}\notin \text{ParQSym} \otimes  PQ(n-1)+PQ(0) \otimes \text{ParQSym}.$$ Therefore $$ PQ(n)=\bigtriangleup^{-1}(\text{ParQSym} \otimes  PQ(n-1)+PQ(0) \otimes \text{ParQSym}).$$
		$\qedsymbol$
	\end{Proof}
	It is easy to see that there are some maps of filtered coalgebras between these filtered coalgebras.
	\begin{Proposition}
		The forgetful maps $$(\text{ParQSym},\left\lbrace PQ_{k}\right\rbrace)\rightarrow (\text{ParQSym},\left\lbrace PQ^{(k)}\right\rbrace),$$ $$(\text{ParQSym},\left\lbrace PQ_{k}\right\rbrace)\rightarrow (\text{ParQSym},\left\lbrace PQ(k)\right\rbrace),$$ $$(\text{ParQSym},\left\lbrace PQ^{(k)}\right\rbrace)\rightarrow (\text{ParQSym},\left\lbrace PQ(k)\right\rbrace),$$ and $$(\text{ParSym},\left\lbrace PS_{k}\right\rbrace)\rightarrow (\text{ParSym},\left\lbrace PS^{(k)}\right\rbrace)$$  are maps of filtered coalgebras.
	\end{Proposition}
	\begin{Proof}
		This follows directly from the definition of maps of filtered coalgebras and that $$ PS_{k}\subseteq PS^{(k)}\subseteq PS(k), PQ_{k}\subseteq PQ^{(k)}.$$ $\qedsymbol$
	\end{Proof}
	Using the following lemma, we can build more coalgebra filtrations.
	\begin{Lemma}[\protect{\cite[Lemma 4.1.3]{David} }]
		Let $(C,\left\lbrace V_{n} \right\rbrace )$ be a filtered coalgebra. If $D$ is a subcoalgebra of $C$, then $\left\lbrace V_{n}\cap D \right\rbrace $ is a coalgebra filtration of $D$.
	\end{Lemma}
	Applying the above lemma to the subcoalgebras of \text{ParQSym} mentioned in this article and the subcoalgebras of \text{ParSym} mentioned in~\cite{ref1}, we can get respective coalgebra filtrations.
	\section{Infinitorial Hopf algebra}\label{Infinitorial Hopf algebra}
	In this section, we will introduce the definition of  {\it Infinitorial Hopf algebra}, and define an  {\it infinitesimal character} $\eta_{\text{ParQSym}} $ to make $(\text{ParQSym},\eta_{\text{ParQSym}})$ an infinitorial Hopf algebra.
	For a composition $\alpha=(\alpha_{1},\alpha_{2},\ldots,\alpha_{n})$, define $lp(\alpha)=\alpha_{n}$  \cite[Definition 2.1]{ref8}. By convention, $lp(()) = 0$.
	Let $H$ be a Hopf algebra. An  {\it infinitesimal character}  of $H$ is a linear map
	$\xi:H\rightarrow\mathbb{K}$ that satisfies $$\xi(ab)=\varepsilon(a)\xi(b)+\xi(a)\varepsilon(b)$$for all $a,b\in H$  \cite{refDominique Manchon}. The set of all infinitesimal characters of $H$ is denoted $\Xi(H)$.
	There is an easy way to determine whether a linear map is an infinitesimal character, using the following equivalent condition.
	\begin{Proposition}
		Let $H$ be a connected graded Hopf algebra. A linear map $\xi:H\rightarrow\mathbb{K}$ is an infinitesimal character of $H$ if and only if $\xi(1_{H}) = 0$ and $\xi(ab) = 0$ for all homogeneous
		$a, b \in H$ of positive degree  \cite{ref8}.
	\end{Proposition}
	Now we define infinitorial Hopf algebras.
	\begin{Definition}\label{iff}\cite[Definition 3.7]{ref8}
		An  {\it infinitorial Hopf algebra} is a pair $(H, \xi)$, where $H$ is a connected graded
		Hopf algebra and $\xi \in \Xi(H)$. A morphism of infinitorial Hopf algebras $(H, \xi)\rightarrow (H', \xi')$ is a
		graded Hopf map $\Phi: H \rightarrow H'$
		that also satisfies $\xi = \xi'\circ \Phi$.
	\end{Definition}
	Define a linear map $\eta:\text{QSym}\rightarrow\mathbb{K}$ by $$\eta(M_{\alpha})=(-1)^{l(\alpha)-1} lp(\alpha).$$Liu and Weselcouch~\cite[Section 3.1]{Ricky Ini Liu and Michael Weselcouch} showed that $\eta$ is an infinitesimal
	character of \text{QSym}. Now we define a linear map $\eta_{\text{ParQSym}}:\text{ParQSym}\rightarrow\mathbb{K}$ by $$\eta_{\text{ParQSym}}=\eta\circ\Psi_{PQ}, $$where $\Psi_{PQ}$ is defined as (\ref{quotient map}).
	\begin{Proposition}
		$(\text{ParQSym},\eta_{\text{ParQSym}})$ is an infinitorial Hopf algebra.
	\end{Proposition}
	\begin{Proof}
		Consider the empty diagram:
		$$ 	\eta_{\text{ParQSym}}(M_{\emptyset})
		=\eta\circ\Psi_{PQ}(M_{\emptyset})
		=\eta(M_{()})
		=0.$$
		The first and the second equations are the definitions of $\eta_{\text{ParQSym}} $ and $\Psi_{PQ}$. The last equation holds since $\eta$ is an infinitesimal character.
		For any non-empty partition diagrams $\pi$ and $\rho$, 
		\begin{equation*}
			\begin{aligned}
				&\eta_{\text{ParQSym}}(M_{\pi} \star  M_{\rho})\\
				=&\eta\circ\Psi_{PQ}(M_{\pi} \star  M_{\rho})\\
				=&\eta(\Psi_{PQ}(M_{\pi})\Psi_{PQ}(M_{\rho}))\\
				=&\eta(M_{\alpha_{\pi}}M_{\alpha_{\rho}})\\
				=&0.
			\end{aligned}
		\end{equation*}
		The first and the third equations are the definitions of $\eta_{\text{ParQSym}} $ and $\Psi_{PQ}$. The second equation holds because $\Psi_{PQ}$ is an algebra map (since it is a map of Combinatorial Hopf Algebras). The last equation holds since $\eta$ is an infinitesimal character.\\
		Every homogeneous
		$a \in \text{ParQSym}$ of positive degree can be seen as a finite linear sum of $M_{\pi}$ with non-empty $\pi$. Therefore, the condition in~\ref{iff} holds.
		
		$\qedsymbol$
	\end{Proof}
	According to  \cite[Section 3]{ref8}, the category of infinitorial Hopf algebras has a terminal object $(Sh,\xi_{s})$, which is defined as follows.
	\begin{Definition} 
		The shuffle algebra $Sh$ is the connected graded Hopf algebra generated from a basis
		${x_{\alpha}}$ indexed by compositions $\alpha$, where the grading is given by $Sh_{n} = \text{span}\left\lbrace x_{\alpha} : \alpha \vDash n\right\rbrace $, whose product is given by shuffling:$$x_{\alpha}x_{\beta} =\mathop{\sum}\limits_{\gamma\in\alpha\shuffle\beta}x_{\gamma},$$and whose coproduct is given by deconcatenation:$$\Delta(x_{\gamma}) =\mathop{\sum}\limits_{\alpha\cdot\beta=\gamma}x_{\alpha}  \otimes  x_{\beta}.$$
		$\xi_{s}:Sh\rightarrow\mathbb{K} $ is defined by 
		\begin{equation*}
			\xi_{s}(x_{\alpha})=\left\{\begin{aligned}
				&1, \text{if } l(\alpha)=1,\\ 
				&0, \text{otherwise}.
			\end{aligned}
			\right.
		\end{equation*}
	\end{Definition}
	Similar to the category of Combinatorial Hopf Algebras, there is a formula for the unique morphism of infinitorial Hopf algebras $(H, \xi)\rightarrow (Sh,\xi_{s})$ as follows:
	\begin{Proposition}\cite[Theorem 3.10]{ref8} 
		For any infinitorial Hopf algebra $(\mathcal{H},\xi)$, there exists a unique morphism of infinitorial Hopf algebras
		$$\Phi:(\mathcal{H},\xi)\rightarrow(Sh,\xi_{s}).$$
		Moreover, $\Phi$ is explicitly given as follows. For $h\in \mathcal{H}_{n}$,$$\Phi(h)=\mathop{\sum}\limits_{\alpha\vDash n}\xi_{\alpha}(h)x_{\alpha},$$where for $\alpha=(\alpha_{1},\ldots,\alpha_{k})$, $\xi_{\alpha} $ is the composite$$\mathcal{H}\stackrel{\bigtriangleup^{(k-1)}}{\longrightarrow}\mathcal{H}^{ \otimes  k}\twoheadrightarrow \mathcal{H}_{\alpha_{1}} \otimes \cdots \otimes  \mathcal{H}_{\alpha_{k}}\stackrel{\xi^{ \otimes  k}}{\longrightarrow}\mathbb{K}, $$where the unlabeled map is the tensor product of the canonical projections onto the
		homogeneous components $\mathcal{H}_{\alpha_{i}} $ for $1\leq i\leq k$.
	\end{Proposition}
	Using the above proposition, we can express the unique morphism of infinitorial Hopf algebras
	$$\Phi_{PS}:(\text{ParQSym},\eta_{\text{ParQSym}})\rightarrow(Sh,\xi_{s})$$  according to the formula$$ \Phi_{PS}(M_{\rho})=\mathop{\sum}\limits_{\alpha_{\rho}\leq\beta }(-1)^{l(\alpha_{\rho})-l(\beta)}\mathop{\prod}\limits_{1\leq i\leq l(\beta)}lp(\alpha_{\rho}^{(i)})x_{\beta},$$ where $(\alpha_{\rho})_{\beta}^{(i)} $ is defined as (\ref{quotient map}) and (\ref{refine part}).
	\section{Deconcatenation basis}\label{Deconcatenation basis}
	In this section, we will define the deconcatenation bases of \text{ParQSym} and give a way to construct deconcatenation bases.
	From the definition of refinement of partition diagrams, we can directly get the following equivalent condition.
	\begin{Proposition}
		For partition diagrams $\pi $ and $\rho=\rho_{1} \otimes \rho_{2} \otimes \cdots \otimes \rho_{l(\rho)}$, where $\rho_{i}$ are $ \otimes $-irreducible, then $\pi\leq\rho$ if and only if $\pi$ can be written as
		\begin{equation}\label{refine condition}
			\pi=\pi_{\rho}^{(1)} \otimes \pi_{\rho}^{(2)} \otimes \cdots \otimes \pi_{\rho}^{(l(\rho))}, 
		\end{equation}where $\pi_{\rho}^{(i)}\leq\rho_{i} $ for each $i$.
	\end{Proposition}
	Let $A:=\mathop{\cup}\limits_{i\geq 0 }A_{i}$ be the set of all partition diagrams, and $A_{>0}:=\mathop{\cup}\limits_{i> 0 }A_{i}$ be the set of all non-empty partition diagrams.
	Given a function on partition diagrams, we can extend it to a function on a pair of partition diagrams $\pi$ and $\rho$, where $\pi\leq\rho$, in the following way, similar to  \cite[Definition 2.2]{Ricky Ini Liu and Michael Weselcouch}.
	\begin{Definition}
		Let $f : A_{>0}\rightarrow \mathbb{K}$ be a function. For partition diagrams $\pi$ and $\rho$, where $\pi\leq\rho$, we define$$f(\pi,\rho)=f(\pi_{\rho}^{(1)})f(\pi_{\rho}^{(2)})\cdots f(\pi_{\rho}^{(l(\rho))}),$$ where $\pi_{\rho}^{(1)} $,$\cdots$,$\pi_{\rho}^{(l(\rho))}$ are given by (\ref{quotient map}). (By convention, $f(\emptyset, \emptyset) = 1$.)
	\end{Definition}
	\begin{Corollary}\label{deconcatation}
		Let $f : A_{>0}\rightarrow \mathbb{K}$ be a function. For partition diagrams $\pi$ and $\rho$, where $\pi\leq\rho$, if $\rho=\rho_{1} \otimes  \rho_{2}$, then there exist $\pi_{1}$ and $\pi_{2} $ such that $ \pi=\pi_{1} \otimes \pi_{2}$, $\pi_{i} \leq\rho_{i}$, and $$f(\pi,\rho)=f(\pi_{1},\rho_{1})f(\pi_{2},\rho_{2}).$$
	\end{Corollary}
	\begin{Proof}
		If $$\rho=\rho_{1}' \otimes \rho_{2}' \otimes \cdots \otimes \rho_{l(\rho)}',$$ where $\rho_{i}'$ is non-empty $ \otimes $-irreducible for each $i$, then there exists $i$ such that $$\rho_{1}= \rho_{1}' \otimes \rho_{2}' \otimes \cdots\rho_{i}'$$ and $$\rho_{2}=\rho_{i+1}' \otimes \rho_{i+2}' \otimes \cdots\rho_{l(\rho)}'.$$ Let $$\pi_{1}=\pi_{\rho}^{(1)} \otimes \pi_{\rho}^{(2)} \otimes \cdots \otimes \pi_{\rho}^{(i)}$$ and $$\pi_{2}=\pi_{\rho}^{(i+1)} \otimes \pi_{\rho}^{(i+2)} \otimes \cdots \otimes \pi_{\rho}^{(l(\rho))},$$ where $\pi_{\rho}^{(1)} $,$\cdots$,$\pi_{\rho}^{(l(\rho))}$ are given by (\ref{quotient map}), which means for each $1\leq j\leq l(\rho)$, $\pi_{\rho}^{(j)}\leq\rho_{j}$. 
		Then $\pi_{i} \leq\rho_{i}$, and 
		$$f(\pi,\rho)=f(\pi_{\rho}^{(1)})\cdots f(\pi_{\rho}^{(i)})f(\pi_{\rho}^{(i+1)})\cdots f(\pi_{\rho}^{(l(\rho))})=f(\pi_{1},\rho_{1})f(\pi_{2},\rho_{2}).$$
		$\qedsymbol$
	\end{Proof}
	Similar to  \cite[Definition 4.1 and Definition 4.2]{Ricky Ini Liu and Michael Weselcouch}, we define deconcatenation bases 
	of \text{ParQSym}.
	\begin{Definition}
		A  {\it deconcatenation basis} of \text{ParQSym} is a graded basis $\left\lbrace X_{\pi} \right\rbrace $ that satisfies $$\Delta(X_{\sigma}) =\mathop{\sum}\limits_{\pi \otimes \rho=\sigma}X_{\pi}  \otimes  X_{\rho},$$for all partition diagrams $\sigma$.
	\end{Definition}
	We also need to define nonsingular function $f : A_{>0}\rightarrow \mathbb{K}$:
	\begin{Definition}
		We call a function $f : A_{>0}\rightarrow \mathbb{K}$  {\it nonsingular} if $f(\pi) \neq 0$ for all non-empty $ \otimes $-irreducible $\pi$.
	\end{Definition}
	Similar to  \cite[Proposition 4.4.]{Ricky Ini Liu and Michael Weselcouch}, we have the following proposition related to nonsingular functions.
	\begin{Proposition}
		Let $\left\lbrace X_{\pi} \right\rbrace $ and $\left\lbrace P_{\pi} \right\rbrace $ be graded bases of $\text{ParQSym}$. Then the following are
		equivalent:\\
		(i) There exists a nonsingular function $f : A_{>0}\rightarrow \mathbb{K}$ such that $$Q_{\pi}= \mathop{\sum}\limits_{\pi\leq\rho}f(\pi,\rho)P_{\rho};$$\\
		(ii) There exists a nonsingular function $g : A_{>0}\rightarrow \mathbb{K}$ such that
		$$P_{\pi}= \mathop{\sum}\limits_{\pi\leq\rho}g(\pi,\rho)Q_{\rho}.$$
	\end{Proposition}
	\begin{Proof}
		By symmetry, it suffices to prove only one direction, so assume that (i) holds. To
		define the function $g$, consider the equations
		\begin{equation}\label{transform}
			\mathop{\sum}\limits_{\pi\leq\rho}f(\pi,\rho)g(\rho)=\begin{cases}
				1 & \text{if }l(\pi)=1,\\ 
				0 &\text{otherwise},
			\end{cases}
		\end{equation}
		for all $\pi\in A_{>0}$. By nonsingularity, $f(\pi, \pi)\neq 0$ for each $\pi$, so it follows from triangularity that (\ref{transform}) uniquely determines $g : A_{>0}\rightarrow \mathbb{K}$. Additionally, when $\pi$ is non-empty $ \otimes $- irreducible, we have
		$f(\pi)g(\pi) = 1$, so $g(\pi) \neq0$ and $g$ is also nonsingular. To show that $g$ satisfies (ii), define $P'_{\pi}= \mathop{\sum}\limits_{\pi\leq\rho}g(\pi,\rho)Q_{\rho}$. Then for fixed $\pi\in A_{>}$,
		\begin{equation*}
			\begin{aligned}
				\mathop{\sum}\limits_{\pi\leq\rho}f(\pi,\rho)P'_{\pi}=&\mathop{\sum}\limits_{\pi\leq\rho}f(\pi,\rho)\mathop{\sum}\limits_{\rho\leq\sigma}g(\rho,\sigma)Q_{\sigma}\\
				=&\mathop{\sum}\limits_{\pi\leq\sigma}(\mathop{\sum}\limits_{\rho:\pi\leq\rho\leq\sigma}f(\pi,\rho)g(\rho,\sigma))Q_{\sigma}.\\
			\end{aligned}
		\end{equation*}
		For each $\sigma=\sigma_{1} \otimes \sigma_{2} \otimes \cdots \otimes \sigma_{l(\sigma)}$, where $\sigma_{i}$ is $ \otimes $-irreducible, $\pi$ can be written as $\pi^{(1)} \otimes \cdots \otimes \pi^{(l(\sigma))} $, where $\pi^{(i)} \leq\sigma_{i}$. Then $\pi\leq\rho\leq\sigma $ if and only if $\rho$ be written as $\rho^{(1)} \otimes \cdots \otimes \rho^{(l(\sigma))} $, where $\pi^{(i)} \leq\rho^{(i)}$. Hence, the sum in parentheses factors as$$ \mathop{\sum}\limits_{\rho:\pi\leq\rho\leq\sigma}f(\pi,\rho)g(\rho,\sigma)=\mathop{\prod}\limits_{1\leq i\leq l(\sigma)}\mathop{\sum}\limits_{\pi^{(i)} \leq\rho^{(i)}}f(\pi^{(i)},\rho^{(i)})g(\rho^{(i)}).$$By (\ref{transform}), this product vanishes unless $\pi^{(i)}=\sigma_{i}$ for each $i$, or equivalently $\pi=\sigma$, in which case it equals 1. Therefore, we obtain$$	\mathop{\sum}\limits_{\pi\leq\rho}f(\pi,\rho)P'_{\rho}=Q_{\pi}. $$
		Combined with (i), this forces $P'_{\pi}=P_{\pi}$ (again by triangularity), completing the proof. 	$\qedsymbol$
	\end{Proof}
	The following proposition gives a way to construct further deconcatenation bases of \text{ParQSym} from a known deconcatenation basis.
	\begin{Proposition}\label{7.3}
		Let $\left\lbrace P_{\pi} \right\rbrace $ be a deconcatenation basis of \text{ParQSym}. Choose elements $\left\lbrace Q_{\pi} \right\rbrace $ of \text{ParQSym} for each $\pi\in A$. Then $\left\lbrace Q_{\pi} \right\rbrace $ is a deconcatenation basis if 
		there exists
		a nonsingular function $f : A_{>0}\rightarrow \mathbb{K}$ such that
		\begin{equation}\label{relation}
			Q_{\pi}= \mathop{\sum}\limits_{\pi\leq\rho}f(\pi,\rho)P_{\rho}.
		\end{equation}
	\end{Proposition}
	\begin{Proof}
		If (\ref{relation}) holds, the above proposition tells us that there exists a nonsingular function $g: A_{>0}\rightarrow \mathbb{K}$ such that $P_{\pi}= \mathop{\sum}\limits_{\pi\leq\rho}g(\pi,\rho)Q_{\rho}$. Hence $\left\lbrace Q_{\pi} \right\rbrace $ is a basis.
		\begin{equation*}
			\begin{aligned}
				\bigtriangleup(Q_{\pi})&=\mathop{\sum}\limits_{\pi\leq\rho}f(\pi,\rho)\bigtriangleup(P_{\rho})\\
				&=\mathop{\sum}\limits_{\pi\leq\rho}f(\pi,\rho)\mathop{\sum}\limits_{\rho=\rho_{1} \otimes  \rho_{2}}P_{\rho_{1}} \otimes  P_{\rho_{2}}\\
				&=\mathop{\sum}\limits_{\substack{\pi_{1}\leq\rho_{1}\\\pi_{2}\leq\rho_{2}\\\pi=\pi_{1} \otimes \pi_{2}}}f(\pi_{1},\rho_{1})P_{\rho_{1}} \otimes  f(\pi_{2},\rho_{2})P_{\rho_{2}}\\
				&=\mathop{\sum}\limits_{\pi=\pi_{1} \otimes \pi_{2}}Q_{\pi_{1}} \otimes  Q_{\pi_{2}}.
			\end{aligned}
		\end{equation*}
		The first and the last equations use (\ref{relation}), and the second one holds because $\left\lbrace P_{\pi} \right\rbrace $ is a deconcatenation basis of \text{ParQSym}. The third one uses Corollary~\ref{deconcatation}.\\
		$\qedsymbol$
	\end{Proof}
	By definition, we have that $\left\lbrace M_{\pi} \right\rbrace $ is a deconcatenation basis, then we can use the above proposition to construct other deconcatenation bases.
	\section{The enriched monomial basis}\label{The enriched monomial basis}
	In this section, we assume 2 is invertible in $\mathbb{K}$.           
	There is another $\mathbb{K}$-basis of \text{QSym} called the  {\it enriched monomial functions} $\left\lbrace \eta_{\alpha}\right\rbrace $  \cite[Definition 6]{enriched monomial basis}:$$\eta_{\alpha}=\mathop{\sum}\limits_{\alpha\leq \beta}2^{l(\beta)}M_{\beta}, $$
	for any composition $\alpha$. 
	The coproduct of this basis is given by $$ \bigtriangleup \eta_{\alpha}=\mathop{\sum}\limits_{ \beta\odot\gamma= \alpha}\eta_{\beta} \otimes  \eta_{\gamma},$$see \cite[Proposition 2]{enriched monomial basis}.
	By analogy with $\left\lbrace \eta_{\alpha}\right\rbrace $, we define the following basis of $\text{ParQSym} $:
	Let $$\eta_{\pi}:=\mathop{\sum}\limits_{\pi\leq \sigma}2^{l(\sigma)}M_{\sigma}, $$then  $\left\lbrace \eta_{\pi}\right\rbrace_{\pi\in A}$ forms a deconcatenation basis of $\text{ParQSym} $  according to the Proposition~\ref{7.3}, (let $f(\pi):=2 $ for all $\pi$), that is $$\bigtriangleup\eta_{\pi}=\mathop{\sum}\limits_{\rho \otimes \sigma=\pi}\eta_{\rho} \otimes \eta_{\sigma},$$for all $\pi$. In this case, applying (\ref{transform}), we have that $$M_{\pi}:=2^{-l(\pi)}\mathop{\sum}\limits_{\pi\leq \sigma}(-1)^{l(\pi)-l(\sigma)}\eta_{\sigma}. $$
		\begin{Lemma}\cite[Lemma 1]{enriched monomial basis}
		Let $S$ and $T$ be two finite sets. Then,\begin{equation*}
			\mathop{\sum}\limits_{I\subseteq S}(-1)^{|I\setminus T|}=\begin{cases}
				2^{|S|} &\text{if }S\subseteq T,\\
				0 &\text{otherwise}.
			\end{cases}
		\end{equation*} 
	\end{Lemma}
	With the above lemma, we can find the relation between the $\eta$-basis and the $L$-basis of \text{ParQSym}.
	\begin{Proposition} The relation between the $\eta$-basis and the $L$-basis of \text{ParQSym} is as follows:\\
		(1) For any partition diagram $\pi$,$$\eta_{\pi}=2\mathop{\sum}\limits_{\pi\sim \rho}(-1)^{|S(\rho)\setminus S(\pi)|}L_{\rho}.$$
		(2) For any $\pi\in I_{n}$,
		$$L_{\pi}=2^{-n}\mathop{\sum}\limits_{\pi\sim \rho}(-1)^{|S(\pi)\setminus S(\rho)|}\eta_{\rho}.$$
	\end{Proposition}
	\begin{Proof}
		(1)	For any partition diagram $\pi$,
		\begin{equation*}
			\begin{aligned}
				2\mathop{\sum}\limits_{\pi\sim \rho}(-1)^{|S(\rho)\setminus S(\pi)|}L_{\rho}&=2\mathop{\sum}\limits_{\rho:\pi\sim \rho}(-1)^{|S(\rho)\setminus S(\pi)|}\mathop{\sum}\limits_{\sigma:\sigma\leq \rho}M_{\sigma}\\
				&=2\mathop{\sum}\limits_{\sigma:\pi\sim \sigma}(\mathop{\sum}\limits_{\rho:\sigma\sim\rho, S(\rho)\subseteq S(\sigma)}(-1)^{|S(\rho)\setminus S(\pi)|})M_{\sigma}.
			\end{aligned}
		\end{equation*}
		Notice that $\mathop{\sum}\limits_{\rho:\sigma\sim\rho, S(\rho)\subseteq S(\sigma)}(-1)^{|S(\rho)\setminus S(\pi)|}$ is nonzero only when $S(\sigma)\subseteq S(\pi)$. In this case, $\pi\leq\sigma$, and $$\mathop{\sum}\limits_{\rho:\sigma\sim\rho, S(\rho)\subseteq S(\sigma)}(-1)^{|S(\rho)\setminus S(\pi)|}=2^{|S(\sigma)|}.$$
		Therefore 	\begin{equation*}
			\begin{aligned}
				2\mathop{\sum}\limits_{\pi\sim \rho}(-1)^{|S(\rho)\setminus S(\pi)|}L_{\rho}
				&=2\mathop{\sum}\limits_{\pi\leq \sigma}2^{|S(\sigma)|}M_{\sigma}\\
				&=\mathop{\sum}\limits_{\pi\leq \sigma}2^{l(\sigma)}M_{\sigma}\\
				&=\eta_{\pi}.
			\end{aligned}
		\end{equation*}
		(2) For any $\pi\in I_{n}$,
		\begin{equation*}
			\begin{aligned}
				&2^{-n}\mathop{\sum}\limits_{\pi\sim \rho}(-1)^{|S(\pi)\setminus S(\rho)|}\eta_{\rho}\\&=	2^{-n}\mathop{\sum}\limits_{\pi\sim \rho}(-1)^{|S(\pi)\setminus S(\rho)|}\mathop{\sum}\limits_{\sigma:\sigma\geq \rho}2^{l(\sigma)}M_{\sigma}\\
				&=2^{-n}\mathop{\sum}\limits_{\sigma:\pi\sim \sigma}(\mathop{\sum}\limits_{\rho:\sigma\geq \rho}(-1)^{|S(\pi)\setminus S(\rho)|})2^{l(\sigma)}M_{\sigma}.\\
			\end{aligned}
		\end{equation*}
		Notice that \begin{equation*}
			\begin{aligned}
				&\mathop{\sum}\limits_{\rho:\sigma\geq \rho}(-1)^{|S(\pi)\setminus S(\rho)|}\\&=\mathop{\sum}\limits_{ ([n-1]\setminus S(\rho))\subseteq ([n-1]\setminus S(\sigma))}(-1)^{|([n-1]\setminus S(\rho))\setminus ([n-1]\setminus S(\pi))|}
			\end{aligned}
		\end{equation*} is nonzero only when $([n-1]\setminus S(\sigma))\subseteq ([n-1]\setminus S(\pi))$that is,$S(\sigma)\supseteq S(\pi)$. In this case, $\pi\geq\sigma$, and $$\mathop{\sum}\limits_{\rho:\sigma\geq \rho}(-1)^{|S(\pi)\setminus S(\rho)|}=2^{|[n-1]\setminus S(\sigma)|}.$$
		Therefore 	\begin{equation*}
			\begin{aligned}
				&2^{-n}\mathop{\sum}\limits_{\pi\sim \rho}(-1)^{|S(\pi)\setminus S(\rho)|}\eta_{\rho}\\
				&=2^{-n}\mathop{\sum}\limits_{\pi\geq \sigma}2^{|[n-1]\setminus S(\sigma)|}2^{l(\sigma)}M_{\sigma}\\
				&=\mathop{\sum}\limits_{\pi\geq \sigma}M_{\sigma}\\
				&=L_{\pi}.
			\end{aligned}
		\end{equation*}
		$\qedsymbol$
	\end{Proof}
	From the definition, we can find $\eta_{\emptyset}=M_{\emptyset}$ is the unit. The next proposition is to tell the product of the $\eta$-basis for non-empty partition diagrams, which is similar to  \cite[Theorem 5]{enriched monomial basis}.
	\begin{Proposition}
		For any non-empty partition diagrams $\rho=\rho_{1} \otimes \rho_{2} \otimes  \cdots  \otimes \rho_{n}$, $\sigma=\sigma_{1} \otimes \sigma_{2} \otimes  \cdots  \otimes \sigma_{m}$, where $ \rho_{i}$, $\sigma_{j}$ are non-empty $  \otimes $-irreducible partition diagrams for $1\leq i\leq n$, $1\leq j\leq m$,$$\eta_{\rho}  \star  \eta_{\sigma}=\mathop{\sum}(-1)^{n(\pi')}\eta_{\pi'}. $$ The sum is over possible expressions $\pi'=\pi_{1} \divideontimes_{1} \pi_{2} \divideontimes_{2} \cdots \divideontimes_{n+m-1} \pi_{n+m} $ where $ \pi_{1}\pi_{2}\cdots\pi_{n+m}$ is a word shuffle of  $\rho_{1}\rho_{2} \cdots \rho_{n} $ and $ \sigma_{1}\sigma_{2} \cdots \sigma_{m}$, and $ \divideontimes_{i} \in\left\lbrace  \otimes , \bullet \right\rbrace$,  and only when $\pi_{l}=\rho_{i},\pi_{l+1}=\sigma_{j}$ or $\pi_{l}=\sigma_{j},\pi_{l+1}=\rho_{i}$ for some $1\leq i\leq n$, $1\leq j\leq m$, $ \divideontimes_{l}$ can be $ \bullet $.
		$$n(\pi')=\#\left\lbrace l|\pi_{l}=\sigma_{j},\pi_{l+1}=\rho_{i}\text{ for some }1\leq i\leq n,1\leq j\leq m\right\rbrace. $$
		
	\end{Proposition}
	\begin{Proof}
		For any non-empty partition diagrams $\rho=\rho_{1} \otimes \rho_{2} \otimes  \cdots  \otimes \rho_{n}$, $\sigma=\sigma_{1} \otimes \sigma_{2} \otimes  \cdots  \otimes \sigma_{m}$, where $ \rho_{i}$, $\sigma_{j}$ are non-empty $  \otimes $-irreducible partition diagrams for $1\leq i\leq n$, $1\leq j\leq m$,
		\begin{equation*}
			\begin{aligned}
				\eta_{\rho}  \star  \eta_{\sigma}&=(\mathop{\sum}\limits_{\rho\leq \rho'}2^{l(\rho')}M_{\rho'}) \star (\mathop{\sum}\limits_{\sigma\leq \sigma'}2^{l(\sigma')}M_{\sigma'})\\
				&=\mathop{\sum}\limits_{\rho\leq \rho',\sigma\leq \sigma'}2^{l(\rho')+l(\sigma')}M_{\rho'} \star  M_{\sigma'}\\
				&=\mathop{\sum}\limits_{\rho\leq \rho',\sigma\leq \sigma'}2^{l(\rho')+l(\sigma')}\mathop{\sum}\limits_{\pi}M_{\pi}\\
				&=\mathop{\sum}\limits_{\rho\leq \rho',\sigma\leq \sigma'}2^{l(\rho')+l(\sigma')}\mathop{\sum}\limits_{\pi}(2^{-l(\pi)}\mathop{\sum}\limits_{\pi\leq \pi'}(-1)^{l(\pi)-l(\pi')}\eta_{\pi'})\\
				&=\mathop{\sum}\limits_{\rho\leq \rho',\sigma\leq \sigma'}\mathop{\sum}\limits_{\pi}2^{l(\rho')+l(\sigma')-l(\pi)}(\mathop{\sum}\limits_{\pi\leq \pi'}(-1)^{l(\pi)-l(\pi')}\eta_{\pi'})\\
			\end{aligned}	
		\end{equation*}
		Each $\pi'$ is obtained by the following steps:\\
		(1) Change some $ \otimes $s in $\rho$ and $\sigma$ into $ \bullet $s, then we get $\rho'$ and $\sigma'$;\\
		(2) Shuffle the $ \otimes $-irreducible part of $\rho'$ and $\sigma'$,  connect $\rho_{i}$ and $\sigma_{j}$ by $ \otimes $ or $ \bullet $, then we get $\pi$;\\
		(3) Change some $ \otimes $s in $\pi$ into $ \bullet $s.\\
		Then $ \pi'$ can also be written as $\pi'=\pi_{1} \divideontimes_{1} \pi_{2} \divideontimes_{2} \cdots \divideontimes_{n+m-1} \pi_{n+m}$, where $ \pi_{1}\pi_{2}\cdots\pi_{n+m}$ is a word shuffle of  $\rho_{1}\rho_{2} \cdots \rho_{n} $ and $ \sigma_{1}\sigma_{2} \cdots \sigma_{m}$, and $ \divideontimes_{i} \in\left\lbrace  \otimes , \bullet \right\rbrace$. The $ \bullet $s in this expression between $\rho_{i}$ and $\rho_{i+1}$ or between  $\sigma_{j}$ and $\sigma_{j+1}$ are from step (1) or (3); the ones between  $\rho_{i}$ and $\sigma_{j}$ are from step (2) or (3); and those between $\sigma_{j}$ and $\rho_{i}$ only are from step (3). Fix $\pi'$, and denote the numbers of the $ \bullet $s of the above three kinds by $n_{1} $,  $n_{2} $ and $n_{3}$ respectively. Then we calculate the coefficient of $ \eta_{\pi'}$.
		Notice that $l(\rho')+l(\sigma')-l(\pi)$ is the number of $ \bullet $s between  $\rho_{i}$ and $\sigma_{j}$ coming from step (2), $l(\pi)-l(\pi')$ is the number of $ \bullet $s coming from step (3). Then the coefficient of $ \eta_{\pi'}$ is 
		\begin{equation*}
			\begin{aligned}
				&\mathop{\sum}\limits_{0\leq j\leq n_{1} }\mathop{\sum}\limits_{0\leq i\leq n_{2} }\binom{n_{1}}{j}\binom{n_{2}}{i}  2^{i}(-1)^{n_{1}+n_{2}+n_{3}-j-i}\\&=(-1)^{n_{2}+n_{3}}(\mathop{\sum}\limits_{0\leq i\leq n_{2} }\binom{n_{2}}{i}(-2)^{i})\mathop{\sum}\limits_{0\leq j\leq n_{1} }\binom{n_{1}}{j}(-1)^{n_{1}-j}
			\end{aligned}	
		\end{equation*}
		Notice that \begin{equation*}
			\mathop{\sum}\limits_{0\leq j\leq n_{1} }\binom{n_{1}}{j}(-1)^{n_{1}-j}=\begin{cases}
				1 & \text{if } n_{1}=0,\\ 
				0 & \text{otherwise}.
			\end{cases}
		\end{equation*}
		When $n_{1}=0 $, the coefficient of $ \eta_{\pi'}$ is 
		\begin{equation*}
			\begin{aligned}
				(-1)^{n_{2}+n_{3}}(\mathop{\sum}\limits_{0\leq i\leq n_{2} }\binom{n_{2}}{i}(-2)^{i})&=(-1)^{n_{2}+n_{3}}(-2+1)^{n_{2}}\\
				&=(-1)^{n_{2}+n_{3}}(-1)^{n_{2}}\\
				&=(-1)^{n_{3}}
			\end{aligned}	
		\end{equation*}
		$\qedsymbol$
	\end{Proof}
	There is another way to express the product:
	For any non-empty partition diagrams $\rho=\rho_{1} \otimes \rho_{2} \otimes  \cdots  \otimes \rho_{n}$, $\sigma=\sigma_{1} \otimes \sigma_{2} \otimes  \cdots  \otimes \sigma_{m}$, where $ \rho_{i}$, $\sigma_{j}$ are non-empty $  \otimes $-irreducible partition diagrams for $1\leq i\leq n$, $1\leq j\leq m$ (with~\ref{length unique}, every non-empty partition diagram can be uniquely written in this form), $$\eta_{\rho}  \star  \eta_{\sigma}=\mathop{\sum}(-1)^{n_{\pi'}}\eta_{\pi'},$$ where the sum is over all expressions $\pi'=\rho'_{1} \bullet \sigma'_{1} \divideontimes_{1} \rho'_{2} \bullet \sigma'_{2} \divideontimes_{2} \cdots \divideontimes_{k-1} \rho'_{k} \bullet \sigma'_{k} $, where $ k\in \mathbb{N}$, $ \divideontimes_{i} \in\left\lbrace  \otimes , \bullet \right\rbrace$, $n_{\pi'}=\#\left\lbrace i|  \divideontimes_{i} = \bullet \right\rbrace $, $\left( \rho'_{1},\ldots,\rho'_{k}\right)\hat{}=\left( \rho_{1},\ldots,\rho_{n}\right) $ and $\left( \sigma'_{1},\ldots,\sigma'_{k}\right)\hat{}=\left( \sigma_{1},\ldots,\sigma_{m}\right) $, and such that $\rho'_{s} \bullet \sigma'_{s}$ is non-empty for all $1\leq s\leq k$, and $\sigma'_{i}$ and $\rho'_{i+1}$ are non-empty if $ \divideontimes_{i} = \bullet $.

	By Takeuchi's formula, we can get an antipode as follows:  $S(\eta_{\emptyset})=\eta_{\emptyset}$. For any non-empty partition diagram $\pi=\pi_{1} \otimes \pi_{2} \otimes \cdots \otimes \pi_{n}$, where $\pi_{i}$ is non-empty $ \otimes $-irreducible partition diagram for $1\leq i\leq n$,
	$$S (\eta_{\pi})=\mathop{\sum}\limits_{k\geq0}(-1)^{k}(\mathop{\sum}\limits_{r\geq 0}(-1)^{n_{\pi'}}\eta_{\pi'}),$$ where the sum is over all expressions $\pi'=\pi_{1}^{1} \bullet \pi_{2}^{1} \bullet \cdots \bullet \pi_{k}^{1} \divideontimes_{1} \pi_{1}^{2} \bullet \pi_{2}^{2} \bullet \cdots \bullet \pi_{k}^{2} \divideontimes_{2} \cdots \divideontimes_{r-1} \pi_{1}^{r} \bullet \pi_{2}^{r} \bullet \cdots \bullet \pi_{k}^{r}$ where $ k\in \mathbb{N}$, $ \divideontimes_{i} \in\left\lbrace  \otimes , \bullet \right\rbrace$, $n_{\pi'}=\#\left\lbrace i|  \divideontimes_{i} = \bullet \right\rbrace$, $(\pi_{1}^{1},\pi_{1}^{2},\ldots,\pi_{1}^{r},\pi_{2}^{1},\ldots,\pi_{2}^{r},\ldots,\pi_{k}^{1},\ldots,\pi_{k}^{r})\hat{}=(\pi_{1},\pi_{2},\ldots,\pi_{n})$ and $\pi_{1}^{s} \bullet \pi_{2}^{s} \bullet \cdots \bullet \pi_{k}^{s}$ and $\pi_{i}^{1} \otimes \pi_{i}^{2} \otimes \cdots \otimes \pi_{i}^{r}$ are non-empty for $1\leq s\leq r$, $1\leq i\leq k $, and $\max\left\lbrace l |\pi_{l}^{i}\neq\emptyset\right\rbrace >\min\left\lbrace m |\pi_{m}^{i+1}\neq\emptyset\right\rbrace$ if $ \divideontimes_{i} = \bullet $.
	
	\section{\texorpdfstring{The enriched $q$-monomial basis and its dual}{The enriched q-monomial basis and its dual}}\label{The enriched q-monomial basis}
	\subsection{\texorpdfstring{The enriched $q$-monomial basis}{The enriched q-monomial basis}}
	The results in the last section can be generalized from 2 to any invertible $r$ in $\mathbb{K}$. For a fixed invertible $r$ in $\mathbb{K}$, let $q:=r-1$. D. Grinberg and E. Vassilieva~\cite{enriched q-monomial basis} defined the enriched $q$-monomial basis of $\text{QSym}$: $$\eta_{\alpha}^{(q)}=\mathop{\sum}\limits_{\alpha\leq \beta}r^{l(\beta)}M_{\beta}, $$
	for any composition $\alpha$. 
	By analogy with $\left\lbrace \eta_{\alpha}^{(q)}\right\rbrace $, we define the following basis of $\text{ParQSym} $:
	Let $$\eta_{\pi}^{(q)}:=\mathop{\sum}\limits_{\pi\leq \sigma}r^{l(\sigma)}M_{\sigma}, $$then  $\left\lbrace \eta_{\pi}^{(q)}\right\rbrace_{\pi\in A}$ forms a deconcatenation basis of $\text{ParQSym} $ since $r$ is invertible, that is $$\bigtriangleup\eta_{\pi}^{(q)}=\mathop{\sum}\limits_{\rho \otimes \sigma=\pi}\eta_{\rho}^{(q)} \otimes \eta_{\sigma}^{(q)},$$for all $\pi$. In this case, applying (\ref{transform}), we have that $$M_{\pi}:=r^{-l(\pi)}\mathop{\sum}\limits_{\pi\leq \sigma}(-1)^{l(\pi)-l(\sigma)}\eta_{\sigma}^{(q)}. $$
	\begin{Lemma}\cite[Lemma 3.13]{enriched q-monomial basis}
		Let $S$ and $T$ be two finite sets. Then,\begin{equation*}
			\mathop{\sum}\limits_{I\subseteq S}(-1)^{|I\setminus T|}q^{|I\cap T|}=\begin{cases}
				r^{|S|} &\text{if } S\subseteq T,\\
				0  &\text{otherwise}.
			\end{cases}
		\end{equation*} 
	\end{Lemma}
	With the above lemma, we can find the relation between the $L$-basis and the $\eta^{(q)}$-basis of \text{ParQSym}, similar to  \cite[Proposition 3.11]{enriched q-monomial basis}.
	\begin{Proposition} The relation between the $L$-basis and the $\eta^{(q)}$-basis of \text{ParQSym} is as follows:\\	
		(1) For any partition diagram $\pi$, $$\eta_{\pi}^{(q)}=r\mathop{\sum}\limits_{\pi\sim \rho}(-1)^{|S(\rho)\setminus S(\pi)|}q^{|S(\rho)\cap S(\pi)|}L_{\rho}.$$
		(2) For any $\pi\in I_{n}$, $$L_{\pi}=r^{-n}\mathop{\sum}\limits_{\pi\sim \rho}(-1)^{|S(\pi)\setminus S(\rho)|}q^{|[n-1]\setminus(S(\rho)\cup S(\pi))|}\eta_{\rho}^{(q)}.$$
	\end{Proposition}
	\begin{Proof}
		(1) For any partition diagram $\pi$,
		\begin{equation*}
			\begin{aligned}
				&r\mathop{\sum}\limits_{\pi\sim \rho}(-1)^{|S(\rho)\setminus S(\pi)|}q^{|S(\rho)\cap S(\pi)|}L_{\rho}\\
				&=r\mathop{\sum}\limits_{\rho:\pi\sim \rho}(-1)^{|S(\rho)\setminus S(\pi)|}q^{|S(\rho)\cap S(\pi)|}\mathop{\sum}\limits_{\sigma:\sigma\leq \rho}M_{\sigma}\\
				&=r\mathop{\sum}\limits_{\sigma:\pi\sim \sigma}(\mathop{\sum}\limits_{\rho:\sigma\sim\rho, S(\rho)\subseteq S(\sigma)}(-1)^{|S(\rho)\setminus S(\pi)|}q^{|S(\rho)\cap S(\pi)|})M_{\sigma}.
			\end{aligned}
		\end{equation*}
		Notice that $\mathop{\sum}\limits_{\rho:\sigma\sim\rho, S(\rho)\subseteq S(\sigma)}(-1)^{|S(\rho)\setminus S(\pi)|}q^{|S(\rho)\cap S(\pi)|}$ is nonzero only when $S(\sigma)\subseteq S(\pi)$. In this case, $\pi\leq\sigma$, and $$\mathop{\sum}\limits_{\rho:\sigma\sim\rho, S(\rho)\subseteq S(\sigma)}(-1)^{|S(\rho)\setminus S(\pi)|}q^{|S(\rho)\cap S(\pi)|}=r^{|S(\sigma)|}.$$
		Therefore 	\begin{equation*}
			\begin{aligned}
				r\mathop{\sum}\limits_{\pi\sim \rho}(-1)^{|S(\rho)\setminus S(\pi)|}q^{|S(\rho)\cap S(\pi)|}L_{\rho}
				&=r\mathop{\sum}\limits_{\pi\leq \sigma}r^{|S(\sigma)|}M_{\sigma}\\
				&=\mathop{\sum}\limits_{\pi\leq \sigma}r^{l(\sigma)}M_{\sigma}\\
				&=\eta_{\pi}^{(q)}.
			\end{aligned}
		\end{equation*}
		(2) For any $\pi\in I_{n}$,
		\begin{equation*}
			\begin{aligned}
				&r^{-n}\mathop{\sum}\limits_{\pi\sim \rho}(-1)^{|S(\pi)\setminus S(\rho)|}q^{|[n-1]\setminus(S(\rho)\cup S(\pi))|}\eta_{\rho}^{(q)}\\&=	r^{-n}\mathop{\sum}\limits_{\pi\sim \rho}(-1)^{|S(\pi)\setminus S(\rho)|}q^{|[n-1]\setminus(S(\rho)\cup S(\pi))|}\mathop{\sum}\limits_{\sigma:\sigma\geq \rho}r^{l(\sigma)}M_{\sigma}\\
				&=r^{-n}\mathop{\sum}\limits_{\sigma:\pi\sim \sigma}(\mathop{\sum}\limits_{\rho:\sigma\geq \rho}(-1)^{|S(\pi)\setminus S(\rho)|}q^{|[n-1]\setminus(S(\rho)\cup S(\pi))|})r^{l(\sigma)}M_{\sigma}.\\
			\end{aligned}
		\end{equation*}
		Notice that \begin{equation*}
			\begin{aligned}
				&\mathop{\sum}\limits_{\rho:\sigma\geq \rho}(-1)^{|S(\pi)\setminus S(\rho)|}q^{|[n-1]\setminus(S(\rho)\cup S(\pi))|}\\&=\mathop{\sum}\limits_{ ([n-1]\setminus S(\rho))\subseteq ([n-1]\setminus S(\sigma))}(-1)^{|([n-1]\setminus S(\rho))\setminus ([n-1]\setminus S(\pi))|}q^{|([n-1]\setminus S(\rho))\cap ([n-1]\setminus S(\pi))|}
			\end{aligned}
		\end{equation*} is nonzero only when $([n-1]\setminus S(\sigma))\subseteq ([n-1]\setminus S(\pi))$that is,$S(\sigma)\supseteq S(\pi)$. In this case, $\pi\geq\sigma$, and $$\mathop{\sum}\limits_{\rho:\sigma\geq \rho}(-1)^{|S(\pi)\setminus S(\rho)|}q^{|[n-1]\setminus(S(\rho)\cup S(\pi))|}=r^{|[n-1]\setminus S(\sigma)|}.$$
		Therefore 	\begin{equation*}
			\begin{aligned}
				&r^{-n}\mathop{\sum}\limits_{\pi\sim \rho}(-1)^{|S(\pi)\setminus S(\rho)|}q^{|[n-1]\setminus(S(\rho)\cup S(\pi))|}\eta_{\rho}^{(q)}\\
				&=r^{-n}\mathop{\sum}\limits_{\pi\geq \sigma}r^{|[n-1]\setminus S(\sigma)|}r^{l(\sigma)}M_{\sigma}\\
				&=\mathop{\sum}\limits_{\pi\geq \sigma}M_{\sigma}\\
				&=L_{\pi}.
			\end{aligned}
		\end{equation*}
		$\qedsymbol$
	\end{Proof}
	From the definition we can find $\eta_{\emptyset}^{(q)}=M_{\emptyset}$ is the unit. The next proposition is to tell the product of $\eta^{(q)}$-basis for non-empty partition diagrams, which is similar to  \cite[Theorem 5.1]{enriched q-monomial basis}.
	\begin{Proposition}
		For any non-empty partition diagrams $\rho=\rho_{1} \otimes \rho_{2} \otimes  \cdots  \otimes \rho_{n}$, $\sigma=\sigma_{1} \otimes \sigma_{2} \otimes  \cdots  \otimes \sigma_{m}$, where $ \rho_{i}$, $\sigma_{j}$ are non-empty $  \otimes $-irreducible partition diagrams for $1\leq i\leq n$, $1\leq j\leq m$,$$\eta_{\rho}^{(q)}  \star  \eta_{\sigma}^{(q)}=\mathop{\sum}(-1)^{n_{2}(\pi')+n_{3}(\pi')}(-q)^{n_{2}(\pi')}\eta_{\pi'}. $$ The sum is over possible expressions $\pi'=\pi_{1} \divideontimes_{1} \pi_{2} \divideontimes_{2} \cdots \divideontimes_{n+m-1} \pi_{n+m}$, and where $ \divideontimes_{i} \in\left\lbrace  \otimes , \bullet \right\rbrace$, and only when $\pi_{l}=\rho_{i},\pi_{l+1}=\sigma_{j}$ or $\pi_{l}=\sigma_{j},\pi_{l+1}=\rho_{i}$ for some $1\leq i\leq n$, $1\leq j\leq m$, $ \divideontimes_{l}$ can be $ \bullet $.
		$$n_{2}(\pi')=\#\left\lbrace l|\pi_{l}=\rho_{i},\pi_{l+1}=\sigma_{j}\text{ for some } 1\leq i\leq n,1\leq j\leq m\right\rbrace, $$and $$n_{3}(\pi')=\#\left\lbrace l|\pi_{l}=\sigma_{j},\pi_{l+1}=\rho_{i}\text{ for some } 1\leq i\leq n,1\leq j\leq m\right\rbrace. $$
	\end{Proposition}
	\begin{Proof}
		For any partition diagrams $\rho$ and $\sigma$,
		\begin{equation*}
			\begin{aligned}
				\eta^{(q)}_{\rho}  \star  \eta^{(q)}_{\sigma}&=(\mathop{\sum}\limits_{\rho\leq \rho'}r^{l(\rho')}M_{\rho'}) \star (\mathop{\sum}\limits_{\sigma\leq \sigma'}r^{l(\sigma')}M_{\sigma'})\\
				&=\mathop{\sum}\limits_{\rho\leq \rho',\sigma\leq \sigma'}r^{l(\rho')+l(\sigma')}M_{\rho'} \star  M_{\sigma'}\\
				&=\mathop{\sum}\limits_{\rho\leq \rho',\sigma\leq \sigma'}r^{l(\rho')+l(\sigma')}\mathop{\sum}\limits_{\pi}M_{\pi}\\
				&=\mathop{\sum}\limits_{\rho\leq \rho',\sigma\leq \sigma'}r^{l(\rho')+l(\sigma')}\mathop{\sum}\limits_{\pi}(r^{-l(\pi)}\mathop{\sum}\limits_{\pi\leq \pi'}(-1)^{l(\pi)-l(\pi')}\eta^{(q)}_{\pi'})\\
				&=\mathop{\sum}\limits_{\rho\leq \rho',\sigma\leq \sigma'}\mathop{\sum}\limits_{\pi}r^{l(\rho')+l(\sigma')-l(\pi)}(\mathop{\sum}\limits_{\pi\leq \pi'}(-1)^{l(\pi)-l(\pi')}\eta^{(q)}_{\pi'}).\\
			\end{aligned}	
		\end{equation*}
		With $\pi'$ fixed, using the same analysis as the proof of the product of the $\eta$-basis, we denote the numbers of the $ \bullet $s of the three kinds by $n_{1} $,  $n_{2} $, and $n_{3}$ respectively. Then we calculate the coefficient of $ \eta^{(q)}_{\pi'}$.
		
		\begin{equation*}
			\begin{aligned}
				&\mathop{\sum}\limits_{0\leq j\leq n_{1} }\mathop{\sum}\limits_{0\leq i\leq n_{2} }\binom{n_{1}}{j}\binom{n_{2}}{i}  r^{i}(-1)^{n_{1}+n_{2}+n_{3}-j-i}\\&=(-1)^{n_{2}+n_{3}}(\mathop{\sum}\limits_{0\leq i\leq n_{2} }\binom{n_{2}}{i}(-r)^{i})\mathop{\sum}\limits_{0\leq j\leq n_{1} }\binom{n_{1}}{j}(-1)^{n_{1}-j}.
			\end{aligned}	
		\end{equation*}
		Notice that \begin{equation*}
			\mathop{\sum}\limits_{0\leq j\leq n_{1} }\binom{n_{1}}{j}(-1)^{n_{1}-j}=\begin{cases}
				1	& \text{if } n_{1}=0,\\ 
				0	& \text{otherwise}.
			\end{cases}
		\end{equation*}
		When $n_{1}=0 $, the coefficient of $ \eta^{(q)}_{\pi'}$ is 
		\begin{equation*}
			\begin{aligned}
				(-1)^{n_{2}+n_{3}}(\mathop{\sum}\limits_{0\leq i\leq n_{2} }\binom{n_{2}}{i}(-r)^{i})&=(-1)^{n_{2}+n_{3}}(-r+1)^{n_{2}}\\
				&=(-1)^{n_{2}+n_{3}}(-q)^{n_{2}}.\\
			\end{aligned}	
		\end{equation*}
		$\qedsymbol$
	\end{Proof}
	Again, by  Takeuchi's formula, we can get an antipode as follows:  $S(\eta^{(q)}_{\emptyset})=\eta^{(q)}_{\emptyset}$. For any non-empty partition diagram $\pi$, 
	$$S (\eta^{(q)}_{\pi})=\mathop{\sum}\limits_{k\geq0}(-1)^{k}(\mathop{\sum}\limits_{r\geq 0}(-1)^{n_{2}(\pi')+n_{3}(\pi')}(-q)^{n_{2}(\pi')}\eta^{(q)}_{\pi'}),$$ where $\pi'$ are the same as those in $S (\eta_{\pi})$, $n_{2}(\pi')$ is the number of $ \bullet $s between  $\rho_{i}$ and  $\sigma_{j}$, and $n_{3}(\pi')$ is the number of $ \bullet $s between $\sigma_{j}$ and $\rho_{i}$.
	\subsection{The dual basis}
	Darij Grinberg and Ekaterina A. Vassilieva~\cite[Definition 4.2]{enriched q-monomial basis} defined the basis of \text{NSym} dual to the basis $\left\lbrace \eta^{(q)}_{\alpha} \right\rbrace $ of \text{QSym}. \\
	For each composition $\alpha$, let $\eta_{\alpha}^{\ast(q)}:=\mathop{\sum}\limits_{\alpha\geq \beta}(-1)^{l(\beta)-l(\alpha)}r^{-l(\beta)}H_{\beta}\in \text{NSym}$. Similarly, we define the dual basis of  $\left\lbrace \eta^{(q)}_{\pi} \right\rbrace $ as follows.
	\begin{Definition}
		For any partition diagram $\pi$, define an element in \text{ParSym} by
		$$\kappa^{(q)}_{\pi}:=\mathop{\sum}\limits_{\pi\geq \sigma}(-1)^{l(\sigma)-l(\pi)}r^{-l(\sigma)}H_{\sigma}. $$
	\end{Definition}
	With the definition, we can rewrite the $H$-basis as follows:
	\begin{Proposition}
		For any partition diagrams $\pi$, $$H_{\pi}:=r^{l(\pi)}\mathop{\sum}\limits_{\pi\geq \sigma}\kappa^{(q)}_{\sigma}. $$
	\end{Proposition}
	\begin{Proof}
		For any partition diagrams $\pi$,
		\begin{equation*}
			\begin{aligned}
				r^{l(\pi)}\mathop{\sum}\limits_{\pi\geq \sigma}\kappa^{(q)}_{\sigma}&=r^{l(\pi)}\mathop{\sum}\limits_{\pi\geq \sigma}\mathop{\sum}\limits_{\sigma\geq \rho}(-1)^{l(\rho)-l(\sigma)}r^{-l(\rho)}H_{\rho}\\
				&=r^{l(\pi)}\mathop{\sum}\limits_{\rho:\pi\geq \rho}(\mathop{\sum}\limits_{\sigma:\pi\geq \sigma\geq\rho}(-1)^{l(\rho)-l(\sigma)})r^{-l(\rho)}H_{\rho}\\
				&=r^{l(\pi)}\mathop{\sum}\limits_{\rho:\pi\geq \rho}\delta_{\rho,\pi}r^{-l(\rho)}H_{\rho}\\
				&=H_{\pi}.
			\end{aligned}
		\end{equation*}
		$\qedsymbol$
	\end{Proof}
	The above proposition also tells us the 	$\left\lbrace \kappa^{(q)}_{\pi} \right\rbrace $ is the basis of \text{ParSym}. Now we check its duality to the $\eta^{(q)}$-basis of \text{ParSym}.
	\begin{Proposition}
		$\left\lbrace \kappa^{(q)}_{\pi} \right\rbrace $ is the basis of \text{ParSym} dual to the basis $\left\lbrace \eta^{(q)}_{\pi} \right\rbrace $ of \text{ParSym}.
	\end{Proposition}
	\begin{Proof}
		For any partition diagrams $\pi$ and $\rho$, similar to the proof of the duality of $\left\lbrace L_{\pi} \right\rbrace $ and $\left\lbrace R_{\pi} \right\rbrace $, we have
		\begin{equation*}
			\begin{aligned}
				\left\langle \eta^{(q)}_{\rho},\kappa^{(q)}_{\pi} \right\rangle &=\left\langle \mathop{\sum}\limits_{\rho\leq \rho'}r^{l(\rho')}M_{\rho'}, \mathop{\sum}\limits_{\pi\geq \pi'}r^{-l(\pi')}(-1)^{l(\pi')-l(\pi)}H_{\pi'} \right\rangle\\
				&=\mathop{\sum}\limits_{\rho\leq \rho',\pi\geq \pi'}r^{l(\rho')-l(\pi')}(-1)^{l(\pi')-l(\pi)}\delta_{\rho',\pi'}\\
				&=\mathop{\sum}\limits_{\rho\leq \pi'\leq\pi}(-1)^{l(\pi')-l(\pi)}\\		
				&=\left\{\begin{aligned}
					&\mathop{\sum}\limits_{S(\rho)\supseteq S(\pi') \supseteq S(\pi)}(-1)^{|S(\pi')\setminus S(\pi)|}, \text{if } \pi\sim\rho,\\ 
					&0, \text{otherwise}
				\end{aligned}
				\right.
				\\
				&=\begin{cases}
					\mathop{\sum}\limits_{S\subseteq S(\rho)\setminus S(\pi)}(-1)^{\#S} & \text{if } \pi\sim\rho,\\ 
					0 & \text{otherwise}\\
				\end{cases}
				\\
				&=\delta_{\rho,\pi}.\\
			\end{aligned}
		\end{equation*} $\qedsymbol$
	\end{Proof}
	Now we give the product of $\kappa^{(q)}$-basis.
	\begin{Proposition}
		For any partition diagrams $\pi$ and $\rho$, we have that$$\kappa^{(q)}_{\pi}\kappa^{(q)}_{\rho}=\kappa^{(q)}_{\pi \otimes \rho}.$$
	\end{Proposition}
	\begin{Proof}
		For any partition diagrams $\pi$ and $\rho$, similar to the proof of the duality of $\left\lbrace L_{\pi} \right\rbrace $ and $\left\lbrace R_{\pi} \right\rbrace $, we have
		\begin{equation*}
			\begin{aligned}
				\kappa^{(q)}_{\pi}\kappa^{(q)}_{\rho}=&(\mathop{\sum}\limits_{\pi\geq \pi'}r^{-l(\pi')}(-1)^{l(\pi')-l(\pi)}H_{\pi'})( \mathop{\sum}\limits_{\rho\geq \rho'}r^{-l(\rho')}(-1)^{l(\rho')-l(\rho)}H_{\rho'})\\
				=&\mathop{\sum}\limits_{\pi\geq \pi',\rho\geq \rho'}r^{-(l(\pi')+l(\rho'))}(-1)^{(l(\pi')+l(\rho'))-(l(\pi)+l(\rho))}H_{\pi'}H_{\rho'}\\
				=&\mathop{\sum}\limits_{\pi\geq \pi',\rho\geq \rho'}r^{-(l(\pi' \otimes \rho'))}(-1)^{(l(\pi' \otimes \rho'))-(l(\pi \otimes \rho))}H_{\pi' \otimes \rho'}\\
				=&\mathop{\sum}\limits_{\pi \otimes \rho\geq \sigma}r^{-(l(\sigma))}(-1)^{(l(\sigma))-(l(\pi \otimes \rho))}H_{\sigma}\\
				=&\kappa^{(q)}_{\pi \otimes \rho}.
			\end{aligned}
		\end{equation*} $\qedsymbol$
	\end{Proof}
	From the above proposition and that \text{ParSym} is a Hopf algebra, we have that $$\bigtriangleup\kappa^{(q)}_{\pi \otimes \rho}=(\bigtriangleup\kappa^{(q)}_{\pi})(\bigtriangleup\kappa^{(q)}_{\rho}),$$
	$$S(\kappa^{(q)}_{\pi \otimes \rho})=S(\kappa^{(q)}_{\rho})S(\kappa^{(q)}_{\pi}),$$for any partition diagrams $\pi$ and $\rho$. 
	So we can get the coproduct and the antipode of $\left\lbrace \kappa^{(q)}_{\pi} \right\rbrace $ recursively by the length.  
	When $l(\pi)=0$, that is, $\pi=\emptyset$, we find that $\kappa^{(q)}_{\emptyset}=H_{\emptyset}$,  so $$S(\kappa^{(q)}_{\emptyset})=S(H_{\emptyset})=H_{\emptyset}=\kappa^{(q)}_{\emptyset}$$ $$\bigtriangleup\kappa^{(q)}_{\emptyset}=\bigtriangleup H_{\emptyset}=H_{\emptyset} \otimes  H_{\emptyset}=\kappa^{(q)}_{\emptyset} \otimes  \kappa^{(q)}_{\emptyset}.$$
	When $l(\pi)=1$,  that is, $\pi$ is non-empty $ \otimes $-irreducible, 
	\begin{equation*}
		\begin{aligned}
			&\bigtriangleup\kappa^{(q)}_{\pi}\\
			&=\mathop{\sum}\limits_{\pi\geq \sigma}(-1)^{l(\sigma)-1}r^{-l(\sigma)}\bigtriangleup H_{\sigma}\\
			&=\mathop{\sum}\limits_{\pi=\pi_{1} \bullet \cdots \bullet \pi_{k}, \pi_{i}\neq \emptyset}(-1)^{k-1}r^{-k}\bigtriangleup H_{\pi_{1} \otimes \cdots \otimes \pi_{k}}\\
			&=\mathop{\sum}\limits_{\pi=\pi_{1} \bullet \cdots \bullet \pi_{k},\pi_{i}\neq \emptyset}(-1)^{k-1}r^{-k}(\bigtriangleup H_{\pi_{1}})\cdots(\bigtriangleup H_{\pi_{k}})\\
			&=\mathop{\sum}\limits_{\substack{\pi=\pi_{1} \bullet \cdots \bullet \pi_{k}\\\pi_{i}\neq \emptyset\\\pi_{i}=\pi_{i}' \bullet  \pi_{i}''}}(-1)^{k-1}r^{-k}( H_{\pi_{1}'} \otimes  H_{\pi_{1}''})\cdots( H_{\pi_{k}'} \otimes  H_{\pi_{k}''})\\
			&=\mathop{\sum}\limits_{\substack{\pi=\pi_{1} \bullet \cdots \bullet \pi_{k}\\\pi_{i}\neq \emptyset\\\pi_{i}=\pi_{i}' \bullet  \pi_{i}''}}(-1)^{k-1}r^{-k}H_{\pi_{1}' \otimes \cdots \otimes \pi_{k}'} \otimes  H_{\pi_{1}'' \otimes \cdots \otimes \pi_{k}''}\\
			&=\mathop{\sum}\limits_{\substack{\pi=\pi_{1} \bullet \cdots \bullet \pi_{k}\\\pi_{i}\neq \emptyset\\\pi_{i}=\pi_{i}' \bullet  \pi_{i}''\\\sigma'\leq\pi_{1}' \otimes \cdots \otimes \pi_{k}'\\\sigma''\leq\pi_{1}'' \otimes \cdots \otimes \pi_{k}''}}(-1)^{k-1}r^{l(\pi_{1}' \otimes \cdots \otimes \pi_{k}')+l(\pi_{1}'' \otimes \cdots \otimes \pi_{k}'')-k}\kappa^{(q)}_{\sigma'} \otimes  \kappa^{(q)}_{\sigma''}.\\
		\end{aligned}
	\end{equation*}
	Notice that $\pi_{i}'$ and $\pi_{j}''$ may be $\emptyset$, so $ l(\pi_{1}' \otimes \cdots \otimes \pi_{k}') $ might be strictly smaller than $k$. For example, if $\pi$ is also $ \bullet $-irreducible, then $\kappa^{(q)}_{\pi}=\frac{1}{r} H_{\pi}$, so 
	\begin{equation*}
		\begin{aligned}
			\bigtriangleup\kappa^{(q)}_{\pi}&=\frac{1}{r}\bigtriangleup H_{\pi}\\&=\frac{1}{r}(H_{\emptyset} \otimes  H_{\pi}+H_{\pi} \otimes  H_{\emptyset})\\&=\kappa^{(q)}_{\emptyset} \otimes \kappa^{(q)}_{\pi}+\kappa^{(q)}_{\pi} \otimes \kappa^{(q)}_{\emptyset},
		\end{aligned}
	\end{equation*}
	$$S(\kappa^{(q)}_{\pi})= \frac{1}{r} S(H_{\pi})=-\frac{1}{r} H_{\pi}=-\kappa^{(q)}_{\pi}.$$
	In general, the coproduct of the $\kappa^{(q)} $-basis may be complicated, but there are some special cases. Recall the injective combinatorial Hopf morphism $\Phi:\text{NSym}\rightarrow\text{ParSym}$, which maps $H_{n}$ and $H_{\alpha}$ to $H_{\pi(n)}$ and  $H_{\pi(\alpha)}$ ($H_{\pi(n)}$ and  $H_{\pi(\alpha)}$ are defined in Section~\ref{Combinatorial Hopf algebra}) respectively. It is easy to see $\Phi(\eta_{\alpha}^{\ast(q)})=\kappa^{(q)}_{\pi(\alpha)} $ from the definitions. Darij Grinberg and Ekaterina A. Vassilieva~\cite[Theorem 4.15]{enriched q-monomial basis} found that $$\bigtriangleup\eta_{n}^{\ast(q)}=\mathop{\sum}\limits_{\substack{|\beta|+|\gamma|=n\\|l(\beta)-l(\gamma)|\leq 1}}(-q)^{\max\left\lbrace l(\beta),l(\gamma)\right\rbrace -1}(q-1)^{\delta_{l(\beta),l(\gamma)}}\eta_{\beta}^{\ast(q)} \otimes \eta_{\gamma}^{\ast(q)},$$where $\bigtriangleup\eta_{n}^{\ast(q)}:= \bigtriangleup\eta_{(n)}^{\ast(q)}$. Since $\Phi$ is a coalgebra morphism, we have that 
	\begin{equation*}
		\begin{aligned}
			\bigtriangleup\kappa^{(q)}_{\pi(n)}=&\bigtriangleup\circ\Phi(\eta_{n}^{\ast(q)})\\=&(\Phi \otimes \Phi)\bigtriangleup(\eta_{n}^{\ast(q)})\\
			=&(\Phi \otimes \Phi)(\mathop{\sum}\limits_{\substack{|\beta|+|\gamma|=n\\|l(\beta)-l(\gamma)|\leq 1}}(-q)^{\max\left\lbrace l(\beta),l(\gamma)\right\rbrace -1}(q-1)^{\delta_{l(\beta),l(\gamma)}}\eta_{\beta}^{\ast(q)} \otimes \eta_{\gamma}^{\ast(q)})\\
			=&\mathop{\sum}\limits_{\substack{|\beta|+|\gamma|=n\\|l(\beta)-l(\gamma)|\leq 1}}(-q)^{\max\left\lbrace l(\beta),l(\gamma)\right\rbrace -1}(q-1)^{\delta_{l(\beta),l(\gamma)}}\kappa^{(q)}_{\pi(\beta)} \otimes \kappa^{(q)}_{\pi(\gamma)}.
		\end{aligned}
	\end{equation*}
	The following  proposition tells the relation between $\kappa^{(q)}$-basis and $R$-basis of \text{ParSym}.
	\begin{Proposition} The relation between $\kappa^{(q)}$-basis and $R$-basis of \text{ParSym} is as follows:\\
		(1) For any partition diagram $\pi$, $$R_{\pi}=r\mathop{\sum}\limits_{\pi\sim \rho}(-1)^{|S(\rho)\setminus S(\pi)|+l(\rho)-l(\pi)}q^{|S(\rho)\cap S(\pi)|}\kappa^{(q)}_{\rho}.$$
		(2) For any $\pi\in I_{n}$,
		$$\kappa^{(q)}_{\pi}=r^{-n}\mathop{\sum}\limits_{\pi\sim \rho}(-1)^{|S(\pi)\setminus S(\rho)|+l(\rho)-l(\pi)}q^{|[n-1]\setminus(S(\rho)\cup S(\pi))|}R_{\rho}.$$
	\end{Proposition}
	\begin{Proof}
		(1) For any partition diagram $\pi$,
		\begin{equation*}
			\begin{aligned}
				&r\mathop{\sum}\limits_{\pi\sim \rho}(-1)^{|S(\rho)\setminus S(\pi)|+l(\rho)-l(\pi)}q^{|S(\rho)\cap S(\pi)|}\kappa^{(q)}_{\rho}\\&=r\mathop{\sum}\limits_{\rho:\pi\sim \rho}(-1)^{|S(\rho)\setminus S(\pi)|+l(\rho)-l(\pi)}q^{|S(\rho)\cap S(\pi)|}\mathop{\sum}\limits_{\sigma:\sigma\leq \rho}(-1)^{l(\sigma)-l(\rho)}r^{-l(\sigma)}H_{\sigma}\\
				&=r\mathop{\sum}\limits_{\sigma:\pi\sim \sigma}(\mathop{\sum}\limits_{\rho:\sigma\sim\rho, S(\rho)\subseteq S(\sigma)}(-1)^{|S(\rho)\setminus S(\pi)|}q^{|S(\rho)\cap S(\pi)|})(-1)^{l(\sigma)-l(\pi)}r^{-l(\sigma)}H_{\sigma}.
			\end{aligned}
		\end{equation*}
		Notice that $\mathop{\sum}\limits_{\rho:\sigma\sim\rho, S(\rho)\subseteq S(\sigma)}(-1)^{|S(\rho)\setminus S(\pi)|}q^{|S(\rho)\cap S(\pi)|}$ is nonzero only when $S(\sigma)\subseteq S(\pi)$. In this case, $\pi\leq\sigma$, and $$\mathop{\sum}\limits_{\rho:\sigma\sim\rho, S(\rho)\subseteq S(\sigma)}(-1)^{|S(\rho)\setminus S(\pi)|}q^{|S(\rho)\cap S(\pi)|}=r^{|S(\sigma)|}.$$
		Therefore 	\begin{equation*}
			\begin{aligned}
				r\mathop{\sum}\limits_{\pi\sim \rho}(-1)^{|S(\rho)\setminus S(\pi)|+l(\rho)-l(\pi)}\kappa^{(q)}_{\rho}
				&=r\mathop{\sum}\limits_{\pi\leq \sigma}r^{|S(\sigma)|}(-1)^{l(\sigma)-l(\pi)}r^{-l(\sigma)}H_{\sigma}\\
				&=\mathop{\sum}\limits_{\pi\leq \sigma}(-1)^{l(\sigma)-l(\pi)}H_{\sigma}\\
				&=R_{\pi}.
			\end{aligned}
		\end{equation*}
		
		(2)	For any $\pi\in I_{n}$,
		
		\begin{equation*}
			\begin{aligned}
				&r^{-n}\mathop{\sum}\limits_{\pi\sim \rho}(-1)^{|S(\pi)\setminus S(\rho)|+l(\rho)-l(\pi)}q^{|[n-1]\setminus(S(\rho)\cup S(\pi))|}R_{\rho}^{(q)}\\&=	r^{-n}\mathop{\sum}\limits_{\pi\sim \rho}(-1)^{|S(\pi)\setminus S(\rho)|+l(\rho)-l(\pi)}q^{|[n-1]\setminus(S(\rho)\cup S(\pi))|}\mathop{\sum}\limits_{\sigma:\sigma\geq \rho}(-1)^{l(\sigma)-l(\rho)}H_{\sigma}\\
				&=r^{-n}\mathop{\sum}\limits_{\sigma:\pi\sim \sigma}(\mathop{\sum}\limits_{\rho:\sigma\geq \rho}(-1)^{|S(\pi)\setminus S(\rho)|}q^{|[n-1]\setminus(S(\rho)\cup S(\pi))|})(-1)^{l(\sigma)-l(\pi)}H_{\sigma}.\\
			\end{aligned}
		\end{equation*}
		Notice that \begin{equation*}
			\begin{aligned}
				&\mathop{\sum}\limits_{\rho:\sigma\geq \rho}(-1)^{|S(\pi)\setminus S(\rho)|}q^{|[n-1]\setminus(S(\rho)\cup S(\pi))|}\\&=\mathop{\sum}\limits_{ ([n-1]\setminus S(\rho))\subseteq ([n-1]\setminus S(\sigma))}(-1)^{|([n-1]\setminus S(\rho))\cap ([n-1]\setminus S(\pi))|}q^{|([n-1]\setminus S(\rho))\cap ([n-1]\setminus S(\pi))|}
			\end{aligned}
		\end{equation*} is nonzero only when $([n-1]\setminus S(\sigma))\subseteq ([n-1]\setminus S(\pi))$that is,$S(\sigma)\supseteq S(\pi)$. In this case, $\pi\geq\sigma$, and $$\mathop{\sum}\limits_{\rho:\sigma\geq \rho}(-1)^{|S(\pi)\setminus S(\rho)|}q^{|[n-1]\setminus(S(\rho)\cup S(\pi))|}=r^{|[n-1]\setminus S(\sigma)|}.$$
		Therefore 	\begin{equation*}
			\begin{aligned}
				&r^{-n}\mathop{\sum}\limits_{\pi\sim \rho}(-1)^{|S(\pi)\setminus S(\rho)|+l(\rho)-l(\sigma)}q^{|[n-1]\setminus(S(\rho)\cup S(\pi))|}R_{\rho}^{(q)}\\
				&=r^{-n}\mathop{\sum}\limits_{\pi\geq \sigma}r^{|[n-1]\setminus S(\sigma)|}(-1)^{l(\sigma)-l(\pi)}H_{\sigma}\\
				&=\mathop{\sum}\limits_{\pi\geq \sigma}(-1)^{l(\sigma)-l(\pi)}r^{-l(\sigma)}H_{\sigma}\\
				&=\kappa^{(q)}_{\pi}.
			\end{aligned}
		\end{equation*}
		$\qedsymbol$
	\end{Proof}
	\section{Competing interests statement}
	The authors declare that they have no known competing financial interests or personal relationships that could have appeared to influence the work reported in this paper.
	

\begin{thebibliography}{99}  
		\bibitem{ref10} Claudia Malvenuto and Christophe Reutenauer,  Duality between quasi-symmetric functions and the Solomon descent algebra,  {\it J. Algebra} {\bf 177} (1995) 967--982, \url{https://doi.org/10.1006/jabr.1995.1336}.
		
		
		
		\bibitem{planar trees}D. Arcis and S. M\'{a}rquez, Hopf algebras on planar trees and permutations, {\it  J.
			Algebra Appl} {\bf 21} (2022) 2250224,	\url{https://doi.org/10.1142/S0219498822502243}.
		\bibitem{David}David E Radford, {\it Hopf Algebras} (Series on Knots and Everything, {\bf 49}. Singapore: World Scientific Publishing Co. Pte. Ltd, 2011) \url{https://doi.org/10.1142/8055}
		\bibitem{enriched monomial basis}
		D. Grinberg and E. Vassilieva, Weighted posets and the enriched monomial basis of QSym, {\it S\'{e}min. Loth. de Comb.} {\bf 85B} (2021), \url{https://doi.org/10.48550/arxiv.2202.04720}.
		
		\bibitem{enriched q-monomial basis}D. Grinberg and E. Vassilieva, The enriched $q$-monomial basis of the quasisymmetric functions,  {\it arXiv e-prints} (2023), \url{https://doi.org/10.48550/arXiv.2309.01118}.
		
		
		
		\bibitem{ref2}D. Grinberg and V. Reiner,  {\it Hopf Algebras in Combinatorics} (2014) {\it arXiv e-prints}, \url{https://doi.org/10.48550/arxiv.1409.8356}.
		
		\bibitem{refDominique Manchon}Dominique Manchon,  {\it Hopf algebras, from basics to applications to renormalization}, (2003) \url{https://doi.org/10.48550/arxiv.math/0408405}.
		\bibitem{diagrams}G. H. E. Duchamp, J.-G. Luque, J.-C. Novelli, C. Tollu, F. Toumazet. (2011). Hopf
		algebras of diagrams.  {\it Int. J. Algebra Comput}. {\bf 21} 889--911, \url{https://doi.org/10.1142/S0218196711006418}.
		
		\bibitem{ref35} Israel M. Gelfand, Daniel Krob, Alain Lascoux, Bernard Leclerc, Vladimir S. Retakh and
		Jean-Yves Thibon,  Noncommutative symmetric functions.  {\it Adv. Math}. {\bf 112} (1995) 218--348, \url{https://doi.org/10.1006/aima.1995.1032}.
		
		\bibitem{ref36} Ira M. Gessel, Multipartite $P$-partitions and inner products of skew Schur functions, {\it Contemp. Math}. {\bf 34} (1984) 289--317, \url{https://doi.org/10.1090/conm/034/777705}.
		\bibitem{planar binary trees}J.-L. Loday and M. O. Ronco, Hopf algebra of the planar binary trees,   {\it Adv.
			Math.} {\bf 139} (1998) 293--309, \url{https://doi.org/10.1006/aima.1998.1759}.
		\bibitem{ref1}	John M. Campbell, A Combinatorial Hopf Algebra on partition diagrams, {\it Commun. Algebra} {\bf 53} (2025) 2723--2744, \url{https://doi.org/10.1080/00927872.2025.2450033}.
		\bibitem{Ordered forests}L. Foissy and J. Unterberger, Ordered forests, permutations, and iterated integrals,  {\it Int. Math. Res. Not.} {\bf 2013} (2013) 846--885, \url{https://doi.org/10.1093/imrn/rnr273}.
		\bibitem{ref4} M. Aguiar, N. Bergeron and F. Sottile,  Combinatorial Hopf algebras and generalized Dehn-Sommerville relations,  {\it Compos. Math}. {\bf 142} (2006) 1--30, \url{https://doi.org/10.1112/S0010437X0500165X}.
		\bibitem{uniform block permutations}M. Aguiar and R. C. Orellana, The Hopf algebra of uniform block permutations,
		{\it J. Algebr. Comb}. {\bf 28} (2008) 115--138, \url{https://doi.org/10.1007/s10801-008-0120-9}.
		
		\bibitem{ref5} Mitsuhiro Takeuchi, Free Hopf algebras generated by coalgebras, {\it J. Math. Soc. Japan} {\bf 23} (1971) 561--582, \url{https://doi.org/10.2969/jmsj/02340561}.
		\bibitem{Sweedler}Moss E. Sweedler, {\it Hopf Algebras} (Mathematics Lecture Note Series.  W. A. Benjamin, Inc., New York, 1969).
		\bibitem{ref9}N. Bergeron and M. Zabrocki, The Hopf algebras of symmetric functions and quasi-symmetric functions in non-commutative variables are free and cofree, {\it J. Algebra Appl.} {\bf 08} (2009) 581--600, \url{https://doi.org/10.1142/S0219498809003485}.
		
		\bibitem{ref8}Ricky Ini Liu and Michael Tang, Shuffle bases and quasisymmetric power sums, {\it arXiv e-prints} (2023), \url{https://arxiv.org/abs/2310.09371}.
		\bibitem{families of trees} R. Grossman and R. G. Larson,  Hopf-algebraic structure of families of trees, {\it J.
			Algebra} {\bf 126} (1989) 184--210, \url{https://doi.org/10.1016/0021-8693(89)90328-1}.
		\bibitem{Ricky Ini Liu and Michael Weselcouch}Ricky Ini Liu and Michael Weselcouch,  P-Partitions and Quasisymmetric Power Sums,  {\it International Mathematics Research Notices} {\bf 2021} (2020) 12707--12747, \url{https://doi.org/10.1093/imrn/rnz375}.
		
		
		\bibitem{Montgomery} Susan Montgomery, {\it Hopf Algebras and Their Actions on Rings} (CBMS Regional Conference
		Series in Mathematics, {\bf 82} Washington, DC: Conference Board of the Mathematical Sciences, 1993)
		
		
		
		
		\bibitem{partition algebra}	T. Halverson, Characters of the partition algebras, {\it J. Algebra} {\bf 238} (2001) 502--533, \url{https://doi.org/10.1006/jabr.2000.8613}.
		
		
		
		
		
		\bibitem{combinatorial structures}W. R. Schmitt, Hopf algebras of combinatorial structures, {\it Can. J. Math.} {\bf 45} (1993) 412--428, \url{https://doi.org/10.4153/CJM-1993-021-5}.
		
		
		
		
		
		
		
		
		
		
		
		
		
		
		
		
		
		
		
		
		
		
		
		
		
		
		
		
		
		
		
		
		\bibitem{subgraph}X. Wang, S.-J. Xu and X. Gao, A Hopf algebra on subgraphs of a graph,  {\it J.
			Algebra Appl.} {\bf 19} (2020) 2050164, \url{https://doi.org/10.1142/S0219498820501649}.
		
		\bibitem{Rota-Baxter}X. Zhang, A. Xu and L. Guo, Hopf algebra structure on free Rota-Baxter algebras by angularly decorated rooted trees,  {\it J. Algebr. Comb.} {\bf 55} (2022) 1331-
		1349, \url{https://doi.org/10.1007/s10801-021-01098-8}.
		\bibitem{solomon2024}
		Y. Solomon, Algebra Structure on Set Partitions, Doctoral dissertation, Graduate Program in Mathematics and Statistics, York University, Toronto, Ontario, 2024.
		
	\end{thebibliography}
\end{document}